\input epsf.tex

\def\C{{\bf C}}
\def\R{{\bf R}}
\def\Z{{\bf Z}}

\def\Q{{\bf Q}}
\def\N{{\bf N}}
\def\ra{{\longrightarrow}}
\def\raa{{\hbox to 25pt{\rightarrowfill} }}

\def\laa{{\hbox to 25pt{\leftarrowfill} }}
\def\g{{\bf g}}
\def\MC{{\cal MC}}

\def\a{{A_{\infty}}}
\def\ga{{\gamma}}
\def\h{{\hbar}}
\def\m{{\bf m}}
\def\p{{\partial}}
\def\L{{L_{\infty}}}
\def\F{{\cal F}}
\def\G{{\Gamma}}
\def\a{{\alpha}}
\def\H{{\cal H}}
\def\OC{{\overline C}}
\def\U{{\cal U}}
\def\V{{\cal V}}
\def\D{{\cal D}}
\def\T{{\cal T}}
\def\W{{\cal W}}
\def\c{{\cal C}}
\def\from{{\leftarrow}}
\def\k{{\bf k}}
\def\s{{\sigma}}
\def\qed{{\quad {\it Q.E.D.}}}
\def\hra{{\hookrightarrow}}
\def\b{{\bullet}}

\centerline{\bf DEFORMATION QUANTIZATION OF POISSON MANIFOLDS, I}\par
\medskip
\medskip
\centerline { {\bf Maxim Kontsevich}}
\vskip 2truecm

\noindent {\bf 0. Introduction}
\vskip 0.5truecm

In this paper it is  proven
 that any finite-dimensional Poisson manifold can be 
canonically quantized (in the sense of deformation quantization).
 Informally, it means that the set of equivalence classes of
  associative algebras close to algebras of functions on manifolds is in 
   one-to-one correspondence with the set of equivalence classes
   of Poisson manifolds modulo diffeomorphisms.
    This is a corollary of a more general statement, which I
   proposed around 1993-1994 (``Formality conjecture'') (see [Ko2], [V]).

  For a long
   time the Formality conjecture resisted all approaches. The solution 
    presented here uses in a essential way ideas of string theory. 
    Our formulas  can be viewed as a perturbation series
     for a topological two-dimensional
      quantum field theory coupled with gravity.
     
  \vskip 0.5truecm
\par\noindent {\bf 0.1. Content of the paper}
\vskip 0.5truecm      
     
      Section 1: an elementary introduction to the deformation
       quantization, and precise formulation of the
        main statement 
       concerning Poisson manifolds.
       
       Section 2: an explicit formula
         for the deformation quantization written in coordinates.
         
          Section 3: an introduction to the deformation theory
           in general, in terms of differential graded Lie algebras.
             The material of this section is basically standard.

            Section 4:  a geometric reformulation
             of the theory introduced in the previous section, in terms
              of odd vector fields on formal supermanifolds.
              In particular, we introduce convenient notions of
               an $\L$-morphism and of a quasi-isomorphism, 
               which gives us a tool
                to identify deformation theories related with two 
                differential graded Lie algebras.
                 Also in this section we state our main result, which
                  is an existence of a quasi-isomorphism between
                   the Hochschild complex of the algebra of polynomials, and
                   the graded Lie algebra of polyvector fields.
                  
                  Section 5: tools
                   for  the explicit construction of the quasi-isomorphism
                    mentioned above.
                   We define compactified
                   configuration spaces related with the Lobachevsky
                   plane, a class of admissible graphs, differential 
                   polynomials on polyvector fields
                related with graphs, and integrals
                     over configuration spaces. Technically
                      the same constructions were used 
                       in   generalizations of the perturbative
                        Chern-Simons theory several years ago (see [Ko1]).
                 Compactifications of the configuration spaces
                  are close relatives of Fulton-MacPherson compactifications
                   of configuration spaces in algebraic geometry (see [FM]).

                     Section 6: it is proven that the machinery introduced
                      in the previous section gives a quasi-isomorphism
                       and establishes the Formality conjecture for affine
                        spaces.
                        The proof is essentially an application of the Stokes
                         formula, and a general result of vanishing
                         of certain integral associated with 
                          a collection of rational 
                        functions on a complex
                        algebraic  variety.
                        
                        Section 7: results of section 6 are extended
                         to the case
                         of general manifolds. In order to do this we 
                      recall  ideas of formal geometry
                       of I.~Gelfand and D.~Kazhdan, and the language
                       of  superconnections. In order to pass
                        from the affine space
                         to  general manifolds we have to 
                        find a non-linear
                        cocycle of the Lie algebra of
                        formal vector fields. It turns
                        out that such a cocycle can be almost directly 
                        constructed from our explicit formulas.
                         In the course of the proof we calculate several
                          integrals and check their vanishing. Also,
                           we introduce a 
                           general notion of direct image
                            for certain bundles of supermanifolds.
                            
                            Section 8: we describe an additional structure
                             present in the deformation theory
                             of associative algebras, the cup-product
                              on the tangent bundle to the super moduli
                               space.
                              The isomorphism constructed in 
                              sections 6, 7 is compatible with this 
                              structure. As a corollary, 
                            we finally   justify the orbit method
                                in the representation theory. 
                                One of new results
                                is the validity of Duflo-Kirillov formulas
                                 for Lie algebras in general rigid tensor 
                                 categories,
                                  in particular for Lie superalgebras.
                                  Another application is an equality
                        between two  cup-products in the context 
                                   of algebraic geometry.

         \vskip 0.5truecm
\par\noindent {\bf 0.2. Plans}
\vskip 0.5truecm 
           
           I am going to write a complement to the present paper.
            It will contain 
  
   1) the comparison  with various known constructions of star-products, 
            the most
            notorious one are by De Wilde-Lecomte and by Fedosov
            for the case of  symplectic manifolds (see [DL], [F]),
             and by Etingof-Kazhdan for Poisson-Lie groups (see [EK]),

    2) a reformulation of the Formality conjecture as an existence
               of a natural construction of a triangulated category
             starting   from an odd symplectic supermanifold,
                
      3) discussion of the arithmetic nature of coefficients in
                 our formulas, and of the possibility to extend
                  main results
                 for algebraic varieties over arbitrary 
                 field of characteristic zero,

     4) an application to the Mirror
                 Symmetry, which was the original motivation for
                 the Formality conjecture (see [Ko4]),

       5) a Lagrangian
                  for a quantum field theory (from [AKSZ])
                  which seems to
                   give our formulas after the perturbation expansion. 
                   
                   Also, I am going to touch other topics (a version
         of formality
                    for cyclic homology, quantization of quadratic 
                    brackets, etc).

         \vskip 0.5truecm
\par\noindent {\bf 0.3. Acknowledgements}
\vskip 0.5truecm

    I am grateful to Y.~Soibelman for many remarks.

\vskip 0.5truecm    
 \noindent {\bf 1. Deformation quantization}    
\vskip 0.5truecm

 \noindent {\bf 1.1. Star-products}    
\vskip 0.5truecm

 Let $A=\G(X,{\cal O}_X)$ be the algebra over $\R$ 
 of smooth functions on a finite-dimensional
  $C^{\infty}$-manifold $X$.
  The star-product on $A$ (see [BFFLS]) is an 
  {\it associative}
   $\R[[\hbar]]$-linear product on 
    $A[[\hbar]]$ given by the following 
  formula for $f,g\in A\subset A[[\h]]$:
  $$(f,g)\mapsto f\star g=fg+\hbar B_1(f,g)+\hbar^2 B_2(f,g)+\dots\in 
  A[[\hbar]]\,,$$
  where $\hbar$ is the formal variable, and $B_i$ are bidifferential 
  operators (i.e. bilinear maps $A\times A\ra A$
   which are differential operators with respect to each argument
    of globally bounded order). 
   The product of arbitrary elements of $A[[\hbar]]$ is defined
  by the condition
   of linearity over  $\R[[\hbar]]$ and $\hbar$-adic continuity:
   $$\left(\sum_{n\ge 0} f_n \,\h^n\right)\star
     \left(\sum_{n\ge 0} g_n\, \h^n\right):=\sum_{k,l\ge 0}
      f_k g_l\,\h^{k+l}+\sum_{k,l\ge 0,\,m\ge 1} 
      B_m(f_k,g_l)\,\h^{k+l+m}\,\,\,.$$

   There is a natural gauge group acting on star-products.
    This group consists of automorphisms of $A[[\hbar]]$
     considered as an $\R[[\hbar]]$-module
     (i.e. linear transformations $A\ra A$ parametrized by $\hbar$), 
     of the following form:
   $$f\mapsto f+\hbar D_1(f)+\hbar^2 D_2(f)+\dots
  ,\,\,\,\,{\rm for }\,\,\,\,f\in A\subset A[[\hbar]] \,,$$
  $$\sum_{n\ge 0} f_n\,\h^n\mapsto \sum_{n\ge 0} f_n\,\h^n
  + \sum_{n\ge 0,m\ge 1} D_m(f_n)\,\h^{n+m},\,\,\,\,\,
  {\rm for\,\,\,general\,\,\,element}\,\,\,f(\h)=\sum_{n\ge 0}
  f_n\,\h^n\in A[[\h]]\,,$$
   where $D_i:A\ra A$ are differential operators. 
   If $D(\h)=1+\sum_{m\ge 1} D_m\,\h^m$ is such an automorphism, it acts
   on the set of star-products as
   $$\star\mapsto\star',\,\,\,f(\h)\star'g(\h):=D(\h)
    \bigl(D(\h)^{-1}(f(\h))\star D(\h)^{-1}(g(\h)\bigr),\,\,\,
    f(\h),g(\h)\in A[[\h]]\,\,\,.$$
    We are interested in star-products up to gauge equivalence.
    
 \vskip 0.5truecm    
 \noindent {\bf 1.2. First approximation: Poisson structures}    
\vskip 0.5truecm
  
It follows    from the associativity of $\star$ 
that the bilinear map $B_1:A\times A\ra A$ satisfies the equation
   $$f B_1(g,h)-B_1(fg,h)+B_1(f,gh)-B_1(f,g)h=0,$$
   i.e. the linear map ${\widetilde B}_1:A\otimes A\ra A$  
   associated 
   with $B_1$ as
   ${\widetilde B}_1(f\otimes g):=B_1(f,g)$, 
   is a $2$-cocycle in the cohomological Hochschild 
   complex of algebra $A$ (the definition 
   of this complex is given in 3.4.2).
   
   Let us decompose $B_1$ into the sum of the 
   symmetric part and of the anti-symmetric part:
   $$B_1=B_1^+ +B_1^-,\,\,\,\,
   B_1^+(f,g)=B_1^+(g,f),\,\,\,B_1^-(f,g)=-B_1^-(g,f)\,\,\,.  $$
   Gauge transformations 
   $$B_1\mapsto B_1',\,\,\,\,\,
   B_1'(f,g)=B_1(f,g)-f D_1(g)+D_1(fg)-D_1(f)g$$
   where $D_1$ is an arbitrary differential operator,
    affect only the symmetric part of $B_1$, i.e. 
     $B_1^-=(B_1')^-$.
    One can show that the symmetric part $B_1^+$ can be killed by
     a gauge transformation 
     (and it is a coboundary in the Hochschild complex).                           
                     
  Also one can show that the skew-symmetric part $B^-_1$ 
   is  a derivation with respect to 
   both functions $f$
    and $g$. Thus, $B_1^-$ comes from a bi-vector field $\a$ on $X$:
     $$B_1^-(f,g)=\langle 
     \a, df \otimes dg\rangle
     ,\,\,\,\, 
     \a\in \G(X, \wedge^2 T_X)\subset \G(X, T_X\otimes T_X)\,\,\,.$$ 
   Analogous fact in algebraic geometry is  
   that the second Hochschild cohomology group of the 
   algebra
    of functions on a smooth affine
     algebraic variety (in characteristic zero) 
  is naturally isomorphic to the space of bi-vector fields (see 4.6.1.1).

    The second term $O(\hbar^2)$ in the associativity equation
     $f\star(g\star h)=(f\star g)\star h$ implies that $\a$ 
     gives a Poisson
      structure on $X$,
      $$\forall f,g,h\,\,\,\,\,\{f,\{g,h\}\}+\{g,\{h,f\}\}+\{h,\{f,g\}\}=0
      ,$$
     $${\rm where}\,\,\,\,
     \{f,g\}:=
     {f\star g-g\star f\over \hbar}_{|\hbar=0}=2\,B_1^-(f,g)=2\langle 
     \a, df \otimes dg\rangle
      \,\,\,.$$
   In  other words, $[\a,\a]=0\in \G(X,\wedge ^3 T_X)$, where the 
   bracket 
   is  
   the Schouten-Nijenhuis bracket on polyvector fields (see 4.6.1
    for the definition of this bracket).

   Thus,  gauge equivalence classes of
   star-products modulo $O(\h^2)$
    are classified by  Poisson structures on $X$. 
    {\it A priori} it is not clear whether there exists a star-product
     with the first term equal to a given Poisson structure, and whether
      there exists a preferred choice of an equivalence class
       of star-products.
      We show in this paper that there is a canonical construction
       of an equivalence class of star-products for any
        Poisson manifold.
        
         \vskip 0.5truecm    
 \noindent {\bf 1.3. Description of quantizations}    
\vskip 0.5truecm

  \proclaim Theorem.  
  The  set of gauge equivalence classes of star
     products on a smooth manifold $X$
      can be naturally identified with
      the set of equivalence classes of
       Poisson structures depending formally on $\hbar$:
       $$\a=\a(\hbar)=\a_1\hbar+\a_2 \hbar^2+\dots\in
       \G(X,\wedge^2 T_X)[[\hbar]],\,\,\,[\a,\a]=0\in 
       \G(X,\wedge^3 T_X)[[\hbar]]$$
        modulo the action of the group of
         formal paths in the diffeomorphism group of $X$, 
         starting at the identity diffeomorphism. 
         \par

    Any given Poisson structure $\a_{(0)}$ gives a path
     $\a(\h):=\a_{(0)} \cdot\h $ and by the Theorem from above, a 
      canonical gauge equivalence class of star products.

 \vskip 0.5truecm    
 \noindent {\bf 1.4. Examples}    
\vskip 0.5truecm

 \noindent {\bf 1.4.1. Moyal product}    
\vskip 0.5truecm

The simplest example of a deformation quantization
 is the Moyal product for the Poisson structure on $\R^d$ with 
 constant coefficients:
 $$\alpha=\sum_{i,j} \alpha^{ij} \p_i \wedge
 \p_j,\,\,\,\a^{ij}=-\a^{ji}\in \R$$ 
 where $\p_i=\p/\p x^i$ is the partial derivative in the direction of 
 coordinate $x^i,\,\,\,i=1,\dots, d$.
 The formula for the Moyal product is
 $$f\star g=fg+\hbar\sum_{i,j} \alpha^{ij} \,\p_i (f) \,\p_j (g)+
 {\hbar^2\over 2} \sum_{i,j,k,l}\alpha^{ij} \alpha^{kl}\, \p_i\p_k(f) \,
 \p_j \p_l (g)+\dots=$$
 $$=\sum_{n=0}^{\infty} {\hbar^n \over n !}\,\,
  \sum_{i_1,\dots,i_n;\,j_1,\dots j_n}\,\, \prod_{k=1}^n
  \a^{i_k j_k}\left(\prod_{k=1}^n \p_{i_k}\right)
   (f)\times \left(\prod_{k=1}^n \p_{j_k}\right)
   (g)\,\,\,.$$
Here and later symbol $\times$ denotes the usual product.
   
  \vskip 0.5truecm    
 \noindent {\bf 1.4.2. Deformation quantization up to the second order}    
\vskip 0.5truecm

   Let $\alpha=\sum_{i,j} \alpha^{ij} \p_i \wedge\p_j$ be a Poisson
      bracket with variable coefficients in an open domain of $\R^d$ 
       (i.e. $\alpha^{ij}$ is not a constant, but a function of
        coordinates), then the following formula
         gives an associative product modulo $O(\h^3)$:
       $$f\star g = fg+\hbar \sum_{i,j}\alpha^{ij} \,\p_i (f) \,\p_j (g)+
 {\hbar^2\over 2} \sum_{i,j,k,l}\alpha^{ij} \alpha^{kl} \,\p_i \p_k(f) \,
 \p_j \p_l (g)+$$
 $$+{\hbar^2\over 3} \left(\sum_{i,j,k,l} \alpha^{ij} \,\p_j(\alpha^{kl})
 \left(
 \p_i \p_k(f)\, \p_l(g)-\p_k(f)\,\p_i \p_l(g)\right)\right)+O(\hbar^3)$$
 The associativity up to the second order
 means that for any 3 functions $f,g,h$ one has
  $$(f\star g)\star h= f\star(g\star h)+ O(\hbar^3)\,\,\,.$$

 \vskip 0.5truecm    
 \noindent {\bf 1.5. Remarks}    
\vskip 0.5truecm

   In general, one  should consider
   bidifferential operators $B_i$  with complex
   coefficients, as we expect to associate by quantization 
  self-adjoint operators in a Hilbert space
 to   real-valued classical observables. In this paper
     we deal with purely formal algebraic properties of the deformation
      quantization and work mainly over the field $\R$ of real numbers.

        Also, it is not clear whether the ``deformation quantization''
        is natural for quantum mechanics. This question
         we will discuss in the next paper.
         A topological open string theory
          seems to be more relevant. 
  
 \vskip 0.5truecm    
 \noindent {\bf 2. Explicit universal formula}    
\vskip 0.5truecm 

Here we propose a formula for the star-product for arbitrary Poisson 
 structure $\a$ in an open domain of the standard
  coordinate space $\R^d$. Terms of our formula modulo $O(\hbar^3)$ 
  are  the same as in the previous section, plus a gauge-trivial
   term of order $O(\hbar^2)$, symmetric in $f$ and $g$.
   Terms of the formula are certain universal polydifferential operators
    applied to coefficients of the bi-vector field $\a$ and 
    to functions $f,g$. All indices corresponding to coordinates
     in the formula 
     appear once as lower indices and once as upper indices, 
     i.e. the formula is invariant under
      affine transformations of $\R^d$.
      
      In order to describe terms proportional to $\hbar^n$
      for any integer $n\ge 0$, we introduce a special class $G_n$ 
      of oriented labeled graphs.
      
      \proclaim Definition. An (oriented) graph $\G$ is a pair $(V_{\G},E_\G)$ 
      of two finite 
      sets such that $E_\G$ is a subset of $V_\G\times V_\G$.
      \par
      
       Elements of $V_\G$ are  vertices of $\G$, elements
        of $E_\G$ are edges of $\G$. If $e=(v_1,v_2)\in E_\G\subseteq V_\G
        \times V_\G$ is an edge then we say that $e$ starts
          at $v_1$ and ends at $v_2$.

         In the usual definition of graphs one
          admits infinite graphs, and also graphs
          with multiple edges. Here we will not meet such structures
           and  use a simplified terminology.
         
       We say that a labeled graph $\G$ belongs to $G_n$ if
       
       1) $\G$ has $n+2$ vertices and $2n$ edges,
       
       2) the set vertices 
         $V_\G$ is $\{1,\dots,n\}\sqcup \{L,R\}$, where
         $L,R$ 
         are just two symbols (capital roman letters, mean Left and
          Right), 
         
       3) edges of $\G$ are labeled by symbols
         $e_1^1,e_1^2,e_2^1,e_2^2,\dots,e_n^1,e_n^2\,,$

        4) for every $k\in \{1,\dots,n\}$ edges labeled by
       $e_k^1$ and $e_k^2$ start
        at the vertex $k$,
       
       5) for any $v\in V_\G$ the ordered pair
         $(v,v)$ is not  an edge of $\G$.

       The 
        set $G_n$ is finite, it has 
       $\bigl(n(n+1)\bigr)^n$ elements for $n\ge 1$ and $1$ element
 for $n=0$.
       
       To each labeled graph $\G\in G_n$ 
         we associate a bidifferential
        operator 
        $$B_{\G,\a}:A\times A\ra A,\,\,\,\,\,
        \,\,\,\,A=C^{\infty}(\V),\,\,\,\V{\rm
        \,\,\,is\,\,\,an\,\,\,open\,\,\,domain\,\,\,in\,\,\,}\R^d$$ 
         which depends  on 
         bi-vector field $\a\in \G(\V,\wedge^2 T_{\V})$, not necessarily a 
         Poisson one. 
        We show one 
        example, from which the general rule should be clear.
         Here $n=3$ and the list of edges is 
         $$
         \bigl(e_1^1,e_1^2,e_2^1,e_2^2,e_3^1,e_3^2\bigr)=
         \bigl((1,L),(1,R),(2,R),(2,3),(3,L),(3,R)\bigr)\,\,\,.$$
         
          \vskip 1cm 
         \centerline{\epsfbox{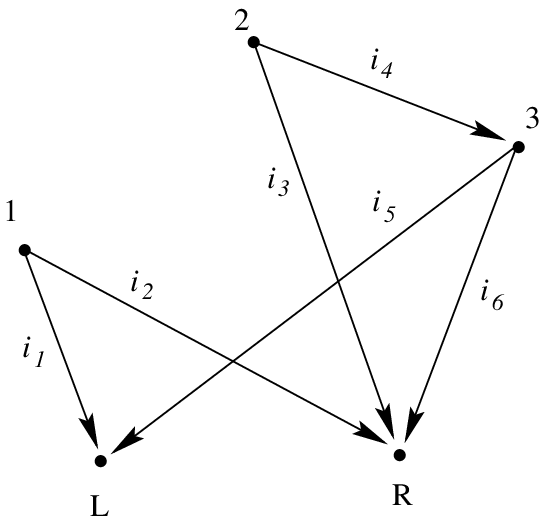}}
       \vskip 11mm

         In the picture of $\G$ we put independent indices
          $\,\,1\le i_1,\dots,i_6\le d\,\,$ on edges,
           instead of labels $e_*^*$. 
         The operator $B_{\G,\a}$ corresponding to this graph is
         $$(f,g)\mapsto \sum_{i_1,\dots, i_6} 
       \a^{i_1 i_2} \a^{i_3 i_4} \p_{i_4}(\a^{i_5 i_6})
       \p_{i_1} \p_{i_5} (f) \p_{i_2}\p_{i_3}\p_{i_6}(g)\,\,\,.$$
       
       The general formula for the operator $B_{\G,\a}$ is
       $$B_{\G,\a}(f,g):=\sum_{
I:E_\G\ra \{1,\dots,d\}}\left[
 \prod_{k=1}^n \left(\prod_{e\in E_\G,\,\,e=(*,k)}\p_{I(e)}\right)
 \a^{I(e_k^1)I(e_k^2)}\right]\times$$
 
  $$\times \left(\prod_{ e\in E_\G,\,\,e=(*,L)} \p_{I(e)}\right)f\times
  \left(\prod_{e\in E_\G,\,\,e=(*,R)} \p_{I(e)}\right)g\,\,\,.$$

       In the next step we associate a weight 
       $W_{\G}\in \R$ with each graph $\G\in G_n$.
        In order to define it we  need an elementary construction 
        from hyperbolic geometry. 
        
        Let $p,q,\,\,p\ne q$ be two points
         on the standard upper half-plane $ \H=\{z\in \C|\,\,Im(z)>0\}$
       endowed with the Lobachevsky metric. We denote by
         $\phi^{h}(p,q)\in \R/2\pi \Z$ the angle at $p$ formed by 
         two  
         lines, $l(p,q)$ and $l(p, \infty)$ passing through $p$ 
         and $q$, and
         through $p$ and the point $\infty$ on the absolute. The direction
          of the measurement of the angle is counterclockwise from
           $l(p,\infty)$ to $l(p,q)$. In the notation $\phi^h (p,q)$
            letter $h$ is
            for {\it harmonic}.
            
            \vskip 1cm
            \centerline{
\epsfbox{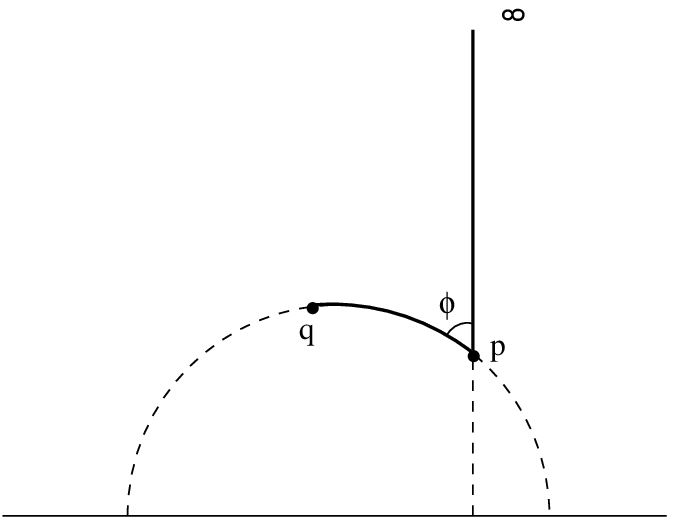}}
\vskip 11mm

           An easy planimetry shows that
        one can express angle $\phi^{h}(p,q)$ in terms of 
            complex numbers:
             $$\phi^{h}(p,q)=Arg((q-p)/(q-{\overline p}))={1\over 2 i}
              Log\left({(q-p)({\overline q}-p)\over
              (q-{\overline p})({\overline q} -
              {\overline p})}\right)\,\,\,.$$
              
          Function $\phi^{h}(p,q)$ can be defined by continuity also in 
          the case
           $p,q\in\H\sqcup \R,\,\,p\ne q$.
           
           Denote by $\H_n$ the  space of configurations
            of $n$ numbered 
           pairwise distinct
           points on $\H$:
           $$\H_n=\{(p_1,\dots,p_n)|p_k\in \H,\,\,\, 
           p_k\ne p_l\,\,\,\,{\rm for}
           \,\,\,k\ne l\}\,\,\,.$$
           $\H_n\subset \C^{n}$ is a non-compact smooth
           $2n$-dimensional manifold.
            We introduce  orientation on $\H_n$ using the natural complex
             structure on it.
            
            If $\G\in G_n$ is a graph as above, and $(p_1,\dots,p_n)\in 
            \H_n$
             is a configuration of points, then we draw a copy of $\G$ 
             on the plane 
          $\R^2\simeq \C$  by assigning point
           $p_k\in \H$ to the vertex  
          $k,\,\,1\le k\le n$, point $0\in \R\subset \C$ to the vertex 
          $L$, 
          and point $1\in \R\subset \C$ to the vertex $R$. Each edge 
          should 
          be drawn
           as a line interval in hyperbolic geometry.
            Every  edge $e$ 
            of the graph $\G$ defines an ordered pair $(p,q)$ 
             of points on $\H\sqcup \R$, thus an angle 
             $\phi^{h}_e:=\phi^{h}(p,q)$.
             If points $p_i$ move around, we get a function
              $\phi^{h}_e$ on $\H_n$ with values in $\R/2 \pi \Z$.
              
              We define the weight  of $\G$ as 
              $$w_\G:={1\over  n! 
              (2\pi)^{2n}}\int\limits_{\H_n}\bigwedge_{i=1}^n 
                  (d\phi^{h}_{e_k^1}\wedge d\phi^{h}_{e_k^2})
                  \,\,\,. $$

     \proclaim Lemma. The integral in the definition  of $w_\G$ is 
     absolutely convergent.
     \par

    This lemma is a particular case of  a more general statement proven in  
     section 6.
    
    \proclaim Theorem. Let $\a$ be a Poisson bi-vector field
     in a domain of $\R^d$. The formula 
    $$f\star g :=\sum_{n=0}^{\infty} {\h^n}
    \sum_{\G\in G_n} w_\G B_{\G,\a}(f,g)$$ defines an associative product.
     If we change coordinates, we obtain a gauge equivalent star-product.
     \par
     
     The proof of this theorem is elementary, it uses only the Stokes
      formula. Again,  this theorem is a corollary
       of a more general statement 
      proven in section 6.

\vskip 0.8truecm
\par\noindent {\bf 3. Deformation theory  via differential graded 
Lie algebras}

\vskip 0.5truecm
\par\noindent {\bf 3.1. Tensor categories $Super$ and $Graded$}
\vskip 0.5truecm

 Here we make a comment about the terminology. This comment
  looks a bit pedantic, but it could help in the struggle with signs
   in formulas.

   The main idea of algebraic geometry is to replace spaces 
    by  commutative associative rings (at least locally). One can further
    generalize this considering commutative associative algebras
   in general 
    tensor categories (see [De]).  In this way
   one can imitate many 
  constructions from algebra and differential geometry.
  
    The fundamental example is supermathematics, i.e. mathematics
  in the tensor category $Super^\k$
  of super vector spaces over a field $\k$ of characteristic zero
   (see Chapter 3 in [M]).
   The category $Super^\k$ is the category
   of $\Z/2\Z$-graded vector spaces over $\k$
   (representations of the group $\Z/2\Z$)
    endowed with the
    standard tensor product, with the standard associativity functor,
     and with a modified
    commutativity functor
    (the Koszul rule of signs). We denote by $\Pi$ the usual functor $Super^\k
    \ra
     Super^\k$
     changing the parity. It is given on objects by the formula
     $\Pi V=V \otimes \k^{0|1}$.
      In the sequel we will consider the standard tensor category $Vect^\k$ of 
       vector spaces over $\k$ as the subcategory of $Super^\k$
        consisting of pure even spaces.

       The basic tensor category which appears everywhere
        in topology and homological algebra is a full subcategory 
         of the tensor category of $\Z$-graded super vector spaces.
         Objects of this category are infinite sums
          ${\cal E}=\oplus_{n\in \Z} 
          {\cal E}^{(n)}$ 
       such  that ${\cal E}^{(n)}$ is pure even for even $n$, 
         and pure odd for odd $n$.
         We will slightly abuse the language  calling this category
         also the category of graded vector spaces,
          and denote it simply by $Graded^\k$. We denote by ${\cal E}^n$ 
         the usual $\k$-vector 
          space underlying the graded component 
          ${\cal E}^{(n)}$. The super vector space obtained  
          if we forget about
          $\Z$-grading on ${\cal E}\in Objects(Graded^\k)$ is
          $\bigoplus_{n\in \Z} \Pi^n({\cal E}^n)$.

          Analogously, we will speak about graded manifolds. They are 
          defined as supermanifolds endowed with $\Z$-grading on the sheaf
           of functions obeying the same conditions on the parity as above.

           The shift functor $[1]:Graded^\k\ra Graded^\k$ 
           (acting from the right)
           is defined
           as the tensor product with graded space $\k[1]$ where
           $\k[1]^{-1}\simeq \k$, $\k[1]^{\ne -1}=0$. Its powers
            are denoted by $[n],\,\,n\in \Z$. Thus, for graded space $\cal E$
            we have
            $${\cal E}=\bigoplus_{n\in \Z} {\cal E}^n[-n]\,\,\,.$$

   Almost all results in the present paper formulated for graded manifolds,
    graded Lie algebras etc.,  
    hold also for supermanifolds, super Lie algebras etc.

\vskip 0.5truecm
\par\noindent {\bf 3.2. Maurer-Cartan equation in differential graded Lie 
algebras}
\vskip 0.5truecm  

This part is essentially standard (see [GM], [HS1], [SS], ...).

Let $\g$ be a differential graded Lie algebra over  field $\k$
of characteristic zero. Below we recall the list of structures and axioms:
$$\g=\bigoplus_{k\in \Z} \g^k[-k],\,\,\,\,\,\,\,
[\,\,,\,\,]:
\g^k\otimes \g^l\ra \g^{k+l},\,\,\,\,\,\,d:\g^k\ra 
\g^{k+1},$$
$$d(d(\ga))=0,\,\,\,d[\ga_1,\ga_2]=[d\ga_1,\ga_2]+(-1)^{\overline{\ga_1}}
[\ga_1,d\ga_2],\,\,\,[\ga_2,\ga_1]=-(-1)^{\overline{\ga_1}
\cdot\overline{\ga_2}}
[\ga_1,\ga_2],$$
$$[\ga_1,[\ga_2,\ga_3]]+(-1)^{\overline{\ga_3}\cdot(\overline{\ga_1}+
\overline{\ga_2})}[\ga_3,[\ga_1,\ga_2]]+
(-1)^{\overline{\ga_1}\cdot(\overline{\ga_2}+
\overline{\ga_3})}[\ga_2,[\ga_3,\ga_1]]=0\,\,\,.$$

 In formulas above
 symbols $\overline{\ga_i}\in\Z$ mean the degrees of  homogeneous elements 
 $\ga_i$, i.e. 
 $\ga_i\in \g^{\overline{\ga_i}}$.
 
  In other words, $\g$ is a Lie algebra in the  tensor 
  category 
  of complexes
   of vector spaces over $\k$. If we forget about the differential 
   and the grading
    on $\g$, we obtain a Lie superalgebra.

 We associate with $\g$ 
a functor $Def_{\g}$ 
on the category of finite-dimensional commutative 
associative
 algebras over $\k$,
  with values in the category of sets. 
 First of all, let us assume that $\g$ is  a nilpotent 
  Lie superalgebra.
  We define set $\MC(\g)$  (the set of solutions of the Maurer-Cartan equation
   modulo the gauge equivalence)
   by the formula
  $$\MC(\g):=\left\{\ga\in \g^1| \,\,d\ga+{1\over  2}[\ga,\ga]=0\right\}\big/
  \Gamma^0$$
  where $\Gamma^0$ is the nilpotent group associated with the nilpotent 
  Lie algebra $\g^0$. The group $\G$ acts by affine transformations 
  of the vector space $\g^1$. The action of $\Gamma^0$ is defined by
  the exponentiation of the infinitesimal action of its Lie algebra:
  $$\alpha\in \g^0\mapsto \left(\dot\ga=d\alpha+[\alpha,\ga]\right)\,\,\,.$$

 Now we are ready to introduce  functor $Def_\g$. Technically, it is convenient
  to define this functor on the category of
  finite-dimensional nilpotent commutative associative algebras 
  {\it without}  unit.
   Let $\m$ be such an algebra, 
  $\m^{dim(\m)+1}=0$. 
  The functor is given (on objects) by the formula
  $$Def_\g(\m)=\MC(\g\otimes \m)\,\,\,.$$
  
  In the conventional approach $\m$ is the maximal ideal in a 
  finite-dimensional
   Artin algebra {\it with} unit $$\m':=\m\oplus 
   \k\cdot {\rm 1}\,\,\,.$$
  In general, one can think about  commutative associative
     algebras without unit as about objects dual to  spaces with base points.
      Algebra corresponding to a space with  base point is the algebra of
       functions vanishing at the base point.
        
        One can  extend the definition of the
         deformation functor to 
         algebras with linear topology
         which are 
       projective limits of nilpotent finite-dimensional algebras. For 
       example, in the deformation quantization we use the following
        algebra  
        over $\R$:  
        $$\m:=\hbar\R[[\hbar]]=\lim_{\from}\left(\h\R[\h]/\h^k\R[\h]\right)
        \rm{\,\,\,\,as\,\,\,\,\,}k\to\infty\,\,\,.$$

  \vskip 0.5truecm
\par\noindent {\bf 3.3. Remark}
\vskip 0.5truecm  

  Several authors, following a suggestion of P.~Deligne,   
        stressed that the set $Def_\g(\m)$
         should be considered as the set of equivalence classes
          of a natural
        groupoid. Almost always 
         in deformation theory, 
          differential graded Lie algebras are supported
           in non-negative degrees, $\g^{<0}=0$. 
           Our principal example in the present
            paper,
           the shifted
           Hochschild complex (see the next subsection), 
           has a non-trivial
            component in degree $-1$, when
            it is considered as a graded Lie algebra. The set $Def_\g(\m)$
             in such a case has a natural structure of the set of equivalence
              classes of a $2$-groupoid.
           In general, if one considers differential graded
            Lie algebras with components in negative degrees, 
            one meets immediately polycategories and nilpotent homotopy
             types. Still, it is only a half of the story
              because one can not
              say 
             anything about $\g^{\ge 3}$ using this language.
              Maybe, the better way is to extend the definition
               of the deformation functor to the category of differential
                graded nilpotent
                commutative associative algebras, see the last remark in
                 4.5.2.

  \vskip 0.5truecm
\par\noindent {\bf 3.4. Examples}
\vskip 0.5truecm

 There are many standard  examples of differential graded Lie algebras and 
 related moduli problems.

\vskip 0.5truecm
\par\noindent {\bf 3.4.1. Tangent complex}
\vskip 0.5truecm

  Let $X$ be a complex manifold. Define $\g$ over $\C$
   as
  $$\g=\bigoplus_{k\in \Z} \g^k[-k];\,\,\,\,
  \g^k=\Gamma(X,\Omega^{0,k}_X\otimes T^{1,0}_X)\,\,\,\,{\rm for }\,\,\,
  k\ge 0,\,\,\,\,\,\g^{<0}=0$$
  with the differential equal to $\overline{\partial}$, and 
  the Lie bracket coming from the cup-product on 
  $\overline{\partial}$-forms
   and the usual Lie bracket on holomorphic vector fields.
   
  The deformation functor related with $\g$ is 
   the usual deformation functor for complex structures on $X$.
   The set
    $Def_\g(\m)$ can be naturally identified with the set of 
    equivalence classes
    of analytic spaces $\widetilde{X}$ endowed with a flat
     map $p:\widetilde{X}\ra Spec(\m')$, and an identification $i:
     \widetilde{X}\times_{Spec(\m')} Spec( \C)\simeq X $ of the 
     special
     fiber of $p$  with $X$.
  
  \vskip 0.5truecm
\par\noindent {\bf 3.4.2. Hochschild complex}
\vskip 0.5truecm
 
  Let $A$ be an associative algebra over  field $\k$ of
  characteristic zero. The graded space of Hochschild cochains
   of $A$ with coefficients in $A$ considered as a bimodule
    over itself is
    $$C^{\bullet}(A,A):=\bigoplus_{k\ge 0}
     C^k(A,A)[-k],\,\,\,C^k(A,A):=Hom_{Vect^\k}(A^{\otimes k},A)\,\,\,.$$
     
     We define graded vector space
  $\g$ over $\k$ by formula $\g:=C^{\bullet}(A,A)[1]$. Thus, we have 
  $$\g=\bigoplus_{k\in \Z} \g^k[-k];\,\,\,\,\g^k:= Hom(A^{\otimes 
  (k+1)},A)\,\,\,{\rm for }\,\,\,k\ge -1,
  \,\,\,\,\,
  \g^{<(-1)}=0\,\,\,.$$

   The differential in $\g$ is shifted by $1$
    the usual differential
  in the Hochschild complex, and the Lie bracket is the Gerstenhaber 
  bracket. The explicit formulas for the differential and for the bracket
   are:
   $$(d\Phi)(a_0\otimes\dots \otimes a_{k+1})=a_0\cdot\Phi(a_1\otimes\dots
   \otimes a_{k+1})-\sum_{i=0}^{k} (-1)^i \Phi(a_0\otimes\dots
   \otimes (a_i\cdot a_{i+1})\otimes\dots \otimes a_{k+1})
   +$$
   $$+(-1)^{k}\Phi(a_0\otimes\dots\otimes a_{k})\cdot a_{k+1}
     ,\,\,\,\,\,\,\,\,\Phi\in \g^k,$$
   and
   
   $$[\Phi_1,\Phi_2]=\Phi_1\circ\Phi_2-(-1)^{k_1 k_2} \Phi_2\circ\Phi_1,
   \,\,\,\,\,\,\,\,\,\Phi_i\in \g^{k_i},$$
   where the (non-associative) product $\circ$ is defined as
   $$(\Phi_1\circ\Phi_2)(a_0\otimes\dots\otimes a_{k_1+k_2})=$$
   $$=
   \sum_{i=0}^{k_1} (-1)^{i k_2}\Phi_1(a_0\otimes\dots\otimes
    a_{i-1}\otimes(\Phi_2(a_i\otimes\dots\otimes a_{i+k_2}))\otimes
    a_{i+k_2+1}\otimes\dots\otimes a_{k_1+k_2})\,\,\,.$$
  
  We would like to give here also an abstract
   definition of the differential and of the bracket on $\g$.
  Let $F$ denote the free coassociative graded 
  coalgebra with counit cogenerated by the graded vector space $A[1]$:
  $$F=\bigoplus_{n\ge 1}
   \otimes^n(A[1])\,\,\,.$$ 
  
   Graded Lie algebra
    $\g$ is the Lie algebra of coderivations of $F$ in the tensor 
   category $Graded^\k$. The associative product 
    on $A$ gives an element $m_A\in \g^1,\,\,\,m_A:A\otimes A\ra A$ 
    satisfying the equation
    $[m_A,m_A]=0$. 
    The differential $d$ in $\g$ is defined as $ad(m_A)$.
    
    Again, the deformation functor related to $\g$ is equivalent to
    the usual deformation functor for algebraic structures. 
    Associative products on $A$ correspond to solutions of the Maurer-
    Cartan
     equation in $\g$. The set $Def_\g(\m)$ is naturally identified with
      the set of equivalence classes of pairs $({\widetilde A},i)$
       where ${\widetilde A}$ is an associative algebra 
       over $\m'=\m\oplus 
   \k\cdot {\rm 1}$ such that
        ${\widetilde A}$ is
         free as an $\m'$-module, and $i$ an isomorphism
        of $\k$-algebras 
        ${\widetilde A}\otimes_{\m'}\k\simeq
        A$.

     The cohomology of the Hochschild complex are 
      $$HH^k(A,A)=Ext^k_{A-mod-A}(A,A),$$
       the $Ext$-groups
       in the abelian category of bimodules over $A$. The 
      Hochschild complex {\it without} shift by $1$ also has a meaning
       in deformation theory, it is responsible for deformations
        of $A$ as a bimodule.

   \vskip 0.5truecm
\par\noindent {\bf 4. Homotopy Lie algebras and quasi-isomorphisms}
\vskip 0.5truecm
   
   In this section we introduce a  language convenient for the homotopy theory
    of 
    differential graded Lie algebras and for the deformation theory.
     The
      ground field $\k$ for linear algebra in our discussion
      is an arbitrary
       field of characteristic zero,
      unless specified.

    \vskip 0.5truecm
\par\noindent {\bf 4.1. Formal manifolds}
\vskip 0.5truecm

Let $V$ be a vector space. We denote by $C(V)$ the cofree cocommutative
coassociative  coalgebra without counit cogenerated by $V$:
$$C(V)=\bigoplus_{n\ge 1}\bigl( \otimes^n V\bigr)^{\Sigma_n}
\subset \bigoplus_{n\ge 1}\bigl( \otimes^n V\bigr)\,\,\,.$$

Intuitively, we think about $C(V)$ as about an object corresponding to
  a formal  manifold, possibly infinite-dimensional, with base point:
  $$\left(V_{formal}, \,\,{\rm base\,\,\,point}\,\right):=
  \,\bigl(\,\,{\rm Formal\,\,\, neighborhood \,\,\,of\,\,\, zero\,\,\, 
  in\,\,\,} V,\,0\,\bigr)\,\,\,.$$
  The reason for this is that if $V$ is finite-dimensional then $C(V)^*$
   (the 
  dual space to $C(V)$)  is the algebra of formal 
  power series on $V$ vanishing at the origin. 
  
  \proclaim Definition. A formal pointed manifold $M$ is an
   object corresponding
   to a coalgebra $\c$ which is
   isomorphic to $C(V)$ for some vector space
    $V$.
   \par
   
   The specific isomorphism between $\c$
    and $C(V)$ is {\it not} considered as a part of data.
    Nevertheless, the vector space $V$ can be reconstructed from $M$
    as the space of primitive elements in coalgebra 
    $\c$.
      Speaking geometrically, $V$ is  the tangent space to $M$
     at the base point.
      A choice of an isomorphism between $\c$
    and $C(V)$ can be considered as a choice of an affine structure on $M$.
  
 If $V_1$ and $V_2$ are two vector spaces then a map $f$ between corresponding
  formal pointed manifolds is defined as a homomorphism of coalgebras
   (the pushforward on distributions supported at zero)
    $$f_*:C(V_1)\ra C(V_2)\,\,\,.$$
    By the universal property of cofree coalgebras any such homomorphism
     is uniquely specified by a linear map
      $$C(V_1)\ra V_2$$ 
      which is the composition of $f_*$ with the canonical projection
       $C(V_2)\ra V_2$. Homogeneous components 
       of this map,
       $$f^{(n)}:
       \bigl(\otimes^n(V_1)\bigr)^{\Sigma_n}\ra V_2,\,\,\,\,n\ge 1$$
         can be considered
       as Taylor coefficients of $f$. More precisely, Taylor coefficients
        are defined as maps 
        $$\partial ^n f: Sym^n(V_1)\ra V_2,
        \,\,\,\,\partial ^n f(v_1\cdot\dots\cdot v_n):=
        {\p^n\over \p t_1\dots \p t_n}
    _{|t_1=\dots=t_n=0} \left(f(t_1 v_1+\dots+t_n v_n)\right)\,\,\,.$$
    Linear map $f^{(n)}$ coincides with $\partial ^n f$ after the
     identification of the {\it subspace}
      $\bigl( \otimes^n V_1\bigr)^{\Sigma_n}\subset \otimes^n V_1$ 
      with the {\it quotient} space 
      $$Sym^n(V_1):= \otimes ^n 
  V_1/\{\,standard\,\,\,relations\,\}\,\,\,.$$

       As in the usual calculus, there is the inverse mapping theorem:
     non-linear map $f$ is invertible iff its first Taylor coefficient 
     $f^{(1)}:V_1\ra V_2$ is invertible.
     
     Analogous definitions and statements can be made in other tensor 
     categories, including  $Super^{\k}$ and $Graded^{\k}$. 
     
      The reader can ask why we speak about  base points
     for formal  manifolds, as such manifolds
      have only one geometric point. The reason is that later we
      will consider formal graded
       manifolds
       depending on formal parameters. In such a situation
         the choice of the base point
        is an essential part of the structure.

   \vskip 0.5truecm
\par\noindent {\bf 4.2. Pre-$\L$-morphisms}
\vskip 0.5truecm
   
Let $\g_1$ and $\g_2$ be two graded vector spaces.

\proclaim Definition. A pre-$\L$-morphism $\F$ from $\g_1$ to $\g_2$ is 
a map of formal pointed graded manifolds
 $$\F:\bigl((\g_1[1])_{formal},0\bigr)\ra \bigl((\g_2[1])_{formal},0\bigr)
 \,\,\,.$$\par

Map $\F$ is defined by its
 Taylor coefficients which are
  linear maps $\p^n \F$ of graded vector spaces:
$$\p^1 \F:\g_1\ra \g_2$$
$$\p^2 \F:\wedge^2(\g_1)\ra \g_2[-1]$$
$$\p^3 \F:\wedge^3(\g_1)\ra \g_2[-2]$$
 $$\dots$$

Here we use the natural isomorphism $Sym^n(\g_1[1])\simeq 
\bigl(\wedge ^n(\g_1)\bigr)[n]$.
     In plain terms, we have a collection of linear maps between
   ordinary vector spaces 
   $$\F_{(k_1,\dots,k_n)}: \g_1^{k_1}\otimes\dots\otimes \g_1^{k_n}\ra 
   \g_2^{k_1+\dots+k_n+(1-n)}$$
   with the symmetry property
   $$\F_{(k_1,\dots,k_n)}(\ga_1\otimes\dots\otimes \ga_n)=
   -(-1)^{k_i k_{i+1}}
   \F_{(k_1,\dots,k_{i+1},k_i,\dots,k_n)}(\ga_1\otimes\dots \otimes\ga_{i+1}
   \otimes\ga_i\otimes\dots\otimes \ga_n)\,\,\,.$$
   One can write (slightly abusing notations) 
   $$\p^n \F(\ga_1\wedge\dots\wedge\ga_n)=
   \F_{(k_1,\dots,k_n)}(\ga_1\otimes\dots\otimes \ga_n)$$
   for $\ga_i\in \g_1^{k_i},\,\,i=1,\dots,n$.
   
    In the sequel we will denote $\p^n \F$ simply by $\F_n$.

 \vskip 0.5truecm
\par\noindent {\bf 4.3. $\L$-algebras and $\L$-morphisms}
\vskip 0.5truecm  

Suppose that we have an odd vector field $Q$ of degree $+1$ on formal
 graded manifold $(\g[1]_{formal},0)$ such that
  the Taylor series for coefficients of 
    $Q$ has
    terms of degree $1$ and $2$ 
      only.
      The first Taylor coefficient $Q_1$ gives a linear  map $\g\ra \g$ of 
      degree $+1$ (or, better, a map
       $\g\ra \g[1]$). The second coefficient $Q^2:\wedge ^2\g\ra \g$ 
       gives a skew-symmetric
       bilinear operation of degree $0$ on $\g$. 

   It is easy to see that if $[Q,Q]=0$ then
    $\g$ is a differential graded Lie algebra, with differential $Q_1$ and
     the bracket $Q_2$, and vice versa.

      In paper [AKSZ] supermanifolds endowed with  an odd
       vector field $Q$ such that $[Q,Q]=0$, are called $Q$-manifolds.
        By analogy, we can speak about 
       formal graded pointed 
    $Q$-manifolds.

   \proclaim Definition. An $\L$-algebra is a  pair $(\g,Q)$ where
    $\g$ is a graded vector space  
   and 
   $Q$ is a differential 
    of degree $+1$ on the graded coalgebra $C(\g)$. \par
    
    Other names for $\L$-algebras are ``(strong) homotopy Lie algebras''
     and ``Sugawara algebras'' (see e.g. [HS2]).
    
    Usually we will denote  $\L$-algebra $(\g,Q)$
    simply by $\g$.

    The structure of an $\L$-algebra on a graded vector space $\g$
    is given by the infinite sequence
    of Taylor coefficients $Q_i$ of the odd vector field $Q$ 
    (coderivation of $C(\g)$):
     $$Q_1: \g\ra \g[1]$$
     $$Q_2: \wedge^2(\g)\ra \g$$
     $$Q_3:\wedge^3(\g)\ra \g[-1]$$
      $$\dots$$
      
      The condition $Q^2=0$ can be translated into an infinite 
      sequence of quadratic constraints on polylinear 
      maps $Q_i$. First of these constraints means  that $Q_1$ 
      is the differential of the graded space $\g$. Thus, $(\g, Q_1)$
      is a complex of vector spaces over $\k$.
       The 
       second constraint means that $Q_2$ is a skew-symmetric 
       bilinear operation on $\g$, for which 
        $Q_1$ satisfies the Leibniz rule. The third constraint means 
        that $Q_2$ satisfies the Jacobi identity
         up to homotopy given by $Q_3$, etc. 
          As we have seen, a differential graded Lie algebra is 
          the same as an $\L$-algebra with 
           $Q_3=Q_4=\dots=0$.
           
            Nevertheless, we recommend to return to the geometric 
            point of view and 
         think in terms of formal graded
          $Q$-manifolds. This naturally 
         leads to the following

\proclaim Definition.  An $\L$-morphism between two $\L$-algebras $\g_1$ and 
$\g_2$ 
is a pre-$\L$-morphism $\F$ such that the associated morphism 
 $\F_*
 :C(\g_1[1])\ra C(\g_2[1])$
 of graded cocommutative coalgebras, is compatible with codifferentials. 
 \par
 
In geometric terms, an $\L$-morphism corresponds to a 
 $Q$-equivariant map between two formal graded manifolds
 with base points.
 
 For the case of differential graded Lie algebras  
 a pre-$\L$-morphism $\F$ is an $\L$-morphism iff it satisfies
  the following equation for any 
 $n=1,2\dots\,\,\,\,$ and homogeneous elements $\ga_i\in \g_1$:
$$ d\F_n(\ga_1\wedge \ga_2\wedge\dots\wedge\ga_n)-\sum_{i=1}^n \pm \F_n(
\ga_1\wedge\dots \wedge d\ga_i\wedge\dots \wedge \ga_n)=$$
$$={1\over 2}\sum_{k,l\ge 1,\,\,k+l=n} {1\over k!l!}\sum_{\sigma\in 
\Sigma_n}
\pm [\F_k(\ga_{\sigma_1}\wedge\dots\wedge\ga_{\sigma_k}),\F_l(\ga_
{\sigma_{k+1}}\wedge
\dots\wedge \ga_{\sigma_n})]+\sum_{i<j}\pm \F_{n-1}([\ga_i,\ga_j]
\wedge\ga_1
\wedge\dots\wedge\ga_n)\,\,\,.$$

Here are first two equations in the explicit form:
$$d \F_1(\ga_1)=\F_1(d\ga_1)\,,$$
$$d\F_2(\ga_1\wedge\ga_2)-\F_2(d\ga_1\wedge\ga_2)-(-1)^{\overline{\ga_1}}\F_2
(\ga_1\wedge d\ga_2)=
\F_1([\ga_1,\ga_2])-[\F_1(\ga_1),\F_1(\ga_2)]\,\,\,.$$

We see that $\F_1$ is a morphism of complexes. The same is true for
 the case of general $\L$-algebras. The graded space $\g$ for 
 an $\L$-algebra $(\g,Q)$ can be considered as the tensor product
 of $\k[-1]$ with the 
  tangent space 
  to the corresponding formal graded manifold at the base point.
   The differential $Q_1$ on $\g$ comes from the action of $Q$
  on the manifold.

      Let us assume that $\g_1$ and $\g_2$ are differential graded Lie 
      algebras, and $\F$ is an $\L$-morphism from $\g_1$ to 
      $\g_2$. Any solution $\ga\in\g_1\otimes \m$
      of the Maurer-Cartan equation  
      where $\m$ is a nilpotent non-unital algebra,
      produces a solution of the Maurer-Cartan equation in $\g_2\otimes\m$:
      $$d\ga+{1\over 2}[\ga,\ga]=0
      \Longrightarrow
       d\widetilde{\ga}+{1\over 2}
     [\widetilde{\ga},\widetilde{\ga}]=0{\rm{\,\,\,\,\,\,where\,\,\,\,}}
      \widetilde{\ga}=\sum_{n=1}^{\infty} {1\over n!} \F_n(\ga
     \wedge\dots\wedge\ga)\in
     \g_2^1\otimes\m\,\,\,.$$
      
      The same formula is applicable to solutions of the Maurer-Cartan
       equation depending 
     formally on  parameter $\hbar$:
     $$\ga(\hbar)=\ga_1\hbar+\gamma_2\hbar^2+\dots\in \g_1^1[[\hbar]]
     ,\,\,\,d\ga(\hbar)+{1\over 2}
     [\ga(\hbar),\ga(\hbar)]=0\Longrightarrow
      d\widetilde{\ga(\h)}+{1\over 2}
     [\widetilde{\ga(\h)},\widetilde{\ga(\h)}]=0 \,\,\,.$$
      
      The reason why it works is that the Maurer-Cartan equation
       in any differential graded Lie algebra $\g$
       is the equation for the subscheme of zeroes of $Q$ in 
        formal manifold $\g[1]_{formal}$.  $\L$-morphisms map
         zeroes of $Q$ to zeroes of $Q$ because they commute with $Q$.
          We will see in 4.5.2
           that $\L$-morphisms induce natural transformations
           of deformation functors.

 \vskip 0.5truecm
\par\noindent {\bf 4.4. Quasi-isomorphisms}
\vskip 0.5truecm

    $\L$-morphisms generalize usual morphisms of differential graded 
    Lie algebras. In particular, the first Taylor coefficient
     of an $\L$-morphism from $\g_1$ to $\g_2$ is a morphism of complexes
      $(\g_1,Q_1^{(\g_1)})\ra(\g_2,Q_1^{(\g_2)})$ where
       $Q_1^{(\g_i)}$ are the first Taylor coefficients of
         vector fields $Q^{(\g_i)}$ (which we denoted before simply by $Q$). 
      
   \proclaim Definition. A quasi-isomorphism is an $\L$-morphism $\F$ such that
    the first 
   component $\F_1$
    induces isomorphism between cohomology groups  
    of complexes $(\g_1,Q_1^{(\g_1)})$ and $(\g_2,Q_1^{(\g_2)})$.
    \par

    The essence of the homotopy/deformation theory is contained in the 
    following
    
    \proclaim Theorem. Let $\g_1,\g_2$ be two $\L$-algebras
    and $\F$ be an $\L$-morphism from $\g_1$ to 
    $\g_2$. Assume that
      $\F$ is a quasi-isomorphism. Then there exists an $\L$-morphism 
      from $\g_2$ to $\g_1$ inducing
      the inverse isomorphism between cohomology of complexes 
      $(\g_i,Q_1^{(\g_i)})\,\,i=1,2$.
       Also, for the case of differential graded algebras,
        $\L$-morphism $\F$ induces an isomorphism between 
        deformation functors associated with $\g_i$.\par
       
       The first part of this theorem shows that if $\g_1$ is 
       quasi-isomorphic to $\g_2$ then
        $\g_2$ is quasi-isomorphic to $\g_1$, i.e. we get an 
        equivalence relation.
       
       The isomorphism between deformation functors
        at the second part of the theorem is given by last formulas
         from 4.3.

  This theorem is essentially standard (see related results in [GM], [HS1],
   [SS]). 
   Our approach consists in
     the translation of all relevant notions to the geometric
      language of formal graded pointed $Q$-manifolds.   
     
 \vskip 0.5truecm
\par\noindent {\bf 4.5. A sketch of the proof}
\vskip 0.5truecm

\par\noindent {\bf 4.5.1. Homotopy classification of $\L$-algebras}
\vskip 0.5truecm   
        
     Any complex of vector spaces can be decomposed into the direct
      sum of a complex with trivial differential and a contractible complex.
       There is an analogous decomposition in the non-linear case.
       
       \proclaim Definition. An $\L$-algebra $(\g,Q)$ 
       is called minimal if the first Taylor
        coefficient $Q_1$ of the coderivation $Q$ vanishes.
        \par
        
        The property of being formal is invariant under $\L$-isomorphisms.
        Thus, one can speak about minimal formal graded pointed
         $Q$-manifolds.
        
         \proclaim Definition. An $\L$-algebra $(\g,Q)$ is called 
         linear contractible
         if higher Taylor coefficients $Q_{\ge 2}$ vanish and
          the differential $Q_1$ has trivial cohomology.
          \par
          
          The property of being linear contractible is not $\L$-invariant.
          One can call formal graded pointed
           $Q$-manifold {\it contractible} iff
           the corresponding 
            differential graded coalgebra is $\L$-isomorphic
            to a linear contractible one.
            
          \proclaim Lemma. Any $\L$-algebra $(\g,Q)$
          is $\L$-isomorphic to
           the direct sum of a minimal and of a linear 
           contractible $\L$-algebras.
           \par
           
           This lemma says that there exists 
           an affine structure on a formal graded 
          pointed  manifold in which the odd vector field $Q$ has
             the form of a direct sum of a minimal and a linear contractible 
             one. This affine structure can be 
             constructed by induction in the degree of the Taylor expansion.
              The base of the induction is the decomposition of
              the complex $(\g,Q_1)$ into the direct sum of
              a complex with vanishing differential and a complex with trivial
               cohomology.  We leave details of the proof of the lemma
                to the reader. \qed

           As a side remark, we  mention  analogy
           between this lemma and a theorem from singularity theory
            (see, for example, the beginning of 11.1 in [AGV]):
            for every germ $f$ of analytic function
            at  critical point
             one can find local coordinates $(x^1,\dots,x^k,y^1,\dots,y^l)$ 
             such that
             $f=constant+Q_2(x)+Q_{\ge 3}(y)$ where $Q_2$ is a nondegenerate
             quadratic form in $x$ and $Q_{\ge 3}(y)$ is a germ of a function 
             in $y$ such that its  Taylor expansion
              at $y=0$ starts at terms of degree at least $3$.
              
            Let $\g$ be an $\L$-algebra and $\g^{min}$ be a minimal
             $\L$-algebra as in the  previous lemma. Then 
             there are two $\L$-morphisms 
             (projection and inclusion)
             $$(\g[1]_{formal},0)\ra (\g^{min}[1]_{formal},0),\,\,\,
             (\g^{min}
             [1]_{formal},0)\ra (\g[1]_{formal},0)$$
              which are both quasi-isomorphisms.
             From this follows that if 
             $$(\g_1[1]_{formal},0)
             \ra (\g_2[1]_{formal},0)$$
              is a quasi-isomorphism
             then there exists a quasi-isomorphism 
             $$(\g_1^{min}
             [1]_{formal},0)\ra (\g_2^{min}[1]_{formal},0)\,\,\,.$$
           Any quasi-isomorphism between minimal $\L$-algebras
           is  invertible, because it induces an isomorphism
           of spaces of  cogenerators (the inverse mapping theorem from 4.1). 
            Thus,  we proved the first
             part of the theorem. Also, we see that
              the set equivalence classes
            of $\L$-algebras up to quasi-isomorphisms can be naturally 
            identified with the set of equivalence classes of minimal 
            $\L$-algebras up to $\L$-isomorphisms.

    \vskip 0.5truecm
\par\noindent {\bf 4.5.2. Deformation functors at fixed points of $Q$}
\vskip 0.5truecm  

               The deformation functor can be defined in terms
                of a formal graded $Q$-manifold $M$ with  base
                point (denoted by $0$). 
                  The set of solutions
                of the Maurer-Cartan equation with coefficients
                 in a finite-dimensional nilpotent non-unital algebra
                  $\m$ is defined as
                   the set of $\m$-points of the formal
                  scheme of zeroes of $Q$:
                  $$Maps\left(\bigl(
                  Spec(\m\oplus \k\cdot 1), \,{\rm base\,\,\,point}\,\bigr),
                  \bigl(Zeroes(Q),0\bigr)\right)\subset
                  Maps\left(\bigl(
                  Spec(\m\oplus \k\cdot 1), \,{\rm base\,\,\,point}\,\bigr),
                  \bigl(M,0\bigr)\right)
                   \,\,\,.$$
                   
                   In terms of the coalgebra $\c$ corresponding to $M$ this
                    set is equal to the set of
                     homomorphisms of coalgebras
                     $\m^*\ra \c$ with the image annihilated by $Q$. 
                       Another way to say is to introduce a {\it global}
                       (i.e. not formal) pointed $Q$-manifold of maps
                       from $\bigl(
                  Spec(\m\oplus \k\cdot 1), \,{\rm base\,\,\,point}\,\bigr)$
                   to $(M,0)$ and consider zeroes of the global vector
                    field $Q$ on it.
                  
                  Two  solutions
                   $p_0$ and $p_1$ of the Maurer-Cartan equation 
                   are called gauge equivalent iff there
                   exists (parametrized by $Spec(\m\oplus\k\cdot 1)$)
                    polynomial family
                    of odd vector fields $\xi(t)$ on $M$ 
                    (of degree $-1$ with respect to 
                    $\Z$-grading) and a polynomial
                    solution of the equation
                    $${d p(t)\over dt}= [Q,\xi(t)]_{|p(t)},\,\,\,
                    p(0)=p_0,\,\,p(1)=p_1,$$
                    where $p(t)$ is a polynomial family of $\m$-points
                    of  formal graded manifold $M$ with  base point.

                    In terms of $\L$-algebras, the set 
                      of polynomial paths  $\{p(t)\}$ is naturally identified
                       with     $\g^1\otimes \m\otimes \k[t]$. Vector fields 
                       $\xi(t)$ depending polynomially on $t$ are not
                       necessarily 
                   vanishing at the base point $0$. 
                   The set of these vector fields
                   is 
                  $$Hom_{Graded^\k}
                  \left( C(\g[1])\oplus (\k\cdot 1)^*, \g\right
                  )\otimes (\m\oplus\k\cdot 1)\,\,\,.$$ 
                  
         One can check that the gauge equivalence defined above is
          an equivalence relation. Alternatively, one can define the 
          equivalence relation as the transitive closure of the relation 
          from above.
           For formal graded pointed manifold $M$
           we define  set $Def_M(\m)$ as the set of gauge equivalence
           classes of solutions of the Maurer-Cartan equation.
            The correspondence $\m\mapsto Def_M(\m)$ extends naturally
            to a functor denoted also by
            $Def_M$. Analogously, for $\L$-algebra $\g$ we denote by
             $Def_\g$ the corresponding deformation functor.

            One can easily prove the following properties:

            1) for a differential graded Lie algebra $\g$ the deformation 
            functor defined as above for $(\g[1]_{formal},0)$, 
            is naturally equivalent to the 
            deformation
             functor defined in 3.2,

            2) any $\L$-morphism gives a natural transformation of functors,

           3) the functor $Def_{\g_1\oplus\g_2}$
             corresponding
              to the direct sum of two $\L$-algebras, is naturally equivalent
               to the product of functors $Def_{\g_1}\times Def_{\g_2}$,

            4) the deformation functor for a linear
               contractible $\L$-algebra $\g$ is
               trivial, $Def_\g(\m)$ is  a one-element set for every $\m$.

              Properties 2)-4) are just trivial, and 1) is easy.
               It follows from  properties 1)-4) that
               if an $\L$-morphism of differential graded Lie algebras
               is a quasi-isomorphism, then it induces an isomorphism
                of deformation functors.
                The theorem is proven. \qed
              
           We would like to notice here that in the definition of
            the deformation functor one can consider just a formal
           pointed  {\it super} $Q$-manifold $(M,0)$ 
             (i.e. not a graded one), and $\m$ could be a finite-dimensional
              nilpotent  differential super commutative 
              associative non-unital algebra.

 \vskip 0.5truecm
\par\noindent {\bf 4.6. Formality}
\vskip 0.5truecm

\par\noindent {\bf 4.6.1. Two differential graded Lie algebras}
\vskip 0.5truecm 

   Let $X$ be a smooth manifold. We associate with it two differential
    graded Lie algebras over $\R$.
 The first differential graded
  Lie algebra $D_{poly}(X)$ is a subalgebra of the shifted Hochschild 
 complex of the algebra $A$ 
 of functions on $X$ (see 3.4.2).
  The space $D^n_{poly}(X),\,\,\,n\ge -1$ consists of local Hochschild 
  cochains $A^{\otimes (n+1)}\ra A$ given 
  by polydifferential
   operators. In local coordinates $(x^i)$ any element of $D^n_{poly}$ 
   can be written as 
   $$f_0\otimes\dots\otimes f_n\mapsto\sum_{(I_0,\dots,I_n)} 
   C^{I_0,\dots,I_n}(x)\cdot 
   \partial_{I_0}(f_0)
   \dots\partial_{I_n}(f_n)$$
   where the sum is finite, $I_k$ denote multi-indices, $\partial_{I_k}$ 
   denote corresponding 
   partial derivatives, and $f_k$ and $C^{I_0,\dots,I_n}$ are 
   functions in $(x_i)$.

   The second differential graded Lie algebra, $T_{poly}(X)$ is the graded
   Lie 
   algebra of polyvector fields on $X$:
   $$T^n_{poly}(X)=\Gamma(X,\wedge^{n+1} T_X),\,\,\,n\ge -1$$
   endowed with the standard Schouten-Nijenhuis bracket and with 
   the differential $d:=0$. We remind here the formula for this bracket:
   $${\rm for\,\,\,}k,l\ge0\,\,\,\,
   \,\,\,[\xi_0\wedge\dots\wedge\xi_k,\eta_0\wedge\dots\wedge\eta_l]=$$
   $$=
   \sum_{i=0}^k\sum_{j=0}^l (-1)^{i+j+k}[\xi_i,\eta_j]\wedge\xi_0\wedge
   \dots\wedge\xi_{i-1}\wedge\xi_{i+1}\wedge\dots\wedge
   \xi_k\wedge\eta_0\wedge\dots\wedge \eta_{j-1}\wedge\eta_{j+1}\wedge\dots
   \wedge\eta_l,\,\,\,\xi_i,\eta_j\in\G(X,T_X)\,,
   $$
   $${\rm for\,\,\,} k\ge 0\,\,\,\,\,[\xi_0\wedge\dots\wedge\xi_k,h]=$$
   $$=
   \sum_{i=0}^k (-1)^i \xi_i(h)\cdot\bigl(\xi_0\wedge
   \dots\wedge\xi_{i-1}\wedge\xi_{i+1}\wedge\dots\wedge
   \xi_k\bigr)\,\,\,h\in \G(X,{\cal O}_X),\,\,\,\xi_i\in \G(X,T_X)\,\,\,.$$
   
 In local coordinates $(x^1,\dots ,x ^d)$, if one replaces $\p/\p x^i$
  by odd variables $\psi_i$ and writes polyvector fields as functions
   in $(x^1,\dots,x^d|\psi_1,\dots,\psi_d)$, the bracket is
   $$[\ga_1,\ga_2]=\ga_1\bullet \ga_2-(-1)^{k_1 k_2} \ga\bullet \ga_1$$
     where introduce the following notation:
     $$\ga_1\bullet \ga_2:= \sum_{i=1}^d
     {\p \ga_1\over\p \psi_i}{\p\ga_2\over \p x^i},\,\,\,\ga_i\in
      T^{k_i}(\R^d) \,\,\,.$$

We have an evident map $\U_1^{(0)}:T_{poly}(X)\ra D_{poly}(X)$:
$$\U_1^{(0)}:(\xi_0\wedge\dots\wedge \xi_n)\mapsto\left(f_0\otimes\dots 
\otimes f_n\mapsto {1\over (n+1)!}
\sum_{\sigma\in\Sigma_{n+1}} sgn(\sigma) \prod_{i=0}^{n} \xi_{\sigma_i}
(f_i)\right)\,,\,\,\,\,{\rm for}\,\,\, n\ge 0,$$
$$h\mapsto \bigl(1\mapsto h\bigr)\,,\,\,\,h\in\G(X,{\cal O}_X)\,\,\,.$$

 \vskip 0.5truecm
\par\noindent {\bf 4.6.1.1. It is a quasi-isomorphism}
\vskip 0.5truecm 

 \proclaim Theorem. $\U_1^{(0)}$ is a quasi-isomorphism of complexes.
 \par
 
 This is a version of Kostant-Hochschild-Rosenberg theorem which says
 that for a smooth affine algebraic variety $Y$ over  field $\k$ of 
 characteristic zero, the Hochschild cohomology of  algebra ${\cal O}(Y)$
 coincides with the space $\oplus_{k\ge 0}
   \G(X,\wedge ^k T_Y)[-k]$ of algebraic polyvector fields on $Y$.
    Analogous statement for $C^{\infty}$ manifolds
     seems to be well known, although we were
     not able to find  it in the literature. In any case, we show here a proof.
   
   {\bf Proof:} First of all, one can immediately check that the image of
    $\U_1^{(0)}$ is annihilated by the differential in $D_{poly}(X)$,
     i.e. that $\U_1^{(0)}$ is a morphism of complexes.
     
      Complex $D_{poly}(X)$ is filtered by the total degree of 
      polydifferential operators. Complex $T_{poly}(X)$ endowed with zero 
      differential
       also carries a very simple  filtration (just by degrees), 
        such that $\U_1^{(0)}$
        is compatible with filtrations. We claim that 
        $$Gr\bigl(\U_1^{(0)}\bigr):\,Gr\bigl(T_{poly}(X)\bigr)
        \ra \,Gr\bigl(D_{poly}(X)\bigr)$$
         is a quasi-isomorphism. In the graded complex 
         $Gr\bigl(D_{poly}(X)\bigr)$ 
         associated with the filtered complex
          $D_{poly}(X)$ all components are sections of some natural
           vector bundles on $X$, and the differential is
          $A$-linear, $A=C^{\infty}(X)$. 
          The same is true by trivial reasons for 
           $T_{poly}(X)$. Thus, we have to check that the map 
           $Gr\bigl(\U_1^{(0)}\bigr)$  
            is a quasi-isomorphism fiberwise.
            
            Let $x$ be a point of $X$ and $T$ be the tangent space
             at $x$. Principal symbols of polydifferential operators
              at $x$ lie in vector spaces
              $$Sym(T)\otimes\dots\otimes Sym(T)\,\,\,\,\,\,(
              n\,\,\,{\rm times},\,\,\,n\ge 0)$$
               where $Sym(T)$ is the free polynomial algebra generated by $T$.
              It better to identify $Sym(T)$ with the cofree
               cocommutative coassociative coalgebra {\it with} counit
                cogenerated by $T$:
                $$\c:=C(T)\oplus (\k\cdot 1)^*\,\,\,. $$
              $Sym(T)$ is naturally isomorphic to the algebra
                 of differential operators on $T$ with constant coefficients.
                If $D$ is such an operator then it defines a linear functional
                 on the algebra of formal power series at $0\in T$:
                 $$f\mapsto \bigl(D(f)\bigr)(0)\,\,\,.$$ 
                 
                 We denote by $\Delta$ the coproduct in coalgebra $\c$. 
                 It is easy to see that differential in the complex 
                 $Gr\bigl(D_{poly}(X)\bigr)$ in the fiber at $x$ is 
                 the following:
                 $$d: \otimes ^{n+1}\c\ra \otimes^{n+2}\c,\,\,\,
                 d=1^*\otimes
                 id_{\otimes^{n+1}\c}
                 -\sum_{i=0}^n \,(-1)^i\,id\otimes\dots\otimes 
                 \Delta_i\otimes
                 \dots\otimes id+(-1)^n id_{\otimes^{n+1}\c}\otimes 1^*
                  $$
                 where $\Delta_i$ is 
                  $\Delta$ applied to the $i$-th argument.

            \proclaim Lemma. Let $\c$ be the
             cofree cocommutative coassociative
             coalgebra with counit cogenerated by vector space $T$. 
             Then the natural 
              homomorphism of complexes
              $$\bigl(\wedge ^{n+1} T,\,{\rm differential\,}=0\bigr)
              \ra \bigl(\otimes ^{n+1} \c,\,{\rm differential\,\,\,
              as\,\,\,above\,}\bigr)$$
              is a quasi-isomorphism. \par
              
              What we consider is one
             of standard complexes in homological algebra. One of possible
              proofs is the
              following:
              
            {\bf Proof:} first of all, we can safely assume that
             $T$ is finite-dimensional. 
             Let us decompose complex $\bigl(\otimes ^{n+1} 
              \c\bigr)$
               into the infinite direct sum of subcomplexes consisting 
               of tensors
                of fixed total degrees
                 (homogeneous components with respect to the action
                  of the Euler vector fields on $T$). Our statement 
                  means in particular
              that for only finitely many
                 degrees these subcomplexes have non-trivial cohomology. 
                 Thus,
                   the statement of the lemma is true iff the analogous 
                   statement
                    holds when infinite sums  are replaced by infinite 
                    products
                    in the decomposition
                     of $\bigl(\otimes ^{n+1} \c\bigr)$. Terms of the 
                     completed complex
                      are spaces  $Hom(A^{\otimes (n+1)},\k)$ where $A$ 
                      is the algebra
                       of polynomial functions on $T$. It is easy to see 
                       that the 
                        completed complex calculates groups 
                        $Ext^{n+1}_{A-mod}(\k,\k)
                    = \wedge^{n+1} T$ where $1$-dimensional space 
                    $\k$ is
                     considered as $A$-module (via values of polynomial at 
                     $0\in T$) and has
                     a resolution 
                     $$\dots\ra A\otimes A\ra A\ra 0\ra\dots$$
            by free $A$-modules. \qed
            
            As a side remark, we notice that the statement of the lemma 
            holds also
             if one replaces $\c$ by $C(T)$ (i.e. the free coalgebra without 
             counit) and removes terms with $1^*$ from the differential.
              In the language of Hochschild cochains it means that the 
              subcomplex of
               {\it reduced} cochains is quasi-isomorphic to the total 
               Hochschild complex.
               
               The lemma implies that $gr\bigl(\U_1^{(0)}\bigr)$  
            is an isomorphism fiberwise, and the theorem is proven. \qed
  
\vskip 0.5truecm
\par\noindent {\bf 4.6.2. Main theorem}
\vskip 0.5truecm 
 
 Unfortunately, the  map $\U_1^{(0)}$ 
  does not commute with Lie brackets, the Schouten-Nijenhuis bracket 
  does not go to the Gerstenhaber bracket. We claim that 
  this defect can be cured:
   
  \proclaim Theorem. There exists an $\L$-morphism $\U$ from 
  $T_{poly}(X)$ to 
 $D_{poly}(X)$ such that 
  $\U_1=\U_1^{(0)}$. \par
  
  In other words, this theorem says that $T_{poly}(X)$ and $D_{poly}(X)$ are 
  quasi-isomorphic differential
  graded Lie algebras. In analogous situation in rational homotopy theory 
  (see [Su]), a differential graded
  commutative algebra is called formal if it is quasi-isomorphic to  its 
  cohomology algebra endowed with zero differential. 
   This explains the name of subsection 4.6.
   
   The quasi-isomorphism $\U$ in the theorem is not canonical. We will 
   construct explicitly a family of quasi-isomorphisms  parametrized
    in certain sense by a contractible space. It means that
    our construction is canonical up to (higher) homotopies.

   Solutions of the Maurer-Cartan equation in $T_{poly}(X)$ are exactly 
   Poisson structures on $X$:
   $$\alpha
      \in T_{poly}^1(X)=\G(X,\wedge^2 T_X),\,\,\,[\alpha,\alpha]=0\,\,\,.$$
      Any such $\alpha$ defines also a solution formally depending on $\h$,
        $$\ga(\hbar):=\alpha\cdot\h\in T_{poly}^1(X)[[\h]]\,,\,\,\,\,
        [\ga(\hbar),\ga(\hbar)]=0\,\,\,.$$
    The gauge group action is the action of the  diffeomorphism group by 
    conjugation. Solutions of the Maurer-Cartan equation in
     $D_{poly}(X)$ formally depending on $\h$ are star-products.
      Thus, we obtain as a corollary that any Poisson structure 
       on $X$ gives a canonical equivalence class of star-products, and the theorem from
        1.3. 
    
    The rest of the paper is devoted to the proof of Theorem 4.6.2, and to the 
    discussion of various applications, corollaries and extensions.
     In the next section (5)
      we will make some preparations for the universal
      formula (section 6) for an $\L$-morphim from $T_{poly}(X)$ to 
      $D_{poly}(X)$ in the case
       of flat space $X=\R^d$. In section 7 we extend our construction
        to general manifolds.
        
        \vskip 0.5truecm
\par\noindent {\bf 4.6.3. Non-uniqueness}
\vskip 0.5truecm 

There are other quasi-isomorphisms between $T_{poly}(X)$ and $D_{poly}(X)$
 which differ essentially  from the quasi-isomorphism $\U$, i.e. not even
  homotopic in a natural sense
  to $\U$. By {\it homotopy} here we mean the following.
   $\L$-morphisms from one $\L$-algebra to another 
   can be identified with fixed points of $Q$ on  infinite-dimensional
    supermanifold of maps. 
 Mimicking constructions and definitions form 4.5.2 one can define an 
 equivalence relation (homotopy equivalence) on the set of $L$-morphisms. 
  
  Firstly, the multiplicative group $\R^{\times}$ acts by automorphisms
   of $T_{poly}(X)$, multiplying elements $\ga\in T_{poly}(X)^k$ by $const ^k$.
   Composing these automorphisms with $\U$ one get a family of 
   quasi-isomorphisms. 
    Secondly, in [Ko2] we constructed an exotic infinitesimal
    $\L$-automorphism of $T_{poly}(X)$
     for the case $X=\R^d$ which probably extends to  general manifolds.
      In particular, this exotic automorphism produces a vector field
       on the ``set of Poisson structures''.
     The evolution with respect to  time $t$  is described by the following 
           non-linear partial differential equation:
        $${d\alpha\over dt}:=\sum_{i,j,k,l,m,k',l',m'}{\partial^3\alpha^{ij}
        \over \partial x^k\partial x^l\partial x^m} 
        {\partial\alpha^{k k'}\over \partial x^{l'}}
        {\partial\alpha^{l l'}\over \partial x^{m'}}
         {\partial\alpha^{m m'}\over \partial x^{k'}}
         \left({ \partial_i}
         \wedge { \partial_j}\right)$$
         where $\alpha=\sum_{i,j} \alpha^{ij}(x){\partial_i}
         \wedge {\partial_j}$ is a bi-vector field on $\R^d$.
         
         {\it A priori} we can guarantee the existence of a solution of the evolution 
         only for small times 
         and real-analytic initial data.
          One can show that 1) this evolution preserves the class
          of (real-analytic) Poisson structures, 2) if two Poisson structures 
          are conjugate by 
          a real-analytic diffeomorphism then the same will hold after the evolution.
           Thus, our evolution operator is essentially intrinsic and does not depend on 
           the choice of
            coordinates.
         
     Combining it with the action of $\R^{\times}$ as above we see that the Lie algebra
   $aff(1,\R)$   of infinitesimal affine transformations of the line $\R^1$
    acts non-trivially on the space of homotopy classes of quasi-isomorphisms between
    $T_{poly}(X)$ and $D_{poly}(X)$. Maybe, there are other
     exotic $\L$-automorphisms, this possibility is not ruled out yet.
      The reader could ask why our quasi-isomorphism
     $\U$ is better than others. Probably, the answer is that only $\U$ (up to homotopy)
      preserves an additional structure present in the deformation quantization, 
      the cup-product on the tangent cohomology (see section 8).

    \vskip 0.5truecm
\par\noindent {\bf 5. Configuration spaces and their compactifications}
\vskip 0.5truecm

\par\noindent {\bf 5.1. Definitions}
\vskip 0.5truecm  
  
  Let $n,m$ be non-negative integers satisfying the inequality
   $2n+m\ge 2$. We denote by $Conf_{n,m}$ the product of the configuration 
   space of the
    upper half-plane with the configuration space of the real line:
   $$Conf_{n,m}=\{(p_1,\dots,p_n;q_1,\dots,q_m)|\, p_i\in\H,\,\,q_j \in 
   \R,\,\,
   p_{i_1}\ne p_{i_2}\,\,\,{\rm for} \,\,\,i_1\ne i_2,\,\,\,\,\,
   q_{j_1}\ne q_{j_2}\,\,\,{\rm for}\,\,\,j_1\ne j_2\}\,\,\,.$$
   
   $Conf_{n,m}$ is a smooth manifold of dimension $2n+m$.
    The group $G^{(1)}$ of holomorphic transformations of $\C P^1$
     preserving the upper half-plane
     and the point $\infty$, acts on $Conf_{n,m}$. This group
      is a $2$-dimensional
      connected Lie group, isomorphic to the group of orientation-preserving
      affine transformations of the real line:
     
     $$G^{(1)}=\{z\mapsto az+b| \,a,b\in \R, \,a>0\}\,\,\,.$$
     
      It follows from
      the condition
      $2n+m\ge 2$ that the action of $G^{(1)}$ on $Conf_{n,m}$ is free. 
      The quotient space $C_{n,m}:=Conf_{n,m}/G^{(1)}$
       is a manifold of dimension $2n+m-2$. If $P=(p_1,\dots,p_n;
       q_1,\dots,q_m)$ is a point of $Conf_{n,m}$ then we denote
        by $[P]$ the corresponding point of $C_{n,m}$.
      
      Analogously, we introduce simpler spaces $Conf_n$ and $C_n$ for 
      any
      $n\ge 2$:
      
      $$Conf_n:=\{(p_1,\dots,p_n)|\,p_i\in \C,\,\,\,p_i\ne p_j\,\,\,
      {\rm for} \,\,\,i\ne j\},$$
      $$C_n=Conf_n/G^{(2)},\,\,\,dim(C_n)=2n-3,$$
      where $G^{(2)}$ is a 3-dimensional Lie group,
      $$G^{(2)}=\{z\mapsto az+b|\, a\in \R, b \in \C, \,a>0\}\,\,\,.$$
      
      We will construct  compactifications ${\overline C}_{n,m}$
      of $C_{n,m}$ ( and  compactifications
      ${\overline C}_n$ of $C_n$)
      which are smooth manifolds with corners.
      
       We remind that a manifold with corners (of dimension $d$) 
       is defined analogously to
       a usual manifold with boundary, 
       with the only difference that the manifold
        with corners looks locally as an open part of  closed
         simplicial cone $(\R_{\ge 0})^d$.
        For example, the closed hypercube $[0,1]^d$ is a manifold
         with corners.
         There is a natural smooth
         stratification by faces of any manifold with corners.
        
        First of all, we give one of possible formal
 definitions of the compactification
${\overline C}_{n,m}$ in the case $n\ge 1$.
 With any point $[(p_1,\dots,p_n;q_1,\dots,q_m)]$
 of $C_{n,m}$  we associate a
    collection of $2n(n-1)+nm$ angles with values in $\R/2\pi\Z$:
    $$\left(Arg(p_i-p_j),Arg(p_i-{\overline p_j}), Arg(p_i-q_j)\right)
    \,\,\,.$$
    It is easy to see that we obtain an embedding 
     of $C_{n,m}$ into the torus $(\R /2\pi \Z)^{2n(n-1)+nm}$.
      The space ${\overline C}_{n,m}$
       is defined as the compactification of the image of this embedding.
        Analogously, the compactification
         ${\overline C}_n$ is obtained using
         angles $Arg(p_i-p_j)$.

       One can show
       that open strata of ${\overline C}_{n,m}$ are naturally isomorphic
        to products of manifolds of type $C_{n',m'}$ and $C_{n'}$. In the next 
        subsection we will describe explicitly ${\overline C}_{n,m}$
         as a manifold with corners.
        
        There is a natural action of the permutation group $\Sigma_n$
         on $C_n$, and also of $\Sigma_n\times \Sigma_m$ on $C_{n,m}$.
         This gives us a possibility to define 
           spaces $C_A$ and $C_{A,B}$ for finite sets $A,B$ such 
          that $\#A\ge 2$ or $2\#A+\#B\ge 2$ respectively.
           If $A'\hra A$ and $B'\hra B$ are inclusions of sets
            then there are natural fibrations (forgetting maps) 
            $C_A\ra C_{A'}$ and   
            $C_{A,B}\ra C_{A',B'}$.

    \vskip 0.5truecm
\par\noindent {\bf 5.2. Looking through a magnifying glass}
\vskip 0.5truecm  
   
The definition of the compactification given in the previous subsection
is not descriptive. 
 We are going to explain
  an intuitive idea underlying a direct construction of
  the  compactification ${\overline C}_{n,m}$ as a manifold with corners.
   For more formal
         treatment of compactifications of configuration spaces we refer 
         the reader to [FM] (for the case of smooth algebraic
          varieties).
         
   Let us try  to look
   through a magnifying glass, or better through
   a microscope with arbitrary magnification, on different parts
    of the picture formed by points on $\H\cup \R\subset \C$, 
    and by the 
     line $\R\subset \C$. 
    Here we  use {\it Euclidean} geometry on
     $\C\simeq \R^2$ instead of Lobachevsky geometry.
    
    Before doing this let us first consider the case of a configuration
     on $\R^2\simeq \C$, i.e. without the horizontal line $\R\subset \C$.
     We say that the configuration $(p_1,\dots,p_n)$ is  
     in  {\it standard position} iff 
     
     1) the diameter  of the set $\{p_1,\dots,p_n\}$
      is equal to $1$, and,
      
      2) the center of the minimal circle containing 
      $\{p_1,\dots,p_n\}$ is $0\in \C$.
      
      It is clear that any configuration of $n$ pairwise distinct points
       in the case $n\ge 2$ can be uniquely put to  standard
        position by an element of group $G^{(2)}$.
         The set of configurations in standard position gives
          a continuous section $s^{cont}$
          of the natural projection map
           $Conf_n\ra C_n$.

       For a configuration in  standard position there could be several
     domains where we will need magnification in order to see details.
      These domains are those where at least two points of the
       configuration come too close to each other.

        After an appropriate magnification of any such  domain we again get
         a stable configuration (i.e.the  number of points there
         is at least $2$). 
         Then we can put it again in  standard position and repeat
         the procedure.
         
         In such a way we get an oriented tree $T$
         with one root, and leaves numbered
          from $1$ to $n$. For example, 
          the configuration on the next figure

          \vskip 1cm
         \centerline{
\epsfbox{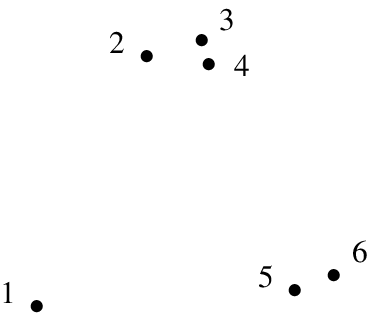}} \vskip 11mm
          gives the tree
          
        \vskip 1cm
         \centerline{
\epsfbox{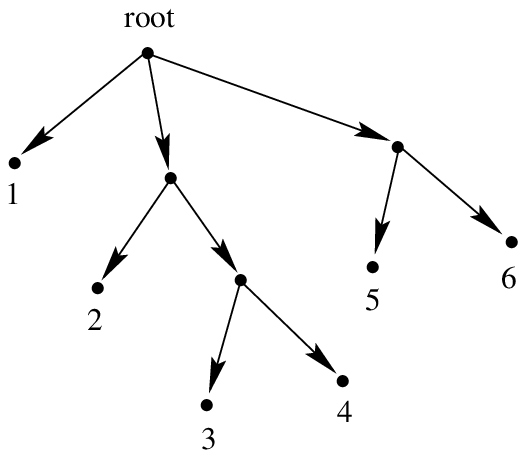}
}
 \vskip 11mm
        
          For every vertex of  tree $T$ except leaves,
          we denote by $Star(v)$ the set of edges starting at $v$.
           For example, in the figure from above the set $Star(root)$
            has three elements, and sets $Star(v)$ for other three
             vertices all have two elements. 
          
          Points in $C_n$ close to one which we consider,
         can be parametrized by the following data:
         
         a) for each vertex $v$ of $T$ except leaves,
          a stable configuration $c_v$ in standard position
          of points
         labeled by the set $Star(v)$,
         
         b) for each vertex $v$ except leaves and
         the root of the tree,
         the scale $s_v>0$ with which we should put a copy
          of $c_v$ instead of the corresponding point $p_v\in \C$ on 
            stable configuration $c_u$ where $u\in V_T$ is such that
            $(u,v)\in E_T$.
            
            More precisely, we act on the configuration $c_v$ by the 
            element $(z\mapsto s_v z+p_v)$ of $G^{(2)}$.
           
           Numbers $s_v$ are small but positive.
            The compactification $\overline{C}_n$ is achieved by 
           formally  permitting  some of scales
           $s_v$ to be equal to $0$. 
           
           In this way we get a compact topological manifold with corners,
            with strata $C_T$
            labeled by trees $T$ (with leaves numbered from $1$
            to $n$). Each stratum $C_T$ is canonically isomorphic to
             the product $\prod_v C_{Star(v)}$ over all vertices $v$ 
             except 
             leaves. In the description as above points of $C_T$
             correspond to collections of configurations with {\it all}
              scales $s_v$ equal to zero. Let us repeat: as a set
              $\overline{C}_n$ coincides with
              $$\bigsqcup_{{\rm trees}\,\,\,T}\,\, \prod_{v\in V_T\setminus
              \{\rm leaves\}} C_{Star(v)}\,\,\,.$$
             
           In order to introduce a smooth structure on 
           $\overline{C}_n$, we should choose a $\Sigma_n$-equivariant
           {\it smooth}
            section $s^{smooth}$ of the projection map $Conf_n\ra C_n$
             instead of the section $s^{cont}$ given by 
             configurations in standard
              position. Local coordinates on $\overline{C}_n$ near a
               given point lying in stratum $C_T$ are scales $s_v\in\R_{\ge 0}$
               close to zero and local coordinates in manifolds 
                $C_{Star(v)}$ for all $v\in V_T\setminus
              \{\rm leaves\}$.
                   The resulting structure of a smooth
               manifold with corners does not depend
                on the choice of section $s^{smooth}$.

              The case of configurations of points on $\H\cup \R$
              is not much harder.
         First of all, we say that a finite non-empty set $S$ 
    of points on $\H\cup\R$
    is in  {\it standard position} iff 
    
    1) the  projection of the convex hull of
     $S$ to the horizontal line $\R\subset \C\simeq \R^2$ is either
     the  one-point set $\{0\}$, or it is an interval
      with the center at $0$,
    
    2) the maximum of the diameter of $S$ and of the distance from
     $S$ to $\R$ is equal to $1$.
     
     It is easy to see that for $2n+m\ge 2$ (the stable case) any 
     configuration
      of $n$ points on $\H$ and $m$ points on $\R$ can be put uniquely
        in standard position by an element of $G^{(1)}$. 
       In order to get a smooth structure, we 
       repeat the same arguments as for the case of manifolds $C_n$.

     Domains where we will need magnification in order to see details,
      are now of two types.
      The first case is when at least two points of the configuration
       come too close to each other. We want to know whether
        what we see is a single point or a collection
         of several points. The second possibility
        is when a point on $\H$ comes too close to $\R$. Here we want also
        to 
        decide whether what we see is a point (or points)
        on $\H$ or on $\R$. 
      
      If the domain which we want to magnify is close to $\R$, then after
        magnification 
       we get again a stable configuration which we can put into the 
       standard
       position. 
        If the domain   is inside $\H$, then after  magnification we get a 
        picture
         without the horizontal line in it, and we are back
          in the situation concerning ${\overline C}_{n'}$ for 
           $n'\le n$.
          
          It is instructional to draw low-dimensional spaces $C_{n,m}$.
           The simplest one, $C_{1,0}={\overline C}_{1,0}$ is just a point.
            The space $C_{0,2}=\OC_{0,2}$ is a two-element set.
            The space $C_{1,1}$ is an open interval, and its closure
             $\OC_{1,1}$ is a closed interval (the real line $\R\subset\C$
              is dashed on the picture):
             
             \vskip 1cm
            \centerline{
\epsfbox{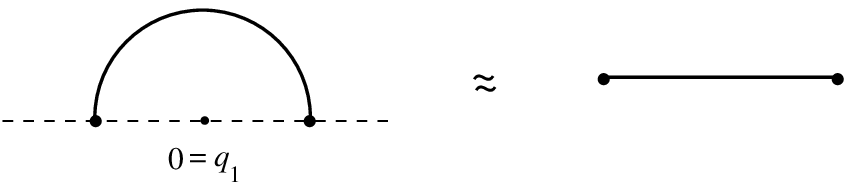}
}
\vskip 11mm
             
             The space $C_{2,0}$ is diffeomorphic 
              to $\H\setminus \{0+1\cdot i\}$. The reason is that
               by  action of $G^{(1)}$ we can put point $p_1$
                to the position $i=\sqrt{-1}\in \H$. The closure
               $\OC_{2,0}$ can be drawn like this:
               
             \vskip 1cm
              \centerline{
\epsfbox{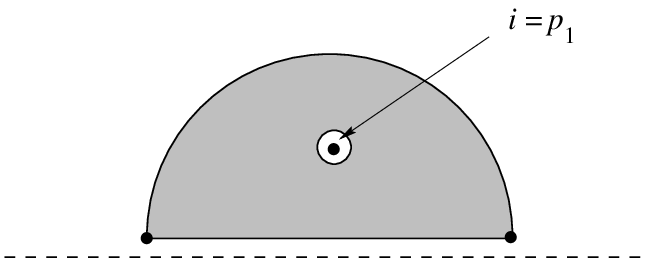}
}  \vskip 11mm
           or like this:

         \vskip 1cm    \centerline{
\epsfbox{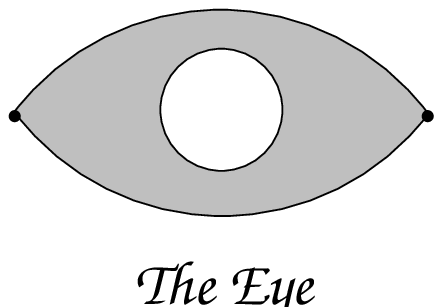}}  \vskip 11mm
              
          Forgetting maps (see the end of 5.1) extend
           naturally to smooth maps of compactified spaces.

            \vskip 0.5truecm
\par\noindent {\bf 5.2.1. Boundary strata}
\vskip 0.5truecm

              We give here the list of all 
              strata in $\OC_{A,B}$ of codimension $1$:
              
              S1) points $p_i\in\H$ for $i\in S\subseteq A$ where
               $\#S\ge 2$, move close to each other but far from $\R$,

                S2) points $p_i\in \H$ for $i\in S\subseteq A$
                 and points $q_j\in\R$ for $j\in S'\subseteq B$
                  where $2\#S+\# S'\ge 2$,
                   all move close to each other
                   and to $\R$, with at least one point left
                    outside $S$ and $S'$, i.e. 
                    $\#S+\#S'\le \#A+\#B-1$.
               
               The stratum of  type S1 is  
               $$\p_S \OC_{A,B}\simeq
                C_{S}\times C_{(A
               \setminus S)\sqcup \{pt\}, B}$$
               where $\{pt\}$ is a one-element set, whose element represents
                the cluster $(p_i)_{i\in S}$ of points in $\H$.
                 Analogously, the stratum of type S2 is
                 $$\p_{S,S'} \OC_{A,B}\simeq
                 C_{S,S'}\times C_{A\setminus S,
                  (B\setminus S')\sqcup
                  \{pt\} }\,\,\,.$$

   \vskip 0.5truecm
\par\noindent {\bf 6. Universal formula}
\vskip 0.5truecm  

In this section we propose a formula for an $\L$-morphism
 $T_{poly}(\R^d)\ra D_{poly}(\R^d)$ generalizing a formula for the star-product
  in section 2.
   In order to write it we need  to make some preparations.

  \vskip 0.5truecm
\par\noindent {\bf 6.1. Admissible graphs}
\vskip 0.5truecm

   \proclaim Definition. Admissible graph $\G$ is 
     an oriented graph with labels
      such that \hfil\break
     1) the set of vertices $V_\G$ is
      $\{1,\dots,n\}\sqcup \{{\overline 1},\dots,{\overline m}
      \}$ where $n,m\in \Z_{\ge 0}$, $2n+2-m\ge 0$;
       vertices from the set
      $\{1,\dots ,n\}$ are called vertices of the first type, vertices
        from $\{ {\overline 1},\dots,{\overline m}\}$ 
        are called vertices of the second type, \hfil\break
       2) every 
       edge $(v_1,v_2)\in E_\G$ starts at a vertex of  
       first type, $v_1\in \{1,\dots,n\}$,\hfil\break
    3) there are no loops, i.e. no edges of the type $(v,v)$,\hfil\break
       4) for every vertex $k\in\{1,\dots,n\}$ of the first type, 
        the set of edges 
        $$Star(k):=\{(v_1,v_2)\in E_\G| \,v_1=k\}$$  
        starting from $s$,
        is labeled by symbols $(e_k^1,\dots,e_k^{\# Star(k)})$. \par

    Labeled oriented graphs 
    considered in section 2 are exactly (after the identifications
     $L=\overline{1},\,\,R=\overline{2}$) admissible graphs
     such that $m$ is equal to $2$,
      and the number of edges starting at every vertex
      of first type is also equal to $2$.
      
  \vskip 0.5truecm
\par\noindent {\bf 6.2. Differential forms on configuration spaces}
\vskip 0.5truecm  

The space $\OC_{2,0}$ (the Eye) is homotopy equivalent to the
 standard circle $S^1\simeq\R/2\pi\Z$. Moreover, one of its 
  boundary components, 
  the space $C_2=\OC_2$, is naturally $S^1$. 
  The other component of the boundary is the union of two
 closed  intervals (copies of $\OC_{1,1}$) with identified end points.
   
  \proclaim Definition. 
  An  angle map is
   a smooth map $\phi: \OC_{2,0}\ra  \R/2\pi\Z\simeq S^1$ 
   such that
    the restriction of $\phi$ to $C_2\simeq S^1$ is the
    angle measured in the anti-clockwise direction from the vertical
     line, and $\phi$ maps the whole upper interval $\OC_{1,1}\simeq [0,1]$ 
     of the Eye, to a point in $S^1$.
     \par
     
     We will denote $\phi([(x,y)])$ simply by $\phi(x,y)$ where
      $x,y\in \H\sqcup \R,\,\,\,x\ne y$.
   It follows from the definition that $d\phi(x,y)
   =0$ if $x$ stays in $\R$.
       
      For example,
      the special map $\phi^{h}$ used in the formula in section 2,
       is an angle map. In the rest of the paper
  (except of the comments 9.5 and 9.8) we can use any $\phi$, not necessarily
   harmonic.

       We are now prepared for the analytic part of the universal formula.
        Let  $\Gamma$ be an admissible graph
        with $n$ vertices of the first type, $m$ vertices
         of the second type and with $2n+m-2$ edges.
           We define the weight of graph $\G$ by the following formula:
        
        $$W_\G:=
       \prod_{k=1}^n {1\over
(\# Star (k))!} {1\over (2\pi)^{2n+m-2}}
        \int\limits_{\OC_{n,m}^+}\,\,\bigwedge _{e\in E_G} 
        d\phi_e\,\,\,.$$
        
        Let us explain what is written here. The domain of integration
         $\OC_{n,m}^+$ is a connected component of $\OC_{n,m}$
          which is the closure of configurations for which
           points $q_j,\,\,1\le j\le m$ on $\R$ are
            placed in the increasing order:
            $$q_1<\dots<q_m\,\,\,.$$
            
           The orientation of $Conf_{n,m}$ is the product of the
            standard orientation on the coordinate space $\R^m\supset
            \{(q_1,\dots,q_m)|\,q_j\in \R\}$, with the product
             of standard orientations on the plane $\R^2$ (for points
              $p_i\in\H
             \subset \R^2$). The group $G^{(1)}$ 
               is even-dimensional and naturally oriented
                because it acts freely and transitively
                 on complex manifold $\H$. Thus, the
              quotient space $C_{n,m}=Conf_{n,m}/G^{(1)}$
              carries again a natural orientation.

                Every edge $e$ of $\G$ defines a map from $\OC_{n,m}$
                 to $\OC_{2,0}$ or to $\OC_{1,1}\subset \OC_{2,0}$
                  (the forgetting map). Here we consider inclusion
                  $\OC_{1,1}$ in $\OC_{2,0}$ as the {\it lower}
                   interval of the Eye. 
                 The pullback of the function $\phi$ 
                  by the  map $\OC_{n,m}\ra \OC_{2,0}$ corresponding to 
                  edge $e$
              is denoted by $\phi_e$.
                  
                  Finally, the ordering in the
                  wedge product of $1$-forms $d\phi_e$
            is fixed by enumeration of the set of
            sources of edges and by the
            enumeration of the set of edges with a given source.
            
            The integral giving $W_\G$ is well-defined because it
             is an integral of a smooth differential form over a compact
             manifold with corners. 
             
   \vskip 0.5truecm
\par\noindent {\bf 6.3. Pre-$\L$-morphisms associated with graphs}
\vskip 0.5truecm            
           
        For any admissible graph $\G$ with $n$ vertices of the first type,
         $m$ vertices of the second type, and $2n+m-2+l$ edges where $l\in \Z$, 
         we define a  
        linear map  $\U_\G:\otimes
        ^n T_{poly}(\R^d)\ra D_{poly}(\R^d)[1+l-n]$.
         This map has only one
          non-zero graded component 
          $(\U_\G)_{(k_1,\dots,k_n)}$ where $k_i=\# Star(i)-1,\,\,
          i=1,\dots,n$. If $l=0$ then from $\U_\G$
           after anti-symmetrization
           we obtain a pre-$\L$-morphism.

            Let $\ga_1,\dots,\ga_n$ be polyvector fields on $\R^d$ of degrees
            $(k_1+1),\dots,(k_n+1)$, and $f_1,\dots, f_m$
              be functions on $\R^d$. We are going to write a  
               formula
             for  function $\Phi$ on $\R^n$:
             $$\Phi:=\left(
             { \U}_\G(\ga_1\otimes\dots\otimes
             \ga_n)\right)(f_1\otimes
             \dots \otimes f_m)\,\,\,.$$

               The formula for $\Phi$ is the sum over all
                configurations of indices running from $1$ to $d$, 
                labeled by $E_\G$:
                $$\Phi=\sum_{I:E_\G\ra \{1,\dots,d\}}\Phi_I\,\,,$$
              where $\Phi_I$ is the product over all $n+m$ vertices of
              $\G$ of certain partial derivatives of functions
               $g_j$ and of coefficients
               of $\ga_i$. 
               
               Namely, with each vertex $i,\,\,1\le i\le n$ of the 
               first type
                we associate function 
               $\psi_i$ on $\R^d$
                which is a coefficient of the polyvector field
                $\ga_i$:
                $$\psi_i=\langle \ga_i, dx^{I(e^1_i)}\otimes\dots
                \otimes dx^{I(e^{k_i+1}_i)}\rangle\,\,\,.$$
                Here we use the identification of polyvector
                 fields with skew-symmetric tensor fields as
                 $$\xi_1\wedge\dots\wedge\xi_{k+1}\ra
                 \sum_{\sigma\in\Sigma_{k+1}} sgn(\sigma)\,
                  \xi_{\sigma_1}\otimes\dots
                  \otimes \xi_{\sigma_{k+1}
                  }\in \G(\R^d,T^{\otimes(k+1)})\,\,\,.$$
                 For each vertex $\overline j$ of second type the associated
                function $\psi_{\overline j}$ is defined as $f_j$.
             
             Now, at each vertex of graph $\G$ we put a function 
             on $\R^d$ (i.e. $\psi_i$
             or $\psi_{\overline j}$). Also, on 
             edges of graph $\G$ there are indices
              $I(e)$ 
              which label coordinates in $\R^d$. In the next step
               we put into each vertex 
              $v$ 
              instead of function $\psi_v$ its partial derivative
               $$\left(\prod_{e\in E_\G,\,e=(*,v)}
               \p_{I(e)}\right) \psi_v,$$
               and then take the product over all vertices $v$ of $\G$.
               The result is by definition the summand $\Phi_I$.
               
              Construction of the function $\Phi$ from the graph
                $\G$, polyvector fields $\ga_i$ and functions $f_j$,
                 is  invariant under the action of the group
                 of affine transformations of $\R^d$ because
                 we contract upper and lower indices.

  \vskip 0.5truecm
\par\noindent {\bf 6.4. Main Theorem for $X=\R^d$, and the proof}
\vskip 0.5truecm      

We define an
$\L$-morphism $\U:T_{poly}(\R^d)\ra D_{poly}(\R^d)$ by the formula for its
$n$-th derivative $\U_n,\,\,\,n\ge 1$ considered as  a skew-symmetric
 polylinear map (see 4.2) from $\otimes
        ^n T_{poly}(\R^d)$ to $D_{poly}(\R^d)[1-n]$:
        
$$\U_n=\sum_{m\ge 0} \,\, 
\sum_{\G\in G_{n,m}} 
W_\G \times \U_\G
\,\,\,.$$
Here $G_{n,m}$ denotes the set of all admissible graphs with $n$ vertices
 of the first type, $m$ vertices in the second group
  and $2n+m-2$ edges, where $n\ge 1, \,m\ge 0\,$ (and automatically
  $2n+m-2\ge 0$).

 \proclaim Theorem. $\U$ is an $\L$-morphism, and also a 
 quasi-isomorphism. \par
 
 {\bf Proof:} 
 first of all, we should check that $\U_n$ is skew-symmetric, i.e. that $\U$
  is a pre-$\L$-morphism. For this see
  subsection 6.5.

   The condition that $\U$ is an $\L$-morhism (see 4.3 and 3.4.2)
   can be written
    explicitly as 
    $$
    f_1\cdot\left(\U_n(\ga_1\wedge\dots\wedge \ga_n)\right)
    (f_2\otimes\dots\otimes f_m)\pm \left(\U_n(\ga_1\wedge\dots\wedge 
    \ga_n)\right)
    (f_1\otimes\dots\otimes f_{m-1})\cdot f_m+$$
    
    $$+\sum_{i=1}^{m-1} \pm\left(\U_n(\ga_1\wedge\dots\wedge \ga_n)\right)
      (f_1\otimes\dots \otimes (f_i f_{i+1})\otimes\dots\otimes f_m)
      +$$
      
    $$+\sum_{i\ne j}\pm 
  \left(\U_{n-1}([\ga_i,\ga_j]\wedge \ga_1\wedge\dots\wedge \ga_n)
  \right)(f_1\otimes\dots\otimes f_m)+$$
  $$+
 {1\over 2}
  \sum_{k,l\ge 1,\,\,k+l=n} {1\over k! l!}\sum_{\sigma\in \Sigma_n}
\pm \left[\U_k(\ga_{\sigma_1}\wedge\dots\wedge\ga_{\sigma_k}),
\U_l(\ga_{\sigma_{k+1}}\wedge
\dots\wedge \ga_{\sigma_n})\right](f_1\otimes \dots\otimes f_m)=0\,\,\,.$$
   
Here $\ga_i$ are polyvector fields, $f_i$ are functions, $\U_n$ are
 homogeneous components of $\U$ (see 4.1).
 There is a way to rewrite this formula. Namely, we define
  $\U_0$ as the map 
  $\otimes ^0(
      T_{poly}(\R^d))\ra D_{poly}(\R^d)[1]$ which maps
      the generator $1$ of 
     $\R\simeq
     \otimes ^0(
      T_{poly}(\R^d))$ to the product $m_A\in D^1_{poly}(\R^d)$ in the algebra
      $A:=C^{\infty}(\R^d)$. Here $m_A:f_1\otimes f_2\mapsto f_1 f_2$ is
      considered as a bidifferential operator.
      
      The condition from above
       for $\U$ to be an $\L$-morphism is equivalent to the following one:
       $$\sum_{i\ne j}\pm 
  \left(\U_{n-1}((\ga_i\bullet \ga_j)\wedge \ga_1\wedge\dots\wedge \ga_n)
  \right)(f_1\otimes\dots\otimes f_m)+$$
  $$+
  \sum_{k,l\ge 0,\,\,k+l=n} {1\over k! l!}\sum_{\sigma\in \Sigma_n}
\pm \left(\U_k(\ga_{\sigma_1}\wedge\dots\wedge\ga_{\sigma_k})\,
\circ \,\U_l(\ga_{\sigma_{k+1}}\wedge
\dots\wedge \ga_{\sigma_n})\right)(f_1\otimes \dots\otimes f_m)=0\,\,\,.$$

Here we use definitions of brackets in $D_{poly}$ and $T_{poly}$ via
 operations $\circ$ (see  3.4.2) and $\bullet$ (see 4.6.1).
We denote the l.h.s. of the expression above as $(F)$. 

$\U+\U_0$ is not a pre-$\L$-morphism because it maps $0$ to a non-zero point
 $m_A$. Still the equation $(F)=0$ makes  sense and means that the map
  $(\U+\U_0)$ from formal $Q$-manifold
   $T_{poly}(\R^d)[1]\bigr)_{formal}$ to the formal neighborhood
    of point $m_A$ in super 
    vector space $D_{poly}(\R^d)[1]$ is $Q$-equivariant,
     where the odd vector field $Q$ on the target is purely quadratic and 
     comes
     from the bracket on $D_{poly}(\R^d)$, forgetting the differential.
     
     Also, the term $\U_0$ comes from the unique graph $\G_0$ which
      was missing in the definition of $\U$. Namely, 
      $\G_0$ has $n=0$ vertices of the first type, $m=2$ vertices of the
      second type, and no edges at all. 
       It is easy to see that $W_{\G_0}=1$ and $\U_{\G_0}=\U_0$.

  We consider the expression $(F)$ simultaneously for all possible
   dimensions $d$. It is clear that one can write
   $(F)$ as a linear combination
  $$\sum_\G c_\G \cdot\U_\G\bigl(\ga_1\otimes\dots\otimes\ga_n)\bigr)(
  f_1\otimes\dots\otimes f_m)$$
  of  expressions  $\U_\G$ for admissible graphs $\G$
   with $n$ vertices of the first type, $m$ vertices
    of the second type, and $2n+m-3$ edges where $n\ge 0,\,m\ge 0,\,
    2n+m-3\ge 0$. 
    We assume that $c_\G=\pm c_{\G'}$ if graph $\G'$ is obtained from $\G$
     by a renumeration of vertices of first type and by a
      relabeling of edges in sets $Star(v)$ (see 6.5 where we discuss 
      signs).
     
     Coefficients $c_\G$ of this linear combination are 
   equal to certain sums
    with signs of weights $\W_{\G'}$ 
      associated with some other graphs $\G'$, and of products
    of two such weights. In particular,
     numbers $c_\G$ do not depend on the dimension
      $d$ in our problem. Perhaps it is better to  
      use here the language of operads, but we will not do it.
       
       We want to check that $c_\G$ vanishes for each $\G$.

   The idea is to identify $c_\G$ with the integral over the boundary
    $\p \OC_{n,m}$ of the closed differential form constructed from $\G$ as
     in 6.2. The Stokes formula gives the vanishing:
     $$
     \int\limits_{\p \OC_{n,m}}\bigwedge_{e\in E_\G}d\phi_e
     =\int\limits_{\OC_{n,m}}d\left(\bigwedge_{e\in E_\G}d\phi_e\right)
     =0\,\,\,.$$
     
      We are going to 
     calculate integrals of the form  $\wedge_{e\in E_\G}d\phi_e$
      restricted to all possible
     boundary  strata of $\p \OC_{n,m}$, and prove that 
       the total integral as above
        is equal to  $c_\G$.
       At the subsection 5.2 we listed two
        groups of boundary strata, denoted by S1 and S2 and labeled by sets
         or pairs of sets. Thus,
      $$0=\,\int\limits_{\p \OC_{n,m}}\bigwedge_{e\in E_\G}d\phi_e
      =\sum_{S}  \int\limits_{\p_S \OC_{n,m}}\bigwedge_{e\in E_\G}d\phi_e
      + \sum_{S, S'}  \int\limits_{\p_{S,S'}
       \OC_{n,m}}\bigwedge_{e\in E_\G}d\phi_e \,\,\,.$$ 
        
   \vskip 0.5truecm
\par\noindent {\bf 6.4.1. Case S1}
\vskip 0.5truecm         
        
        Points $p_i\in\H$
         for $i$ from  subset $S\subset\{1,\dots,n\}$ where
         $\#S\ge 2$, move close to each other. 
           The integral over the stratum $\p_S \OC_{n,m}$
           is equal to the product 
             of an integral over $C_{n_1,m}$ with an integral over
             $C_{n_2}$  where
              $n_2:=\#S,\,\,n_1:=n-n_2+1$. The integral
              vanishes by dimensional
               reasons unless
                the number of edges of $\G$ connecting vertices
              from $S$ is equal to $2n_2-3$.

         There are several possibilities:

 \vskip 0.5truecm
\par\noindent {\bf 6.4.1.1. First subcase of S1: $n_2=2$}
\vskip 0.5truecm  

           In this subcase
             two vertices from $S_1$ are connected exactly by one
           edge, which we denote by $e$. 
           The integral over $C_2$ here gives number $\pm 1$ (after division 
           by $2\pi$ coming from the of the formula for weights $W_{\G}$).
           The total integral over the boundary stratum
            is equal to the integral of a new graph $\G_1$ obtained
             from $\G$ by the contraction of  edge $e$.
             It is easy to see  (up to a sign) that this corresponds
              to the first line in our expression
              $(F)$, the one where the 
               operation $\bullet$ on polyvector fields appears.
               
        \vskip 1cm
         \centerline{
\epsfbox{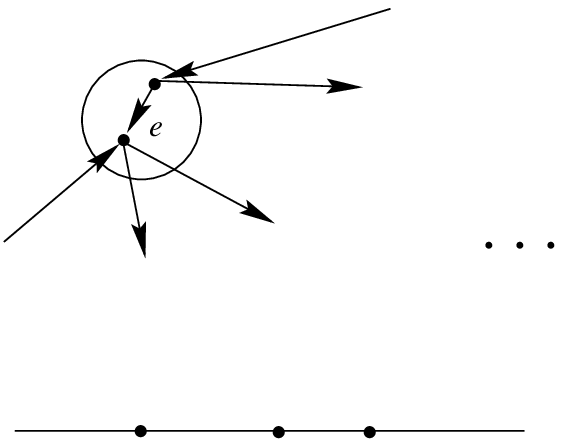}
}   \vskip 11mm

  \vskip 0.5truecm
\par\noindent {\bf 6.4.1.2. Second subcase of S1: $n_2\ge 3$}
\vskip 0.5truecm         

 This is the most non-trivial case.
         The integral corresponding to the corresponding boundary
          stratum
          vanishes because the integral of any product of $2n_2-3$
           angle
           forms over $C_{n_2}$ where $n_2\ge 3$ vanishes, as
            is proven later in 6.6.

      \vskip 1cm   
         \centerline{
\epsfbox{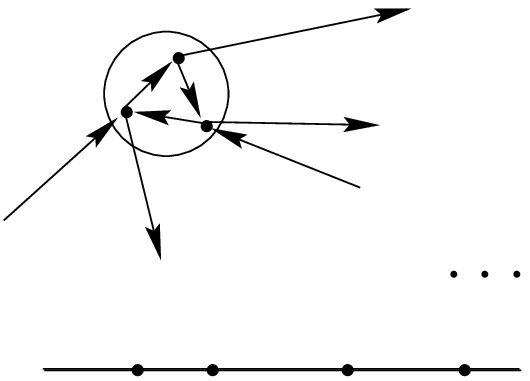}
}\vskip 11mm

   \vskip 0.5truecm
\par\noindent {\bf 6.4.2. Case S2}
\vskip 0.5truecm

           Points $p_i$ for $i\in S_1\subset\{1,\dots,n\}$ and points
                $q_j$ for ${\overline j}\in S_2\subset\{{\overline 1},\dots,
                {\overline m}\}$
                 move close to each other and to the horizontal line $\R$.
                 The condition is that $2n_2+m_2-2\ge 0$ and
                  $n_2+m_2\le n+m-1$ where $n_2:=\#S_1,\,\,m_2:=\#S_2$.
                   The corresponding stratum is isomorphic to
                    $C_{n_1,m_1}\times C_{n_2,m_2}$ where
                     $n_1:=n-n_2,\,\,m_1=m-m_2+1$.
                      The integral of this stratum decomposes into the 
                      product of two integrals. It vanishes if the number of
                      edges of
                      $\G$  connecting vertices from $S_1\sqcup S_2$ is not 
                      equal to $2n_2+m_2-2$.
                      
    \vskip 0.5truecm
\par\noindent {\bf 6.4.2.1. First subcase of S2: no bad edges}
\vskip 0.5truecm                 
             
             In this subcase we assume  that there is no edge $(i,j)$
                in $\G$ such that $i\in S_1, j\in \{1,\dots,n\}
                \setminus S_1$.
                
                The integral over the boundary stratum
                 is equal to the product $W_{\G_1}\times
                 W_{\G_2}$ where $\G_2$ is the restriction of $\G$ to
                  the subset $S_1\sqcup S_2\subset
                  \{1,\dots, n\}\sqcup\{{\overline 1},\dots,
                  {\overline m}\}=V_\G$, and 
                  $\G_1$ is obtained 
                  by the 
                   contraction of all vertices in this set to a new vertex
                   of the second type. 
                   Our condition guarantees that
                    $\G_1$ is an admissible graph.
                     This corresponds to the second line in $(F)$, where
                      the product $\circ$ on polydifferential
                       operators appears.
                      
                     \vskip 1cm  
              \centerline{
\epsfbox{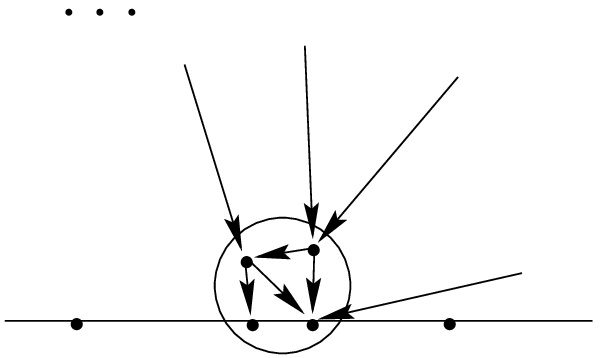}
} \vskip 11mm

       \vskip 0.5truecm
\par\noindent {\bf 6.4.2.2. Second subcase of S2: there is a bad edge}
\vskip 0.5truecm                  
                 
                 Now we assume that there is an edge $(i,j)$
                in $\G$ such that $i\in S_1, j\in \{1,\dots,n\}\setminus 
                S_1$.   In this case the integral is zero because
                      of the condition $d\phi(x,y)=0$ if $x$ stays on the line 
                  $\R$.

           \vskip 1cm    \centerline{
\epsfbox{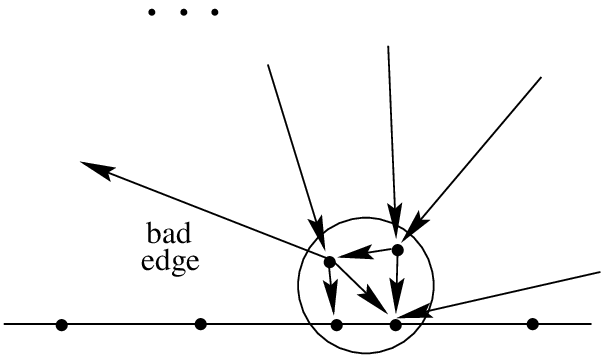}
}    \vskip 11mm

                 The reader can wonder about what happens if after
                  the collapsing the graph will have multiple edges.
                  Such terms do not appear in $(F)$. Nevertheless,
   we ingore them because in this case 
                 the differential form which we integrate
vanishes
     as it contains the square of $1$-form.

   Thus, we see that we exhausted all
    possibilities and get contributions of all terms in the formula $(F)$.
    We proved that $c_{\G}=0$ for any $\G$, and
    that $\U$ is an $\L$-morphism.

    \vskip 0.5truecm
\par\noindent {\bf 6.4.3. We finish the proof of the theorem from 6.4}
\vskip 0.5truecm 

   In order to check that $\U$
    it is a quasi-isomorphism, we should 
   show that its component $\U_1$ coincides with $\U_1^{(0)}$
    introduced in 4.6.1. It follows from definitions that every 
   admissible 
    graph with $n=1$ vertex of first type and $m\ge 0$ vertices of the second 
    type, and with $m$ edges, is the following tree:
       
   \vskip 1cm   \centerline{
\epsfbox{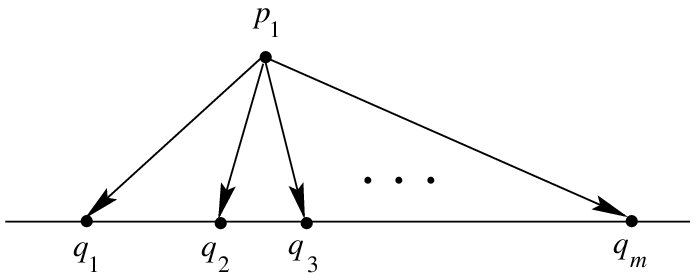}
} \vskip 11mm
       
       The integral corresponding to this graph is $(2\pi)^m /m!$. 
        The map ${\widetilde U}_\G$
         from polyvector fields to polydifferential operators
         is the   one which appears in 4.6.1:
         $$\xi_1\wedge\dots\wedge\xi_m\ra 
         {1\over m!}
         \sum_{\sigma\in \Sigma_m}
          sgn(\sigma)\cdot \xi_{\sigma_1}\otimes\dots \otimes \xi_{\sigma_m},
          \,\,\,\,\xi_i\in\G(\R^d,T)\,\,\,.$$

          The theorem is proven. \qed
          
      \vskip 0.5truecm
\par\noindent {\bf 6.4.4. Comparison with the formula from section 2 }
\vskip 0.5truecm 

The weight $w_\G$ defined in section 2 differ from $W_\G$ defined in 6.2
 by factor $2^n/n!$. On the other hand, the bidifferential operator
  $B_{\G,\a}(f,g)$  is $2^{-n}$ times $\U_\G(\a\wedge\dots
  \wedge \a)(f\otimes g)$. The inverse factorial $1/n!$ appears
   in the Taylor series (see the end of 4.3). Thus, we obtain the formula from
    section 2.

   \vskip 0.5truecm
\par\noindent {\bf 6.5. Grading, orientations, factorials,
 signs}
\vskip 0.5truecm 

Homogeneous components of  $\U$ are maps of graded spaces
 $$Sym^n((\bigoplus_{k\ge 0} \G(\R^d,\wedge^k T)[-k])[2])\ra (\underline{Hom}
 (A[1]^{\otimes m}, A[1]))[1]$$
 where $\underline{Hom}$ denotes the internal $Hom$ in the tensor
 category $Graded^{\k}$. We denote the expression from above by $(E)$. 
First of all,
 in the  expression $(E)$ 
 each polyvector field $\ga_i\in \G(\R^d,\wedge^{k_i}T)$ 
 appears with the shift
 $2-k_i$. In our formula for $\U$ the same $\ga_i$ 
 gives $k_i$ edges of the graph,
 and thus $k_i$ $1$-forms which we have to integrate. Also, it gives $2$
 dimensions for the integration domain $\OC_{n,m}$. Secondly, every
  function $f_j\in A$ appears with shift $1$ in $(E)$ 
   and gives  $1$ dimension to the integration domain.
    We are left with two shifts by $1$ in $(E)$
     which are accounted for $2$ dimensions of the group $G^{(1)}$.
      From this it is clear that our  formula for $\U$ is
compatible with $\Z$-grading. 

Moreover, it is also clear that
 things responsible for various signs in our formulas:\hfil\break
 1) the orientation of $\OC_{n,m}$,\hfil\break
 2) the order in which we multiply $1$-forms $d\phi_e$,\hfil\break
 3) $\Z$-gradings of vector spaces in $(E)$,\hfil\break
 are naturally decomposed into pairs. This implies that the enumeration
  of the set of vertices of $\G$, and also the enumeration of edges
   in sets $Star(v)$ for vertices $v$ of the first type are not really
    used. Thus, we see that $\U_n$ is skew-symmetric.

Inverse factorials $1/(\# Star(v)!)$ kill the summation over enumerations
 of sets $Star(v)$. The inverse factorial $1/n!$ is the final formula
  does not appear
 because we consider higher derivatives which  are already multiplied
  by $n!$.

 The last thing to check is that in our derivation of the fact
  that $\U$ is an $\L$-morphism using Stokes formula we did not loose
   anywhere a sign. This is really hard to explain. How, for example,
    can one compare the standard orientation on $\C$ with shifts by 2
     in $(E)$?
    As a hint to the reader we would like to mention that it is
     very convenient to ``put'' the resulting expression
      $$\Phi:=\left(
             { \U}_\G(\ga_1\otimes\dots\otimes
             \ga_n)\right)(f_1\otimes
             \dots \otimes f_m)$$
  to the point $\infty$ on the absolute. Also, we can not 
  guarantee that we did not make  errors
   when we  write formulas
     in {\it lines}.

 \vskip 0.5truecm
\par\noindent {\bf 6.6. Vanishing of integrals}
\vskip 0.5truecm

 In this subsection we consider the space $C_n$ of $G^{(2)}$-equivalence 
 classes 
 of configurations of points on 
 the Euclidean plane. Every two indices $i,j,\,\,i\ne j, \,1\le i,j\le n$
  give a forgetting map $C_n\ra C_2\simeq S^1$. We denote by
   $d\phi_{i,j}$ the closed $1$-form on $C_n$ which is the pullback of the
    standard $1$-form $d(angle)$ on the circle. We use the same notation
     for the pullback of this form to $Conf_n$.
    
 \proclaim Lemma. Let $n\ge 3$ be an integer. The integral
  over $C_n$ of the product of any $2n-3=dim(C_n)$ closed
   $1$-forms $d\phi_{i_\a,j_\a}$,
   $\a=1,\dots, 2n-3$, 
   is equal to zero. \par

   {\bf Proof:} First of all, we  identify $C_n$ 
    with the subset $C_n'$ of $Conf_n$
     consisting of configurations such that the point $p_{i_1}$
      is $0\in \C$ and $p_{j_1}$ is on the unit circle $S^1\subset \C$.
       Also, we rewrite the form which we integrate as
       $$\bigwedge_{\a=1}^{2n-3}d\phi_{i_\a,j_\a}
       =d\phi_{i_1,j_1}\wedge \bigwedge_{\a=2}
       ^{2n-3} d(\phi_{i_\a,j_\a}-\phi_{i_1,j_1})\,\,\,.$$
      
      Let us map the  space  $C_n'$ onto the space $C''_n\subset Conf_n$
       consisting 
      of
       configurations with $p_{i_1}=0$ and $p_{j_1}=1$,  applying rotations
        with the center at $0$. Differential forms 
        $d(\phi_{i_\a,j_\a}-\phi_{i_1,j_1})$ on $C_n'$ are pullbacks of
         differential forms $d\phi_{i_\a,j_\a}$ on $C_n''$. 
        The integral of a product of $2n-3$ closed $1$-forms
         $d\phi_{i_\a,j_\a},\,\,\a=1,\dots ,2n-3$ over $C'_n$ is equal 
         to $\pm 2\pi$ times
          the integral of the product
          $2n-4$ closed $1$-forms
           $d\phi_{i_\a,j_\a},\,\,\a=2,\dots ,2n-3$ over $C''_n$.
       
       The space $C''_n$ is a complex 
       manifold. We are calculating 
       the absolutely
       converging integral of the type
       $$\int\limits_{C''_n} \prod_{\alpha} d Arg(Z_{\alpha})$$
       where $Z_{\alpha}$ are holomorphic invertible functions on $C''_n$
        (differences between complex coordinates of points of the 
        configuration). We claim that it is zero, because of the 
        general result proven in 6.6.1. \qed

 \vskip 0.5truecm
\par\noindent {\bf 6.6.1. A trick using logarithms}
\vskip 0.5truecm

        \proclaim Theorem. Let $X$ be a complex algebraic variety
         of dimension $N\ge 1$, and $Z_1,\dots,Z_{2N}$ be rational 
         functions on $X$, not equal identically to zero.
          Let $U$ be any Zariski open subset of $X$ such that 
          functions $Z_\a$ are defined and non-vanishing
             on $U$, and
           $U$ consists of smooth points.
          Then the integral 
          $$\int\limits_{U(\C)} \wedge_{\a=1}^{2N} d(Arg\,Z_\a)$$
           is absolutely convergent, and equal to zero.  \par
           
           This result seems to be new,
           although the main trick used in the proof is well-known.
            A.~Goncharov told me that he also came to the same result
             in his study of mixed Tate motives.
           
       {\bf Proof:} First of all, we claim that
         the differential
        form $\wedge_{\a=1}^{2N} d Arg(Z_\a)$ on $U(\C)$ 
         coincides with the form $\wedge_{\a=1}^{2N}
         d Log\,|Z_\a|$ (this is the trick). 
         
  We can replace $d Arg(Z_\a)$
         by the difference of a holomorphic an anti-holomorphic form
  $${1\over  2 i}\left(d(Log\,Z_{\a})-d(Log\,{\overline Z}_\a)
  \right)\,\,\,.$$
   Thus, the form which we integrate over $U(\C)$ is a sum of products of 
    holomorphic and of anti-holomorphic forms.
     The summand corresponding to a product of a non-equal number
      of holomorphic and of anti-holomorphic forms, vanishes
       identically because $U(\C)$ is a complex manifold.
       The conclusion is that the number of anti-holomorphic
        factors in non-vanishing summands is the same for all of them,
         it coincides
        with the complex dimension $N$ of $U(\C)$. The same products of
         holomorphic and of anti-holomorphic forms survive
         in the product 
         $$\bigwedge_{\a=1}^{2N} d\,Log\,{|Z_\a|}=
         \bigwedge_{\a=1}^{2N} {1\over  2}
         \left(d(Log\,Z_{\a})+d(Log\,{\overline Z}_\a)
         \right)\,\,\,.$$
         
         Let us choose a compactification ${\overline U}$ of $U$ such
          that ${\overline U}\setminus U$ is a divisor with normal crossings.
          If $\phi$ is a smooth differential form on $U(\C)$ such that
           coefficients of $\phi$ are locally integrable on ${\overline U}(\C)$,
            then we denote by ${\cal I} (\phi)$ corresponding differential
             form on ${\overline U}(\C)$ with coefficients in the space of
              generalized
              functions.
              
            \proclaim Lemma. Let $\omega$ be a form on $U(\C)$ which is a
             linear combination of products of functions $Log\,|Z_\a|$
              and of $1$-forms $d\,Log\,|Z_\a|$ where $Z_\a\in
              {\cal O}^{\times}(U) $ are
               regular invertible functions on $U$.
                Then coefficients of
                $\omega$ and of $d\omega$
                are locally $L^1$ functions on ${\overline U}(\C)$.
                 Moreover,
                  ${\cal I} (d \omega)=d({\cal I}(\omega))$. Also,
                   the integral
                    $\int\limits_{U(\C)}\omega$ is absolutely convergent
                     and equal to the integral $
                     \int\limits_{{\overline U}(\C)}{\cal I}(\omega)$.
                   \par
                   
                The lemma is an elementary exercise
                 in generalized functions, after  passing
                 to local coordinates on 
                ${\overline U}(\C)$. We leave  details of the proof to 
                the reader. Also, the statement of the lemma
                 remains true without
                 the condition that 
                ${\overline U}\setminus U$ is a divisor with normal crossings.
                 \qed 
                
                The vanishing of the integral in the theorem  
                is clear now by the Stokes formula:
                
               $$\int\limits_{U(\C)}\bigwedge_{\a=1}^{2N}
                d\,Arg\,(Z_\a)
                =\int\limits_{U(\C)}\bigwedge_{\a=1}^{2N}
                d\,Log\,|Z_\a|=\int\limits_{{\overline U}(\C)} {\cal I}\left(
                d\left(Log\,|Z_1|\,
                \bigwedge_{\a=2}^{2N}
                d\,Log\,|Z_\a|\right)\right)=$$
                $$=\int\limits_{{\overline U}(\C)}
                 d\left({\cal I}\left(
               Log\,|Z_1|\,
                \bigwedge_{\a=2}^{2N}
                d\,Log\,|Z_\a|\right)\right)=0\,\,\,.\,\,\,\qed$$

    \vskip 0.5truecm
\par\noindent {\bf 6.6.2. Remark}
\vskip 0.5truecm     

The vanishing of the integral in the lemma from 6.6 has 
 higher-dimensional analogue  which is crucial in the perturbative
  Chern-Simons theory in the dimension $3$, and its generalizations
   to dimensions $\ge 4$ (see [Ko1]). 
   However, the vanishing of integrals in dimensions $
   \ge 3$ follows from a much simpler fact which is
    the existence of a geometric
    involution
    making the integral equal to minus itself. In the present paper
     we will use many times
     such kind of arguments involving involutions.

    \vskip 0.5truecm
\par\noindent {\bf 7. Formality conjecture for general manifolds}
\vskip 0.5truecm

        In this section we establish the formality conjecture for general
         manifolds, not only for open domains
         in $\R^d$. It turns out that 
         that essentially all work is already done. The only new 
         analytic result is  vanishing of certain integrals over  
          configuration spaces, analogous to the lemma from 
          6.6.

         One can treat $\R^d_{formal}$,
             the formal completion of vector space 
           $\R^d$ at zero,  
              in many respects as  usual manifold.
              In particular, we can define differential graded
               Lie algebras $D_{poly}(\R^d_{formal})$ and $T_{poly}
               (\R^d_{formal})$.
                The Lie algebra  
                $W_d:=Vect(\R^d_{formal})$ is the standard
                Lie algebra of  formal 
                vector fields. We consider
                 $W_d$ as a differential graded Lie algebra
                  (with the trivial grading and the differential equal to 0).
                 There are natural homomorphisms
                  of differential graded Lie algebras:
                  $$m_T: W_d\ra T_{poly}
               (\R^d_{formal}),\,\,\,\,\,m_D:W_d\ra D_{poly}
               (\R^d_{formal})\,,$$
            because  vector fields can be considered
            as polyvector fields and as differential operators.

          We will use the following properties of the quasi-isomorphism
           $\U$ from 6.4:
           
           P1) $\U$ can be defined for $\R^d_{formal}$ as well,
           
           P2) for any $\xi\in W_d$ we have the equality
            $$\U_1(m_T(\xi))=m_D(\U_1(\xi))\,,$$

           P3) $\U$ is $GL(d,\R)$-equivariant,

           P4) for any $k\ge 2$, $\xi_1,\dots,\xi_k\in W_d$ 
           we have the equality 
            $$\U_k(m_T(\xi_1)\otimes\dots\otimes
            m_T(\xi_k))=0$$,
           
           P5) for any $k\ge 2$, $\xi\in gl(d,\R)\subset W_d$, and for 
           any $\eta_2,\dots,\eta_k\in T_{poly}
               (\R^d_{formal})$ we have
                $$\U_k(m_T(\xi)\otimes\eta_2
                \otimes\dots\otimes\eta_k)=0\,\,\,.$$
           
          We will construct quasi-isomorphisms from $T_{poly}(X)$
           to $D_{poly}(X)$ for arbitrary $d$-dimensional manifold $X$
            using only properties P1-P4 of the map $\U$.
             Properties P1,P2 and P3 are evident, and the
             properties P4,P5 will be
             established later (subsections 7.3.1.1 and 7.3.3.1).
           
            It will be convenient to
            use in this section the geometric language
           of formal graded manifolds, instead of the algebraic language
            of $\L$-algebras.
           Let us fix the dimension $d\in \N$. We introduce three
            formal graded $Q$-manifolds {\it without} base points:
            $$\T,\D,\W\,\,\,.$$ 
            These formal graded $Q$-manifolds are obtained in the usual way
            from  
            differential graded Lie algebras $T_{poly}
               (R^d_{formal})$, $D_{poly}
               (R^d_{formal})$ and $W_d$ forgetting base points.
           
                In next two subsections (7.1 and 7.2) we present
                two general geometric constructions, which  will used
                 in 7.3 for the proof of formality of $D_{poly}(X)$.
                 
        \vskip 0.5truecm
\par\noindent {\bf 7.1. Formal geometry (in the sense of I.~Gelfand 
and D.~Kazhdan) }
\vskip 0.5truecm             
                 
        Let $X$ be a smooth manifold of dimension $d$.
         We associate with $X$ two infinite-dimensional manifolds,
         $X^{coor}$ and $X^{aff}$. The manifold $X^{coor}$
          consists of pairs $(x,f)$ where $x$ is a point of $X$
          and $f$ is an infinite germ of a coordinate system
           on $X$ at $x$,
           $$f:(\R^d_{formal},0)\hookrightarrow (X,x)\,\,\,.$$
        We consider $X^{coor}$ as a projective limit
         of finite-dimensional manifolds (spaces of finite germs
          of coordinate systems). There is an action 
            on $X^{coor}$ 
           of the (pro-Lie) group $G_d$
           of formal diffeomorphisms of $\R^d$ preserving base point
           $0$. 
           The natural projection 
            map $X^{coor}\ra X$ is a principal $G_d$-bundle.
            
             The manifold $X^{aff}$ is defined as the quotient space 
              $X^{coor}/GL(d,\R)$. It can be thought as the
               space of formal affine structures at points of $X$.
               The main reason to introduce
                $X^{aff}$ is that  fibers of the natural projection
               map $X^{aff}\ra X$ are contractible.

             The Lie algebra of the group $G_d$
              is a subalgebra of codimension $d$ in $W_d$. It consists
               of formal vector fields vanishing at zero. Thus,
                $Lie(G_d)$ acts on $X^{coor}$. It is easy to see
                 that in fact the whole Lie algebra $W_d$
                  acts on $X^{coor}$ and is isomorphic
                   to  the tangent space to $X^{coor}$ at each point.
                    Formally, the 
            infinite-dimensional manifold $X^{coor}$ looks as 
             a principal homogeneous space of the non-existent
              group with the Lie algebra $W_d$.
              
               The main idea of formal geometry (se [GK]) is to
                replace $d$-dimensional manifolds
                by ``principal homogeneous spaces'' of $W_d$.
            Differential-geometric constructions on $X^{coor}$ can be
             obtained from Lie-algebraic constructions for $W_d$.
             For a while we will work only with $X^{coor}$, and then at the end
             return to $X^{aff}$. In terms of Lie algebras it corresponds
              to the difference between absolute and relative cohomology.

        \vskip 0.5truecm
\par\noindent {\bf 7.2. Flat connections and $Q$-equivariant maps}
\vskip 0.5truecm        

Let $M$ be a $C^{\infty}$-manifold (or a complex analytic manifold, or
 an algebraic manifold, or a projective limit of manifolds,...). 
Denote by $\Pi T M$ the supermanifold which is the
total space of the tangent bundle of $M$ endowed with the reversed parity.
              Functions on the $\Pi T M$ are differential forms
              on $M$. The de Rham differential $d_M$ on forms can be considered
              as an odd vector field on $\Pi T M$ with the square equal to $0$.
            Thus, $\Pi TM$ is a $Q$-manifold. 
             It seems that the accurate notation for $\Pi TM$ considered
              as a graded manifold 
               should be $T[1] M$ 
               (the total space of the graded vector bundle
                $T_M[1]$ considered as a graded
                 manifold).

           Let $N\ra M$ be a bundle over a manifold $M$ whose fibers
               are manifolds, or vector spaces, etc., endowed with a flat
               connection $\nabla$.  Denote by $E$ the pullback
              of this bundle to $B:=\Pi TM$. The connection $\nabla$
               gives a lift of the vector field $Q_B:=d_M$ on $B$
               to
                the vector field $Q_E$ on $E$.
                 This can be done for arbitrary connection, and only for flat 
                 connection the identity $[Q_E,Q_E]=0$ holds. 
             
            A generalization of a (non-linear) bundle with a flat
             connection is a  $Q$-equivariant bundle whose
              total space and the base are  $Q$-manifolds.
               In the case of graded vector bundles over $T[1] M$
                this notion was introduced
           Quillen under the name of a  superconnection (see [Q]). 
               A generalization of the notion of a covariantly flat morphism
     from one bundle to another is the notion of a $Q$-equivariant map.
       
       \proclaim Definition. A  flat family over $Q$-manifold $B$ 
     is a pair $(p:E\ra B,\sigma)$ where
        $p:E\ra B$ is a
       $Q$-equivariant
        bundle whose fibers are formal  manifolds, and a $\sigma
        :B\ra E$ is a 
        $Q$-equivariant 
         section of this bundle.
         \par

          In the case $B=\{point\}$ a flat family over $B$ is the same
          a formal  $Q$-manifold with  base point. 
          It is clear that flat families 
          over a given 
        $Q$-manifold form a category.
        
        We  apologize for the terminology. More precise name
         for ``flat families'' would be ``flat families of pointed formal
         manifolds'', but it is too long.
        
        One can define 
         analogously  flat graded families over graded $Q$-manifolds.

      We refer the reader to a discussion of further examples of 
      $Q$-manifolds    in [Ko3].
      
      \vskip 0.5truecm
\par\noindent {\bf 7.3.  Flat families in deformation quantization}
\vskip 0.5truecm    
      
      Let us return to our concrete situation. We construct in this section
      two
       flat families over $\Pi TX$ (where $X$ is a $d$-dimensional
        manifold),
        and a morphism between them. This  will be done
         in several steps.
      
   \vskip 0.5truecm
\par\noindent {\bf 7.3.1.  Flat families over $\W$}
\vskip 0.5truecm

         The first bundle over $\W$ is trivial as a $Q$-equivariant
         bundle,
            $$\T\times\W\ra \W$$
            but with a non-trivial section ${\sigma}_\T$. This section is not 
            the zero section, but the graph of the $Q$-equivariant map
             $\W\ra \T$ coming from the  homomorphism
             of differential graded Lie algebras $m_T:W_d\ra 
             T_{poly}(\R^d_{formal})$.
            Analogously, the second bundle is the trivial $Q$-equivariant
            bundle 
            $$\D\times\W\ra \W$$
           with the section $\s_\D$ coming from the homomorphism
           $m_D:W_d\ra D_{poly}(\R^d_{formal})$.

          Formulas from  6.4  
          give    a $Q$-equivariant map
           $\U:\T\ra \D$.

             \proclaim Lemma. The morphism $(\U\times id_{\W}):\T\times\W\ra
             \D \times \W$ is a morphism of flat families over $\W$.
             \par
             
             {\bf Proof:}
             We have to check that $(\U\times id_{\W})$ maps one section 
              to another, i.e. that
               $$(\U\times id_{\W})\circ \s_\T=
              \s_D\in Maps(\W,\D\times \W)\,\,\,.$$

                                We compare Taylor coefficients.
                                 The linear
              part $\U_1$
              of $\U$ maps a vector field (considered as a polyvector
               field) to itself, considered as a differential operator
                (property P2).
                 Components $\U_k(\xi_1,\dots
                 , \xi_k)$ for $k\ge 2,\,\,\,\xi_i\in T^0(\R^d)=
                 \G(\R^d,T)$ 
                 vanish, which is the property P4. \qed
                 
                \vskip 0.5truecm
\par\noindent {\bf 7.3.1.1.  Proof of the property P4}
\vskip 0.5truecm 

Graphs appearing in the calculation of $\U_k(\xi_1,\dots
                 , \xi_k)$ have $k$ edges, $k$ vertices of the first type,
                  and $m$ vertices of the second type, where
                  $$2k+m-2=k\,\,\,.$$
                  Thus, there are no such graphs for $k\ge3$ as $m$
                  is non-negative. The only intersting case is $k=2,m=0$.
                 The graph is looking as
                       
           \vskip 1cm     \centerline{
\epsfbox{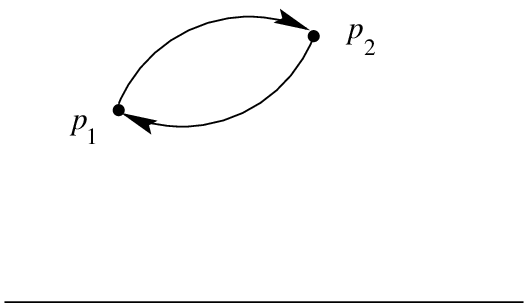}
} \vskip 11mm
                 
                  By our construction, $\U_2$ restricted to vector fields
                  is equal
                to the non-trivial quadratic map 
                $$\xi\longmapsto \sum_{i,j=1}^d 
                \p_i(\xi^j)\p_j(\xi^i)\in\G(\R^d,
                {\cal O}),\,\,\,\,
                \xi=\sum_i \xi^i\p_i\in\G(\R^d,T) $$
                with the weight
                $$\int\limits_{C_{2,0}} d\phi_{(12)} d\phi_{(21)}=\int
                \limits_{\H\setminus
                 \{z_0\}} d\phi(z,z_0)\wedge d\phi(z_0,z)$$
                 where $z_0$ is an arbitrary point of $\H$.
                 
                 \proclaim Lemma. For arbitrary angle map the integral
                  $\int\limits_{\H\setminus
                 \{z_0\}} d\phi(z,z_0)\wedge d\phi(z_0,z)$ is equal to zero.
                 \par
                 
                 {\bf Proof:}
                  We have a map $\OC_{2,0}\ra S^1\times S^1,\,\,
                 \, [(x,y)]\mapsto (\phi(x,y),\phi(y,x))$. We calculate
                 the  integral of the 
                 pullback of the standard volume element on two-dimensional
                  torus.
                 It is easy to see that the integral does not depend on the 
                 choice of map $\phi:\OC_{2,0}\ra S^1$. The reason
                  is that the image of the boundary of the integration domain
                   $\p \OC_{2,0}$ 
                   in $S^1\times S^1$ cancels with the reflected
                    copy of itself under the involution $(\phi_1,\phi_2)
                    \mapsto (\phi_2,\phi_1)$ of the torus $S^1\times S^1$.  
                  Let us assume that $\phi=\phi^h$ and $z_0=0+1\cdot i$.
                   The integral vanishes because the involution
                    $z\mapsto -{\overline z}$ reverses the orientation
                     of $\H$ and preserves the  form 
                     $d\phi(z,z_0)\wedge d\phi(z_0,z)$. \qed
   
  \vskip 0.5truecm
\par\noindent {\bf 7.3.2.  Flat families over $\Pi T(X^{coor})$}
\vskip 0.5truecm                   

If $X$ is a $d$-dimensional manifold, then there is a natural map
 of  $Q$-manifolds (the Maurer-Cartan form)
   $$\Pi T(X^{coor})\ra \W\,\,\,.$$
   It follows from  following general reasons.
   If $G$ is a Lie group, then it acts freely
    by left translations on itself, and also on $\Pi TG$.
     The quotient $Q$-manifold $\Pi TG/G$ 
     is equal to $\Pi \g$ where $\g=Lie(G)$. 
     Thus, we have a $Q$-equivariant map
       $$\Pi TG\ra \Pi \g\,\,\,.$$
      Analogous construction works for any principal
       homogeneous space over $G$. We apply it
        to $X^{coor}$ considered as a principal homogeneous
         space for a non-existent group with the Lie algebra $W_d$.
         
         The pullbacks of flat families of formal manifolds
          over $\W$ constructed in 7.3.1,
           are  two flat families over $\Pi T ( X^{coor})$.
           As $Q$-equivariant bundles
             these families are trivial bundles
             $$\T\times \Pi T (X^{coor})\ra \Pi T (X^{coor}),\,\,\,\,
             \D\times \Pi T (X^{coor})\ra \Pi T (X^{coor})\,\,\,.$$
             Pullbacks of sections $\sigma_\T$ and $\sigma_\D$
              gives sections in the bundles above. These
               sections we denote again by $\sigma_\T$ and $\sigma_\D$. 
            The pullback of the morphism $\U\times id_\W$ is also a morphism
            of flat families. 
              
   \vskip 0.5truecm
\par\noindent {\bf 7.3.3.  Flat families over $\Pi T (X^{aff})$}
\vskip 0.5truecm

       Recall  that $X^{aff}$ is the quotient space of $X^{coor}$
       by the action of $GL(d,\R)$. Thus, from functorial
        properties of operation $\Pi T$ ($=\underline{Maps}(\R^{0|1},\cdot)$)
         follows
         that $\Pi T (X^{aff})$ is the quotient
             of $Q$-manifold $\Pi T (X^{coor})$ by the action of $Q$-group
              $\Pi T (GL(d,\R))$. We will construct an
               action of $\Pi T (GL(d,\R))$ on flat families 
               $\T\times \Pi T (X^{coor})$ and $\D\times \Pi T (X^{coor})$
                over $\Pi T (X^{coor})$. We claim that 
                the morphism between these
                 families is invariant under the action of $\Pi T (GL(d,\R))$.
                Flat families over $\Pi T (X^{aff})$ will be defined
                 as quotient families. The morphism between them will
                  be the quotient morphism.

                The action of $\Pi T (GL(d,\R))$ on $\T$ and on $\W$
                 is defined as follows. First of all, if $G$ is a Lie group
                 with the Lie algebra $\g$, then $\Pi T G$ acts
                 $Q$-equivariantly 
                  on $Q$-manifold $\Pi \g$, via the identification
                   $\Pi \g=\Pi TG/G$. Analogously, 
                    if $\g$ is a subalgebra of a larger Lie algebra
                     $\g_1$, and an action of $G$ on $\g_1$ is given in a way 
                     compatible with the inclusion $\g\hra\g_1$, then
                      $\Pi T G$
                    acts on $\Pi \g_1$. We apply this
                     construction to the case $G=GL(n,\R)$ and
                      $\g_1=T_{poly}(\R^d_{formal})$ or
                       $\g_1=D_{poly}(\R^d_{formal})$.

              One can check easily that sections $\sigma_\T$ and $\sigma_\D$
               over $\Pi T (X^{coor})$ 
               are $\Pi T (GL(d,\R))$-equivariant.
               Thus, we get two flat families over $\Pi T( X^{aff})$.
               
               The last thing we have to check is that the morphism $\U\times 
               id_{\Pi T (X^{coor})}$ 
                of flat families 
                $$
            \T \times  \Pi T (X^{coor})\ra
             \D\times \Pi T (X^{coor})$$
               is $\Pi T (GL(d,\R))$-equivariant.
             After the translation
            of the problem to the language of Lie algebras, we
             see that we should check that $\U$
              is $GL(d,\R)$-invariant
                 (property P3,
                 that is clear by our construction),
                  and that if we substitute an element of $gl(d,\R)\subset
                   W_d$ in $\U_{\ge 2}$, we get zero (property P5, see
                    7.3.3.1).

     \proclaim Conclusion. We constructed two flat families over $
              \Pi T (X^{aff})$ and a morphism between them. Fibers of
              these families are isomorphic to $\T$ and to $\D$. \par               
                    
                   \vskip 0.5truecm
\par\noindent {\bf 7.3.3.1. Property P5}
\vskip 0.5truecm

               This is again reduces to the calculation of an integral.
                Let $v$ be a vertex 
                of $\G$ to which we put element of $gl(d,\R)$. There
                 is exactly one edge starting at $v$ because we
                 put a vector field here.
            If there are no edges ending at $v$, 
            then the integral is zero because
                  the domain of integration is foliated by lines
                  along which all forms vanish. These
                    lines are 
                    level sets of the function $\phi(z,w)$ where 
                    $w\in \H\sqcup
                    \R$
                    is fixed and $z$ is the point on $\H$ corresponding to
                     $v$.

           \vskip 1cm       \centerline{
\epsfbox{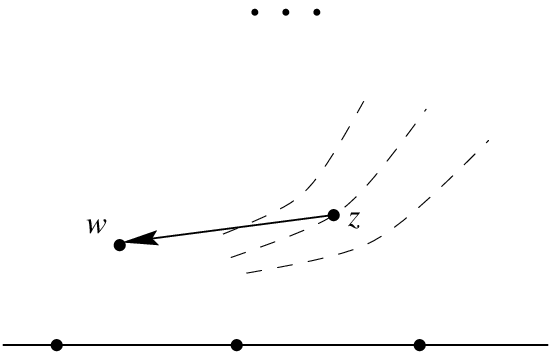}
} \vskip 11mm

 If there are at least $2$ edges ending at $v$, then the corresponding
  polydifferential operator is equal to zero, because second
   derivatives of coefficients
   of a linear vector field vanish.
                   
                   The only relevant case is when there is only one edge
                    starting at $v$, and only one edge ending
                    there. If these two edges connect our vertex with
                    the same vertex of $\G$, then the vanishing follows
                    from the lemma in the section 7.3.2.1. 
                    If our vertex is connected
                     with two different vertices,

                 \vskip 1cm    \centerline{
\epsfbox{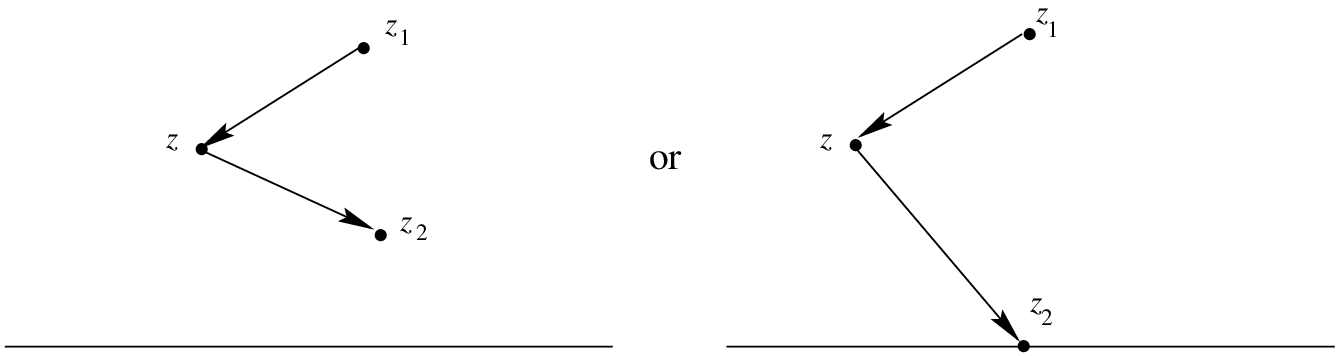}
} \vskip 11mm \hfill\break
    then we apply the following two
                     lemmas:
                       
             \proclaim Lemma. Let $z_1\ne z_2\in\H$ be two distinct points
             on $\H$. Then the integral
              $$\int\limits_{z\in \H\setminus\{z_1,z_2\}} d\phi(z_1,z)
              \wedge d\phi(z,z_2)$$
              vanishes. \par
              
              \proclaim Lemma. Let $z_1\in \H,\, z_2\in\R$ be 
              two  points
             on $\H\sqcup \R$. Then the integral
              $$\int\limits_{z\in \H\setminus\{z_1,z_2\}} d\phi(z_1,z)
              \wedge d\phi(z,z_2)$$
              vanishes. \par
              
            {\bf Proof:}
            One can prove analogously to the lemma in 7.3.1.1 
             that the integral does not depend 
             on the choice of an angle map, and also on points $z_1,\,z_2$.
            In the case of $\phi=\phi^h$ and both points $z_1,z_2$ are pure
             imaginary, the vanishing follows from the anti-symmetry 
             of the integral under the involution $z\mapsto -{\overline z}$.
              \qed

                   \vskip 0.5truecm
\par\noindent {\bf 7.3.4.  Flat families over $X$}
\vskip 0.5truecm 

Let us choose a section $s^{aff}$ of the bundle
           $X^{aff}\ra X$. Such section always exists because
            fibers of this bundle are contractible. For example, any
          torsion-free   connection $\nabla$
          on the tangent bundle to $X$
           gives a section $X\ra X^{aff}$. Namely, the exponential map
            for $\nabla$ gives an identification of a neighborhood of
             each point $x\in X$ with a neighborhood of zero in the vector 
             space $T_x X$, i.e. an affine structure on $X$ near $x$, and
              a point of $X^{aff}$ over $x\in X$.
           
             The section $s^{aff}$ defines a map of formal graded
             $Q$-manifolds $\Pi TX\ra \Pi T (X^{aff})$.
             After taking the pullback
               we get  two flat families $\T_{s^{aff}}$ and
                $\D_{s^{aff}}$ 
                over $\Pi TX$ and an morphism $m_{s^{aff}}$
                from one to another.
                
               We claim that these two flat families admit definitions
                independent of $s^{aff}$. Only the morphism  $m_{s^{aff}}$ 
                 depends on $s^{aff}$. 
                   
                   Namely, let us consider 
                   infinite-dimensional bundles of differential
                    graded Lie algebras
                   $jets_{\infty} T_{poly}$
                    and  $jets_{\infty} D_{poly}$
                    over $X$ whose fibers
                   at $x\in X$ are spaces of infinite jets
                   of  polyvector fields or 
                   polydifferential operators at $x$
                   respectively. These two bundles carry natural
                    flat connections (in the usual sense, not as in 7.2)
                     as any bundle of infinite jets.
                    Thus, we have two flat families (in generalized sense) 
                     over $\Pi TX$.
                   
            \proclaim Lemma.  Flat families $\T_{s^{aff}}$ and
                $\D_{s^{aff}}$ are canonically isomorphic to flat families 
                described just above. \par
                
        {\bf Proof:} it follows from definitions that  pullbacks
   of bundles  $jets_{\infty} T_{poly}$
                    and  $jets_{\infty} D_{poly}$  from $X$
                         to $X^{coor}$ 
                    are canonically
                     trivialized. The Maurer-Cartan $1$-forms on $X^{coor}$ 
                      with values in graded Lie algebras 
                      $T_{poly}(\R^d_{formal})$ or
                       $D_{poly}(\R^d_{formal})$
                       come from  pullbacks of flat connections
                        on bundles of infinite jets. Thus, we identified
                       our  flat families over $\Pi T(X^{coor})$
                        with pullbacks. The same is true for 
                        $X^{aff}$. \qed
              
     \vskip 0.5truecm
\par\noindent {\bf 7.3.5. Passing to global sections}
\vskip 0.5truecm

          If in general $(p:E\ra B,\sigma)$ is a 
           flat family, then one can make a new 
             formal pointed $Q$-manifold:
             $$\left( \Gamma(E\ra B)_{formal},\sigma\right)\,\,\,.$$
             This is an infinite-dimensional formal
       super manifold, the formal completion of the space of sections   
        of the bundle $E\ra B$ at the point $\sigma$.
          The structure of $Q$-manifold on $\Gamma(E\ra B)$
          is evident because the 
          Lie supergroup $\R^{0|1}$ acts on $E\ra B$.

                \proclaim Lemma. Formally completed spaces 
                of global sections of flat families
                 $\T_{s^{aff}}$ and
                $\D_{s^{aff}}$ 
                a naturally quasi-isomorphic to $T_{poly}(X)$ and
                 $D_{poly}(X)$ respectively.
                \par
                
               {\bf Proof:}
                It is well-known that if $E\ra X$ is a vector bundle
                then de Rham cohomology of $X$ with coefficients in  formally
                 flat infinite-dimensional bundle $jets_{\infty} E$
                  are concentrated in degree $0$ and canonically isomorphic
                  to the vector space $\G(X,E)$. Moreover, the
                   natural homomorphism of complexes
                   $$\bigl(\G(X,E)[0],{\rm\,differential\,}=0\bigr)
                   \ra \bigl(\Omega^*(X,jets_{\infty}(E))
                   ,{\rm \, de\,\,\,Rham
                   \,\,\,differential\,}\bigr)$$
                   is quasi-isomorphism.
                   
                   Using this fact, the lemma from the previous subsection, 
                    and appropriate filtrations (for spectral sequences)
                     one sees that 
                that the natural $Q$-equivariant
               map from the formal
                $Q$-manifold $( T_{poly}(X)_{formal}[1],0) $ to 
                $(\G( \T_{s^{aff}}\ra T[1]X)_{formal},\sigma_\T)$
                 (and analogous map for 
                $D_{poly}$)
 is a  quasi-isomorphism. \qed

          	It follows from the lemma above and the result of 
                    4.6.1.1 that we have a chain of 
          	quasi-isomorphisms
          	$$T_{poly}(X)[1]_{formal} \ra
          	\G( \T_{s^{aff}}\ra T[1]X)_{formal} \ra 
          	\G( \D_{s^{aff}}\ra T[1]X)_{formal} \longleftarrow 
          	T_{poly}(X)[1]_{formal}\,\,\,.$$
          	
          	Thus, differential graded Lie algebras $T_{poly}(X)$
          	 and $D_{poly}(X)$ are quasi-isomorphic.
          	  The theorem from 4.6.2. is proven. \qed 
          	 
          	 The  space of sections of the bundle $X^{aff}\ra X$
          	  is contractible. From this fact one can
          	   conclude that the quasi-isomorphism constructed above
          	    is well-defined homotopically. 
         
         \vskip 0.5truecm
\par\noindent {\bf 8. Cup-products}
\vskip 0.5truecm

\par\noindent {\bf 8.1. Cup-products on tangent cohomology}
\vskip 0.5truecm    

Differential graded Lie algebras $T_{poly}$, $D_{poly}$ and  (more
 generally)
 shifted by $[1]$ Hochschild complexes of arbitrary associative algebras,
 all carry an
 additional structure.  We do not know at the moment a  definition,
  it should be something close to so called homotopy Gerstenhaber algebras
   (see [GV], [GJ]), although definitely not precisely this.
 At least, a visible part of this structure 	is a commutative
  associative product of degree $+2$ on cohomology of
  the tangent space to any solution of the
  Maurer-Cartan equation. Namely, if $\g$ is one of differential
   graded Lie algebras listed above and $\ga \in (\g\otimes \m)^1$
    satisfies $d\ga+{1\over 2}[\ga,\ga]=0$ where $\m$ is a finite-dimensional
     nilpotent non-unital differential
      graded commutative associative algebra,
      the tangent space $T_{\ga}$ is defined as 
      complex $\ga\otimes\m[1]$ endowed with the differential
       $d+[\ga,\cdot]$. Cohomology space $H_{\ga}$
       of this differential  
        is a graded module over graded algebra $H(\m)$
         (the cohomology space of $\m$ as a complex). If $\ga_1$ and $\ga_2$ 
        are two gauge equivalent solutions, then $H_{\ga_1}$ and
         $H_{\ga_2}$ are (non-canonically) equivalent $\m$-modules.

        We define now  cup-products for all three differential graded Lie 
        algebras listed at the beginning of this section.
         For $T_{poly}(X)$ the cup-product is defined as the usual
          cup-product of polyvector fields (see 4.6.1). One can check
           directly that this cup-product is compatible with the
            differential  $d+[\ga,\cdot]$, and is
            a graded commutative associative product.
             For the Hochschild complex of an associative algebra
              $A$ the cup-product on $H_{\ga}$
              is defined in a more 
             tricky way. It is defined on the complex by  the formula
             $$(t_1\cup t_2)(a_0\otimes\dots\otimes a_n):=$$
            $$ \sum_{0\le k_1\le k_2\le k_3\le k_4\le n}
             \pm \ga^{n-(k_2-k_1+k_4-k_3)}(a_0\otimes\dots
             \otimes t_1(a_{k_1}\otimes\dots)\otimes a_{k_2}\otimes
             \dots\otimes t_2(a_{k_3}\otimes\dots)\otimes a_{k_4}
             \otimes\dots)$$
             where $\ga^{l}\in (\k[0]
             \cdot 1\oplus\m)^{1-l}\otimes Hom(A^{\otimes (l+1)},A)$
              is homogeneous component of $(\ga+1\otimes m_A)$.
              
              It is not a trivial check that the cup-product on the Hochschild 
              complex is 
             compatible with differentials, and also is
              commutative, associative
              and gauge-equivariant
              on the level of cohomology. Formally, we will not use this fact.
              The proof is a direct calculation with Hochschild cochains.
               Even if one replaces formulas by appropriate pictures
                the calculation is still quite long, about 4-5 pages
                 of tiny drawings. Alternatively, there is a simple 
                 abstract explanation using the interpretation
                  of the deformation theory related with the shifted Hochschild
                  complex as a deformation theory of triangulated categories
                  (or, better, $A_{\infty}$-categories, see [Ko4]).
                  We will discuss it in more details in the sequel
                   to the present paper.

              We define the cup-product for $D_{poly}(X)$ by the restriction
               of formulas for the cup-product in $C^{\bullet}(A,A)$.

              \vskip 0.5truecm
\par\noindent {\bf 8.2. Compatibility of $\,\U$ with cup-products }
\vskip 0.5truecm  
              
              \proclaim Theorem. The quasi-isomorphism $\U$ constructed
               in section 6 maps the cup-product for $T_{poly}(X)$
                to the cup-product for $D_{poly}(X)$.
                \par
                
              {\bf Sketch of the proof:}
                we translate  the statement of the theorem
                 to the language of graphs and integrals.
                 The tangent map is given by integrals where one of vertices
                  of the first type is marked. This is the vertex where we put 
               a representative $t$ for
               the    tangent element $[t]\in H_{\gamma}$. We put 
       copies of  $\ga$  (which is   a  polyvector field with values in
       $\m$) into all other vertices of the first type. The rule which we 
       just described follows directly from the Leibniz formula
        applied to the Taylor series for $\U$.

               Now we are interested in the behavior of the tangent  map
               with respect to a bilinear operation on the tangent space.
                It means that we have now {\it two} marked
                vertices of the first type.
         
 \vskip 0.5truecm
\par\noindent {\bf 8.2.1. Pictures for the cup-product in polyvector fields}
\vskip 0.5truecm                  

 We claim that
                  the cup-product for the case $T_{poly}(X)$
                   corresponds to pictures where two
                    points (say, $p_1,p_2$)
                     where we put representatives of
                     elements of $H_{\ga}$ which we want
                     to multiply, are infinitely close points on $\H$.
                     Precisely, it means that we integrate over
                    preimages $P_{\alpha}$
                    of some point $\alpha$ in $\R/2\pi \Z
                        \simeq C_2\subset \OC_{2,0}$
                    with respect to the forgetting map
                    $$\OC_{n,m}\ra\OC_{2,0}\,\,\,.$$
                    It is easy to see that $P_{\alpha}$ has codimension $2$
                     in $\OC_{n,m}$ and contains no strata $C_T$ of
                      codimension $2$. It implies that as a singular chain
                    $P_{\alpha}$ is equal to the sum  of closures
                    of non-compact hypersurfaces
                     $$P_{\alpha}\cap \p_S(\OC_{n,m}),\,\,
                     P_{\alpha}\cap \p_{S_1,S_2}(\OC_{n,m})$$
                      in boundary
                      strata of $\OC_{n,m}$. It is easy to see that
                       intersections $P_{\alpha}\cap \p_{S_1,S_2}(\OC_{n,m})$
                        are empty, and intersection 
                        $P_{\alpha}\cap \p_S(\OC_{n,m})$ is non-empty
                         iff $S\supseteq \{1,2\}$.
                         The picture is something like
                         
               \vskip 1cm      \centerline{
\epsfbox{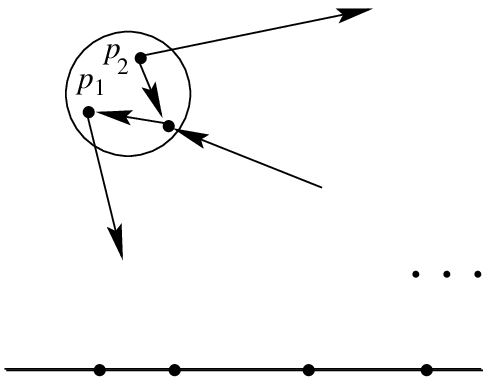}
}   \vskip 11mm
 
     Points $p_1$ and $p_2$ should not be connected by an edge
     because otherwise the integral vanishes, there is no directions
     over which we can integrate form $d\phi(p_1,p_2)$. Also, if $\#S\ge 3$
      then the integral vanishes by lemma from 6.6.
       The only non-trivial case which is left is when $S=\{1,2\}$ and
        points $p_1,\,p_2$ are not connected:
        
  \vskip 1cm \centerline{
\epsfbox{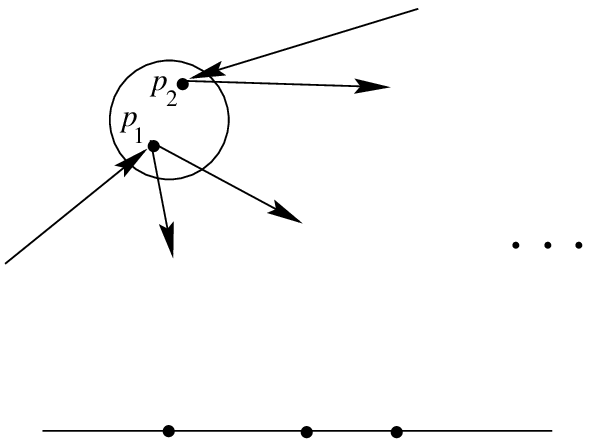}
}  \vskip 11mm                
  
This figure exactly corresponds to the cup-product in  $T_{poly}(X)$.

 \vskip 0.5truecm
\par\noindent {\bf 8.2.2. Pictures for the cup-product in 
the Hochschild complex}
\vskip 0.5truecm    
                     
                     The cup-product for $D_{poly}(X)$ is given by pictures
                     where these two points are separated and infinitely close
                      to $\R$. Again, the precise definition
                       is that we integrate of the preimage $P_{0,1}$ of point
                       $[(0,1)]\in \OC_{0,2}\subset \OC_{2,0}$.
                       Analysis analogous to the one from the previous 
                       subsection shows that $P_{0,1}$ does not intersect
                        any boundary stratum of $\OC_{n,m}$. Thus, as
                    a chain of codimension $2$ this preimage $P_{0,1}$
                     coincides with the union of closures
                      of strata $C_T$ of codimension
                     $2$ such that $C_T\subseteq P_{0,1}$. 
                     It is easy to see that any such stratum give pictures like
                     the one below where there is no arrow going from
  circled regions outside (as in the picture in 6.4.2.2),
   \vskip 1cm      \centerline{
\epsfbox{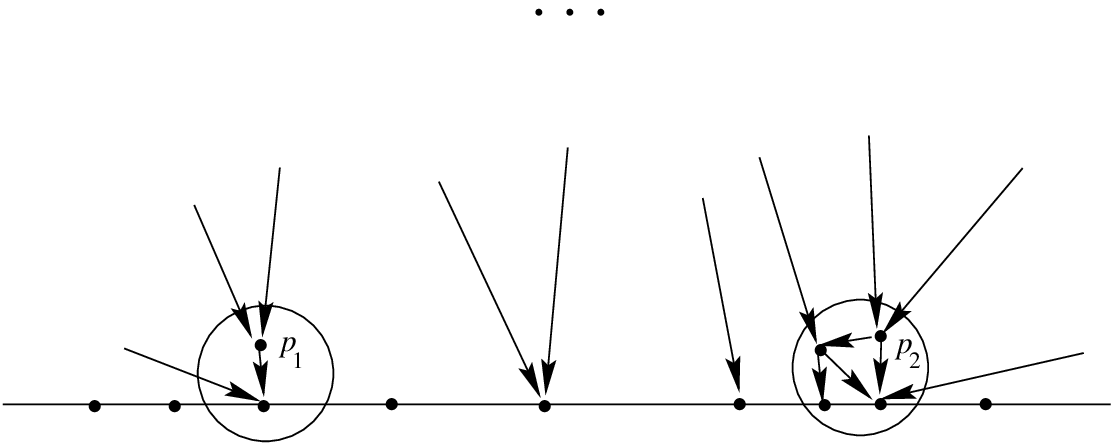}
}  \vskip 11mm \hfill\break
 and we get exactly the cup-product in the tangent cohomology
    of
    the Hochschild complex 
   as was described
 above.

\par\noindent {\bf 8.2.3. Homotopy between two pictures}
\vskip 0.5truecm                        
                       Choosing a path
                       form one (limiting) configuration of two
                       points
                        on $\H$ to another configuration,
                         \vskip 1cm     \centerline{
\epsfbox{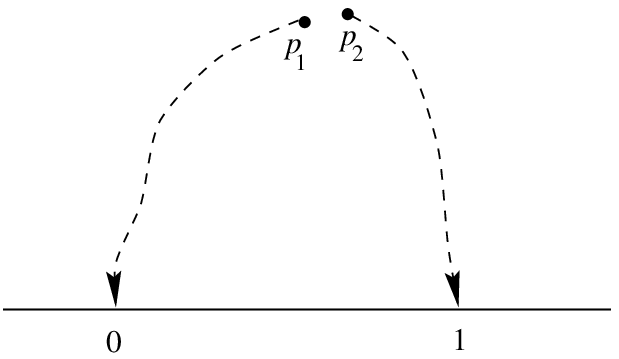}
} \vskip 11mm   \hfill\break   
                          we see that two products coinside
                         on the level of cohomology. \qed

\vskip 0.5truecm            
        
\par\noindent {\bf 8.3. First application: Duflo-Kirillov isomorphism}
\vskip 0.5truecm

\par\noindent {\bf 8.3.1. Quanization of the Kirillov-Poisson bracket}
\vskip 0.5truecm 
            
            Let $\g$ be a {\it finite-dimensional} Lie algebra over $\R$.
            The dual space to $\g$ endowed with the Kirillov-Poisson
             bracket  is naturally a Poisson manifold (see [Ki]). 
             We remind here the formula for
              this bracket: if $p\in \g^*$ is 
              a point and $f,g$ are two functions
             on $\g$ then the value $\{f,g\}_{|p}$ is defined as
             $\langle p,[df_{|p},dg_{|p}]\rangle$ where differentials
              of functions $f,g$ at $p$ are considered as elements of $
              \g\simeq (\g^*)^*$.
               One can consider $\g^*$ as an algebraic Poisson manifold
               because coefficients of the Kirillov-Poisson bracket
                are linear functions on $\g^*$.
                
                \proclaim Theorem.  The canonical quantization
                 of the Poisson manifold $\g^*$ is isomorphic to
                 the  family of algebras $\U_{\hbar}(\g)$ defined as 
                 universal enveloping algebras of $\g$ endowed with the
                  bracket $\hbar[\,,\,]$. 
                  \par

                  {\bf Proof:} 
               in  6.4 we constructed a canonical star-product
                 on the algebra of functions on a affine space
                  for a given Poisson structure.
                   Thus, we have canonical star-product on $C^{\infty}(\g^*)$.
                   We claim that the product of any two {\it polynomials}
                   on $\g^*$ is a polynomial in $\h$ with coefficients
                    which are polynomials on $\g^*$.
                    The reason is that the star-product is constructed
                     using contraction of indices. Let us denote by
                      $\beta \in \g^*\otimes\g^*
                      \otimes \g$ the tensor giving the Lie bracket on $\g$.
                      All natural operations
                       $Sym^k(\g)\otimes Sym^l(\g)\ra Sym^m(\g)$ which 
                       can be defined
                        by contractions of indices with several copies of
                         $\beta $, exist only for $m\le k+l$, and for every 
                         given $m$
                          there are only finitely many ways to contract 
                          indices.
                          Thus, it makes sense to put $\h$ equal to $1$ 
                          and obtain a product on $Sym(\g)=\oplus_{k\ge 0}
                          Sym^k(\g)$.
                           We denote this product also by $\star$.
                     
                      It is easy to see that for $\ga_1,\ga_2\in\g$ the 
                         following identity holds:
                         $$\ga_1\star\ga_2-\ga_2\star\ga_1=
                         [\ga_1,\ga_2]\,\,\,.$$
                         Moreover, the top component of $\star$-product
                          which maps $Sym^k (\g)\otimes Sym^l(\g)$
                           to $Sym^{k+l}(\g)$,
                           coincides with the product on $Sym(\g)$.
                            From this two facts one concludes that
                           there exists a unique isomorphism of algebras
                           $$I_{alg}:(\U\g,\cdot)\ra (Sym(\g),\star)$$
                            such that $I_{alg}(\ga)=\ga$ for $\ga\in \g$,
                             where $\cdot$ denotes the universal enveloping
                              algebra of $\g$ with the standard product.

                               One can easily recover variable $\hbar$ in this
                              description and get the statement of the theorem.
                              \qed
                              
                              \proclaim Corollary.
                               The center of the universal enveloping algebra
                                is canonically isomorphic as an algebra
                                 to the algebra 
                  $\bigl(Sym(\g)\bigr)^{\g}$ 
                                of $\g$-invariant polynomials
                                 on $\g^*$. \par
         {\bf Proof:} The center of $\U \g$    is $0$-th cohomology for 
         the (local) Hochschild 
                          complex of $\U\g$ endowed with the 
                          standard cup-product. The algebra 
                       $\bigl(Sym(\g)\bigr)^{\g}$  is the 
                  $0$-th cohomology of
                            the algebra of polyvector fields on $\g^*$
                             endowed with the differential
                              $[\a,\cdot]$ where $\a$ is the Kirillov-Poisson
                              bracket.  From the theorem
                               8.2 we conclude that applying the tangent
                               map  to $\U$
                                we get an isomorphism of algebras.

                  \vskip 0.5truecm
\par\noindent {\bf 8.3.2. Three isomorphisms}
\vskip 0.5truecm

                  	In the proof of theorem 8.3.1 we introduced
                  	 an isomorphism $I_{alg}$ of algebras.
                        
                              We denote by $I_{PBW}$ the isomorphism of
                              vector spaces
                              $$Sym(\g)\ra \U
                              \g$$ 
                              (subscript from the 
                              Poincar\'e-Birkhoff-Witt theorem, see 8.3.5.1),
                              which is defined as
                              $$\ga_1\ga_2\dots\ga_n\ra 
                              {1\over n!}\sum_{\sigma\in\Sigma_n}
                              \ga_
                              {\sigma_1}
                              \cdot\ga_{\sigma_2}
                              \cdot\dots\cdot \ga_{\sigma_n}\,\,\,.$$
                              
                            Analogously to arguments from above, one can see
                             that the tangent map from polyvector fields
                          on $\g^*$ to the Hochschild complex
                           of the quantized algebra can be defined for 
                           $\h=1$ and
                           for polynomial coefficients. We denote by $I_T$ its
                           component which maps polynomial
                            $0$-vector fields on $\g^*$
                            (i.e. elements of $Sym(\g)$) to $0$-cochains 
                            of the 
                            Hochschild complex of the algebra $(Sym(
                            \g),\star)$.
                            Thus, $I_T$ is an isomorphism of vector spaces
                            $$I_T:Sym(\g)\ra Sym(\g)$$
                         and the restriction of $I_T$ to
                         the algebra of $ad(\g)^*$-invariant polynomials 
                         on $\g^*$
                          is an isomorphism of algebras
                          $$Sym(\g)^\g\ra Center((Sym(\g),\star))\,\,\,.$$
                          
                          Combining all facts from above we get a sequence of 
                          isomorphisms of vector spaces:
                        $$Sym(\g){\buildrel I_T\over  \raa }
                         Sym(\g )
                         {\buildrel I_{alg} \over \laa}\U\g
                        {\buildrel I_{PBW} \over\laa }
                        Sym(\g)\,\,\,.$$
                        
                        These isomorphisms are $ad(\g)$-invariant.
                         Thus, one get isomorphisms
                        $$ (Sym(\g))^\g\buildrel {I_T}_{|\dots}
                        \over \raa 
                        Center(Sym(\g),\star)\buildrel
                        {I_{alg}}_{|\dots}\over \laa
                        Center(\U\g)\buildrel {I_{PBW}}_{|\dots}\over
                        \laa (Sym(\g))^\g\,,$$
                        where the subscript $|\dots$ denotes the restriction
                         to subspaces of $ad(\g)$-invariants .
                        Moreover, first two arrows are isomorphism of algebras.
                         Thus, we proved the following
                        
                        \proclaim Theorem. The restriction of the map 
                        $$\bigl(I_{alg}\bigr)^{-1}\circ I_T:Sym(\g)\ra \U\g$$
                         to $(Sym(\g))^\g$ is an isomorphism 
                          of algebras $(Sym(\g))^\g\ra \, Center\,(\U\g)$.
                          \par
                        
                         \qed
                        
          \vskip 0.5truecm
\par\noindent {\bf 8.3.3. Automorphisms of $Sym(\g)$}
\vskip 0.5truecm                
                    
                        Let us calculate automorphisms 
                        $I_T$ and $I_{alg}\circ I_{PBW}$
                        of the vector space 
                        $Sym(\g)$. We claim that
                        both these automorphisms are translation invariant
                        operators on the space $Sym(\g)$
                        of polynomials on $\g^*$.
                        
                        The algebra of translation invariant
                operators on the space of
                polynomials on a vector space  $V$ is canonically
                isomorphic to the algebra of formal power series generated
                 by  $V$. Generators of this algebra
                  acts as derivations along constant vector fields in $V$.
                  Thus, any such operator can be seen as a formal power
                   series at zero on the dual vector space $V^*$. We apply 
                   this formalism to the case $V=\g^*$.
                   
                   \proclaim Theorem. Operators 
                   $I_T$ and $I_{alg}\circ I_{PBW}$ respectively
                    are translation 
                   invariant
                operators associated with formal power series 
                $S_1(\ga)$ and $S_2(\ga)$
                at zero in
                $\g$ of the form
                $$S_1(\ga)=exp\left(\sum_{k\ge 1}
                 c_{2k}^{(1)}\,Trace\,(ad( \ga)^{2k})\right),\,\,
             S_2(\ga)=exp\left(\sum_{k\ge 1} 
             c_{2k}^{(2)}\,Trace\,(ad( \ga)^{2k})\right)  $$
           where $c_2^{(1)},c_4^{(1)},\dots$ and  $c_2^{(2)},c_4^{(2)},\dots$ 
           are two infinite sequences of real numbers indexed by even natural
            numbers.
            \par
            
            {\bf Proof:} we will study separately two cases.

          \vskip 0.5truecm
\par\noindent {\bf 8.3.3.1. Isomorphism $I_T$}
\vskip 0.5truecm  

            The isomorphism $I_T$ is given by the sum over terms
             corresponding to admissible graphs $\G$ 
             with no vertices of the second type,
             one special vertex $v$ of the first type such that no edge start 
             at $v$,
             and such that at any other vertex  start two 
             edges and ends no more than one edge. Vertex $v$ is 
             the marked vertex where we put an element of $Sym(\g)$
             considered as an element of tangent cohomology.
              At other vertices we put
                    the Poisson-Kirillov bi-vector field on $\g^*$, i.e. the
                     tensor of commutator operation in $\g$.
                      As the result we get $0$-differential operator,
                       i.e. an element of algebra $Sym(\g)$.

   It is easy to see that any such graph is isomorphic to a union of 
    copies of ``wheels'' $Wh_n,\,\,\,n\ge 2$:
    
   \vskip 1cm  \centerline{
\epsfbox{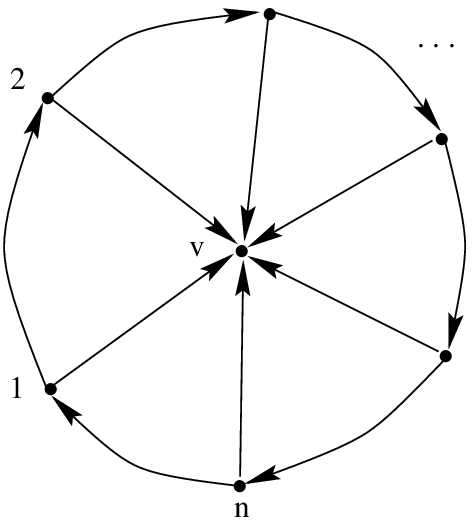}
} \vskip 11mm\hfill\break
     with identified central vertex $v$. The following picture shows a typical
     graph:
     
\vskip 1cm
\centerline{ \epsfbox{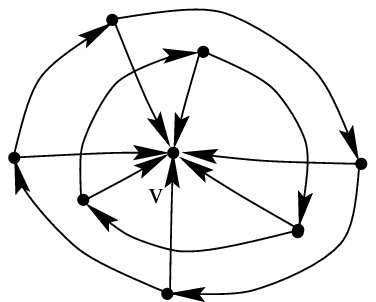}  }  
\vskip 11mm

             In the integration we may assume that the point corresponding 
             to $v$
              is fixed, say that it is $i\cdot 1+0\in\H$,
               because group $G^{(1)}$ acts simply transitively on $\H$. 
               First of all, the 
              operator
              $Sym(\g)\ra Sym(\g)$
               corresponding to the individual wheel $Wh_n$ is 
              the differential operator on $\g^*$ with constant coefficients,
               and it corresponds to the polynomial $\ga\mapsto 
               \,Trace\,(ad(\ga)^n)$
                on $\g$. The operator corresponding to the joint of several 
                wheels
                 is the product of operators associated with individual wheels.
                  Also, the integral corresponding to the joint is the 
                  product of
                  integrals. Thus, with the help of  symmetry factors, we get
                  that
                  the total operator is equal to the exponent of the sum of 
                  operators associated with
                   wheels $Wh_n,\,\,n\ge 2$ with weights
                  equal to  corresponding integrals. By the symmetry 
                  argument
                  used several times before ($z\mapsto -\overline{z}$), we see
                   that integrals corresponding to wheels with odd $n$ vanish.
                   We proved the first statement of our theorem. \qed

          \vskip 0.5truecm
\par\noindent {\bf 8.3.3.2. Isomorphism $I_{alg}\circ I_{PBW}$}
\vskip 0.5truecm  

           The second case, for the operator $I_{alg}\circ I_{PBW}$, 
                   is a bit more tricky.
                   Let us write a formula for this map:
                   $$I_{alg}\circ I_{PBW}:
                   \ga^n\mapsto \ga\star\ga\star\ga\dots\star\ga\,\,\,\,\,
                   (n\,\,{\rm\, copies\,\,\,of\,\,\,}\ga)\,\,\,.$$
                   This formula defines the map unambiguously because
                    elements $\ga^n,\,\,\ga\in\g,\,\,\,n\ge 0$ generate
                     $Sym(\g)$ as a vector space.
                     
                   In order to multiply several (say, $m$, where $m\ge 2$)
                    elements of the quantized algebra we should put these
                     elements at $m$ {\it fixed} points in increasing order
                     on $\R$ and take the sum
                      over all possible graphs with $m$ vertices of the
                       second type of corresponding expressions with 
                       appropriate 
                       weights. The result does not depend on the position
                        of fixed points on $\R$ because the star-product
                        is associative. Moreover, if we calculate a power
                         of a given element with respect to the 
                         $\star$-product,
                     we can put all these points in arbitrary order.
                     It follows that we can take an average over 
                     configurations
                      of $m$ points on $\R$ where each point is random, 
                      distributed independently from other points, with 
                      certain
                      probability density on $\R$. We choose a probability 
                      distribution on $\R$ with a smooth symmetric (under
                       transformation $x\mapsto -x$) density $\rho(x)$.
                        We assume also that $\rho(x)dx$ is the restriction
                         to $\R\simeq C_{1,1}$ of a smooth $1$-form on
                         $\OC_{1,1}\simeq \{-\infty\}\sqcup \R\sqcup 
                         \{+\infty\}$.
                       With probability $1$ our $m$ points will be 
                       pairwise
                       distinct. One can check easily that the 
                       interchanging
                      of order of integration (i.e. for the taking 
                      mean value 
                      from the probability theory side, and for the 
                      integration
                       of differential forms over configuration
                        spaces) is valid operation in our case.
                    
                    The conclusion is that the $m$-th  power of an element 
                    of quantized algebra can be 
                    calculated as a sum over all graphs with $m$ vertices of 
                    the second type, with weights equal to integrals
                     over configuration spaces where we integrate 
                     products of forms $d\phi$ and $1$-forms $\rho(x_i)dx_i$
                      where $x_i$ are points moving along $\R$.
                           
                           The basic element of pictures in our case
                            are ``wheels without axles'':
                           
 \vskip 1cm \centerline{ \epsfbox{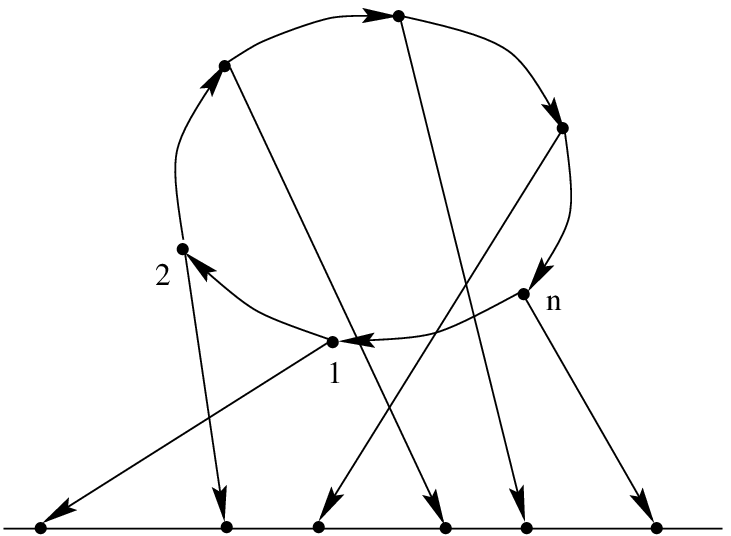}  }      \vskip 11mm

 and the $\Lambda$-graph (which gives $0$ by symmetry reasons):
 
    \vskip 1cm \centerline{ \epsfbox{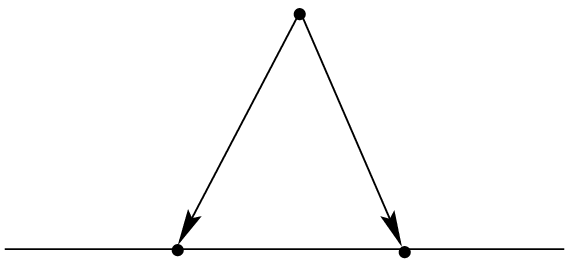}  }      \vskip 11mm
                  
               The typical total picture is something like (with $m=10$):
               
\vskip 1cm   \centerline{ \epsfbox{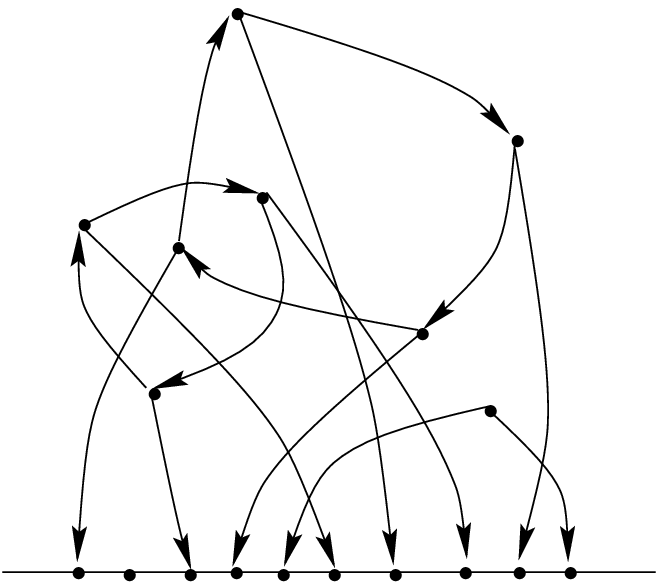}  }    \vskip 11mm             
                        
                  Again, it is clear from all this that the operator
                $I_{alg}\circ I_{PBW}$ is a differential operator
                 with constant coefficients on $Sym(\g)$, equal to
                  the exponent of the sum of operators corresponding to
                   individual wheels. These operators are again proportional
                    to operators associated with power series on $\g$
                     $$\ga \ra \,Trace\,(ad(\ga)^n)\,\,\,.$$    
                      By the same symmetry
                      reasons as above we see that integrals
                         corresponding to odd $n$ vanish.
                     The second part of the theorem is proven. \qed

          \vskip 0.5truecm
\par\noindent {\bf 8.3.4. Comparison with the Duflo-Kirillov isomorphism}
\vskip 0.5truecm

                           For the case of semi-simple $\g$ 
                           there is so called Harish-Chandra isomorphism
                            between algebras $\bigl(Sym(\g)\bigr)^{\g}$
                             and $Center(\U\g)$.  
                             A.~Kirillov realized that there is a way to
                             rewrite the Harish-Chandra isomorphism
                              in a form which has sense for arbitrary
                               finite-dimensional Lie algebra, i.e.
                                without 
                               using  the Cartan and
                                Borel subalgebras,  the Weyl group etc.
                                 Later  M.~Duflo (see [D]) proved that
                                 the map proposed by Kirillov is an isomorphism
                                 for all finite-dimensional Lie algebras.
                                 
                          The explicit formula for the Duflo-Kirillov
                           isomorphism
                           is the following: 
                            
                          $$I_{DK}
                          :\bigl(Sym(\g)\bigr)^{\g}
                         \simeq \,Center( \U(\g)),
                          \,\,\,\, I_{DK}= {I_{PBW}}_
                          {|(Sym(\g))^{\g}}  
                          \circ {I_{strange}}_{
                        |(Sym(\g))^{\g}  },
                          $$
                           where 
                          $I_{strange}$ 
                              is an invertible translation invariant
                              operator on
                             $Sym(\g)$ 
              associated with the following
               formal power series on $\g$ at zero,
               reminiscent of the 
                      square root of the Todd class:
                     $$\ga\mapsto \,{\rm exp}\,
                       \left(\sum_{k\ge 1} {B_{2k}\over 4 k (2k)!}\, Trace \,
                     (ad(\ga)^{2k})\right) $$
                     where $B_2,B_4,\dots$ are Bernoulli numbers.
                      Formally, one can write
                       the r.h.s. as
                       $det(q(ad(\ga)))$ where
                       $$q(x):=\sqrt{ {e^{x/2}-e^{-x/2}\over x}}\,\,\,.$$

                     The fact that the Duflo-Kirillov
                      isomorphism is an isomorphism
                      of algebras is highly non-trivial. All proofs known
                       before (see [Du],[Gi]) used 
                        certain facts about finite-dimensional Lie algebras
                        which follow only from the classification theory.
                        In particular, the fact that the analogous isomorphism
                         for Lie superalgebras
                        is compatible with products, was not known.
                         
                         We claim that our isomorphism 
                          coincides with the 
                         Duflo-Kirillov 
                          isomorphism.
                      Let us just sketch the argument. 
                          In fact, we claim that
                          $$I_{alg}^{-1}\circ I_T=I_{PBW}\circ I_{strange}
                          \,\,\,.$$
                       If it is not true then we get a non-zero
                        series $Err\in t^2\R[[t^2]]$ such that
                         the translation invariant operator on $Sym(\g)$ 
                         associated with
                          $\ga\mapsto I_{det(exp(ad\ga)))}$ gives an
                          {\it automorphism} of algebra $(Sym(\g))^\g$.
                           Let $2k>0$ be the degree of first non-vanishing
                            term in the expansion of $Err$.
                        Then it is easy to see that
                         the operator on $Sym(\g)$ 
                         associated with the polynomial $\ga\mapsto
                          Trace(ad(\ga)^
                         {2k}$ is a {\it derivation} when restricted
                          to $(Sym(\g))^\g$.
                       One can show that it is not true
                            using Lie algebras  $\g=gl(n)$ for large $n$.
                            Thus, we get a contradiction and proved
                             that $Err=0$.
                           \qed

                            As a remark we would like to mention that
                             if one replaces series $q(x)$
                              above just by
                         $$\left({x\over 1-e^{-x}}\right)^{-{1\over 2}}
                         $$ 
                         then one still get an isomorphism of algebras.
                          The reason is that the one-parameter group
                          of automorphisms of $Sym(\g)$ associated with series
                          $$\ga\ra {\rm exp} (const \cdot Trace(ad(\ga)))$$
                           {\it preserves} the  structure 
                            of Poisson algebra on $\g^*$.
                            This  one-parameter group also acts by 
                            automorphisms of $\U\g$. It is analogous to 
                            the  Tomita-Takesaki flow of weights for
                          von Neumann  factors.

                                \vskip 0.5truecm
\par\noindent {\bf 8.3.5. Results in rigid tensor
 categories}
\vskip 0.5truecm 

 Many proofs from this paper  can be transported to a more
  general context of rigid $\Q$-linear tensor categories (i.e. abelian
   symmetric monoidal categories with the duality functor imitating 
   the behavior of finite-dimensional vector spaces).
We will be very brief here. 

 First of all, one can formulate and prove
 the Poincar\'e-Birkhoff-Witt theorem in a great generality,
  in $\bf Q$-linear additive
  symmetric monoidal categories with infinite sums and kernels of projectors.
 For example, it holds in the category of $A$-modules where $A$
 is  arbitrary commutative associative algebra over $\Q$.
  Thus, we can speak about universal enveloping algebras and the isomorphism
   $I_{PBW}$. 

One can define  Duflo-Kirillov morphism for a Lie algebra in
a  $\k$-linear rigid tensor
 category where $\k$ is a field of characteristic zero, because
  Bernoulli numbers are rational. Our result from 8.3.4 saying
   that it is a morphism of algebras, 
   holds in this generality as well. It {\it does not} hold for 
   infinite-dimensional Lie algebras because we use traces of products
    of operators
    in the adjoint representation.
    
  In [KV] a conjecture
   was made in the attempt to prove that that Duflo-Kirillov formulas
    give a morphism of algebras. It seems that using our
      result one can  prove this conjecture. Also, there is 
       another related conjecture concerning two
        products in the algebra of chord diagrams (see [BGRT])
       which seems to follow from our results too.

              \vskip 0.5truecm
\par\noindent {\bf 8.4. Second application: algebras of $Ext$-s.}
\vskip 0.5truecm                         
                
                   Let $X$ be complex manifold, or a smooth algebraic variety
                    of field $\k$ of characteristic zero.
                    We associate with it two graded vector spaces.
                     The first space $HT^\b(X)$ is the direct sum $
                     \bigoplus_{k,l}
                      H^k(X,\wedge^l T_X)[-k-l]$.
                       The second space $HH^\b(X)$
                       is the space $\bigoplus_k
                        Ext^k_{Coh(X\times X)} ({\cal O}_{diag}, 
                        {\cal O}_{diag})[-k]$
                        of $Ext$-groups in the category
                         of coherent sheaves on $X\times X$ from the sheaf of 
                         functions on the diagonal to itself.
                          The space $HH^\b(X)$ can be thought as the
                          Hochschild cohomology of the space $X$.
                           The reason is that the Hochschild
                           cohomology of any algebra $A$ can be also defined
                            as $Ext^\b_{A-mod-A}(A,A)$ in the category
                             of bimodules.

                          Both spaces, $HH^\b(X)$ and $HT^\b(X)$
                           carry natural products.
                            For $HH^\b(X)$ it is the Yoneda composition, and
                             for $HT^\b(X)$ it is the cup-product
                             of cohomology and of polyvector fields.
                             
                             \proclaim Claim. 
                             Graded algebras $HH^\b(X)$ and $HT^\b(X)$
                              are canonically isomorphic. The isomorphism
                               between them is functorial with respect to 
                               \'etale maps. \par
                              
                              This statement is again a corollary of
                              the theorem from 8.2. 
                                We will give the proof of it, and explain 
                              an application 
                               to the Mirror Symmetry (see [Ko4])  
                                in the next paper.
   
   \vskip 0.5truecm
   
   \centerline{\bf Bibiliography}\par
   
   \vskip 0.5truecm
   
   \item{[AKSZ]}
   M.~Alexandrov, M.~Kontsevich, A.~Schwarz, O.~Zaboronsky,
 {\sl The Geometry of the Master Equation and 
 Topological  Quantum Field Theory}, Intern.\ Jour.\ of Mod.\ Phys.,
  {\bf 12} (1997), no. 7, 1405 - 1429, and hep-th/9502010.
  
  \item{[AGV]} V.~I.~Arnold, S.~M.~Gusein-Zade, A.~N.~Varchenko,
   {\sl Singularities of Differentiable Maps}, Vol. I,
    Birkh\"auser 1985.
    
    \item{[BGRT]}, D.~Bar-Natan, S.~Garoufalidis, L.~Rozansky,
     D.~Thurston, {\sl Wheels, wheeling, and the Kontsevich integral of the
      unknot}, q-alg/9703025.

   \item{[BFFLS]} F.~Bayen, M.~Flato, C.~Fr\o nsdal, A.~Lichnerowicz, 
   D.~Sternheimer,
 {\sl Deformation theory and quantization. I. Deformations of symplectic
  structures}, Ann.\ Physics {\bf 111} (1978), no. 1, 61 - 110.
   
   \item{[De]} P.~Deligne, {\sl Cat\'egories tannakiennes}, The Grothendieck 
   Festschrift,
    Vol. II, Progress in Mathematics 87, Birkh\"auser 1990, 111 - 195.
    
    \item{[DL]} M.~De Wilde, P.~B.~A.~Lecomte,  {\sl Existence of star-products
 and of formal deformations in Poisson Lie algebra of arbitrary symplectic
  manifolds}, Let.\ math.\ Phys., {\bf 7} (1983), 487 - 496.

\item{[Du]} M.~Duflo, {\sl Caract\'eres des alg\`ebres de Lie r\'esolubles},
 C.\ R.\ .Acad.\ Sci., {\bf 269} (1969), s\'erie a, p. 437 - 438.
 
 \item{[EK]} P.~Etingof, D.~Kazhdan, {\sl Quantization of Lie Bialgebras, I},
 Selecta Math., New Series, {\bf 2} (1996), no. 1, 1 - 41.

   \item{[Fe]} B.~Fedosov, {\sl A simple geometric construction of deformation 
 quantization}, J.\ Diff.\ Geom., {\bf  40} (1994) 2, 213 - 238.
  
  \item{ [FM] } W.~Fulton, R.~MacPherson, {\sl Compactification of configuration
    spaces}, Ann.\ Math. {\bf 139} (1994), 183 - 225.
   
   \item {[GK]} I.~M.~Gelfand, D.~A.~Kazhdan, {\sl Some problems of differential
    geometry and the calculation of cohomologies of Lie algebras of vector
     fields}, Soviet Math.\ Dokl.,  {\bf 12} (1971), no. 5, 1367 - 1370.

  \item{ [GV] } M.~Gerstenhaber, A.~Voronov, {\sl
  Homotopy $G$-algebras and moduli space
    operad}, Intern.\ Math.\ Res.\ Notices (1995), No. 3, 141 - 153.
   
  \item{ [GJ] }E.~Getzler, J.~D.~S.~Jones, {\sl
  Operads, homotopy algebra and iterated
         integrals for double loop spaces},  hep-th/9403055.
    
    \item{[Gi]} V.~Ginzburg, {\sl Method of Orbits in the Representation
     Theory of Complex Lie  Groups},  Funct.\ Anal.\ and appl.,
     {\bf 15} (1981), no. 1, 18 - 28.

\item{[GM]} W.~Goldman, J.~Millson, {\sl The homotopy invariance
 of the Kuranishi space}, Ill.\ J.\ Math., {\bf 34} (1990), no. 2, 337 - 367.
  
\item{[HS1]} V.~Hinich, V.~Schechtman, {\sl Deformation theory and Lie algebra
 homology}, alg-geom/9405013.
 
 \item {[HS2]} V.~Hinich, V.~Schechtman, {\sl Homotopy Lie algebras},
  I.~M.~Gelfand Seminar, Adv.\ Sov.\ Math., {\bf 16} (1993), part 2, 1 - 28.

    \item {[KV]} M.~Kashiwara, M.~Vergne, {\sl The Campbell-Hausdorff
     formula and invariant hyperfunctions}, Invent.\ Math., {\bf 47}
      (1978), 249 - 272.
    
    \item{[Ki]} A.~Kirillov, {\sl Elements of the Theory of Representations},
    Springer-Verlag 1975.

   \item{[Ko1]} M.~Kontsevich, {\sl Feynman diagrams and low-dimensional topology},
    First European Congress of Mathematics (Paris, 1992), Vol. II,
     Progress in Mathematics 120, Birkh\"auser (1994), 97 - 121.
     
  \item { [Ko2]} M.~Kontsevich, {\sl
  Formality Conjecture}, D.~Sternheimer et al. (eds.),
    Deformation Theory and Symplectic Geometry, Kluwer 1997, 139 - 156.  
   
  \item{ [Ko3] } M.~Kontsevich, {\sl
  Rozansky-Witten invariants via formal geometry},
    dg-ga/9704009, to appear in Compos.\ Math.
    
    \item{[Ko4]} M.~Kontsevich, {\sl Homological algebra of mirror symmetry},
 Proceedings of ICM, Z\"urich 1994, vol.\ I, Birkh\"auser (1995), 120 - 139.

   \item{[M] }Yu.~I.~Manin, {\sl Gauge Field theory and Complex Geometry},
    Springer-Verlag 1988.
   
  \item{ [Q] } D.~Quillen, {\sl Superconnections and the
   Chern character}, Topology {\bf 24} (1985), 89 - 95.

\item{[SS]} M.~Schlessinger, J.~Stasheff, {\sl The Lie algebra structure 
on tangent cohomology and deformation theory}, J.\ Pure Appl.\ Algebra {\bf 89}
 (1993), 231 - 235.

\item{[Su]} D.~Sullivan, {\sl Infinitesimal computations in topology}, 
I.\ H.\ E.\ S.\ Publ.\ Math.,
  no. 47 (1977), 269 - 331.

 \item{[V]} A.~Voronov, {\sl Quantizing Poisson Manifolds}, q-alg/9701017.

\vskip 2cm

I.H.E.S., 35 Route de Chartres,

Bures-sur-Yvette 91440, FRANCE

\medskip

email: {\it maxim@ihes.fr}

          \end{document}